\documentclass[11pt,reqno]{amsart}
\usepackage[top=1in, bottom=0.9in, left=1in, right=1in]{geometry} 
\usepackage{amsmath, amssymb, graphicx, mathrsfs, mathtools}
\usepackage{amsthm}
\usepackage{palatino,bbm}
\usepackage[makeroom]{cancel}
\usepackage[usenames,dvipsnames,svgnames]{xcolor}
\usepackage{bm}%for boldmath
\usepackage{latexsym,comment}
\usepackage{upref, eucal}
\usepackage[normalem]{ulem} %%for a strikethrough font 
\usepackage[all]{xy}
\usepackage{xcolor}
\usepackage{palatino, mathrsfs}
\usepackage[makeroom]{cancel}
\usepackage{url}

\usepackage[hypertexnames=false]{hyperref}

\usepackage{tikz-cd}

\usepackage{autonum}
\usepackage{enumitem}%%for controlling spaces in enumerate
\def\namedlabel#1#2{\begingroup
	#2%
	\def\@currentlabel{#2}%
	\phantomsection\label{#1}\endgroup
}
\makeatother
\usepackage{tikz,color,soul}
\usetikzlibrary{matrix,arrows}

\usetikzlibrary{calc}

\usepackage{tikz-cd}

\allowdisplaybreaks

\newcommand {\nc} {\newcommand}
\newcommand {\enm} {\ensuremath}

\def \d{\delta}
\def \uuu{u}
\nc {\bdm} {\begin{displaymath}}
	\nc {\edm} {\end{displaymath}}

\numberwithin{equation}{section}
\definecolor{labelkey}{rgb}{1,0,0}

\newtheorem{theorem}{Theorem}[section]

\theoremstyle{plain} % just in case the style had changed
\theoremstyle{plain} % just in case the style had changed
\newcommand{\thistheoremname}{}
\newtheorem{genericthm}[theorem]{\thistheforemname}

\newtheorem*{genericthm*}{\thistheoremname}
\newenvironment{namedthm*}[1]
{\renewcommand{\thistheoremname}{#1}%
	\begin{genericthm*}}
	{\end{genericthm*}}

\newtheorem{lemma}[theorem]{Lemma}
\newtheorem{proposition}[theorem]{Proposition}
\newtheorem{corollary}[theorem]{Corollary}
\newtheorem{remark}{Remark}
\newtheorem*{definition*}{Definition}

\newtheorem*{conjecture*}{Conjecture}
\newtheorem*{theorem*}{Theorem}
\newtheorem*{remark*}{Remark}
\newtheorem{example}[theorem]{Example}
\numberwithin{equation}{section}
\DeclareMathOperator{\Gal}{Gal}

\makeatletter
\renewcommand{\pmod}[1]{\mkern4mu({\operator@font mod}\mkern 6mu#1)}
\makeatother

\newcommand\fp{\mathfrak{p}} 
\newcommand\fm{\mathfrak{m}}
\newcommand\fa{\mathfrak{a}}
\newcommand\ft{\mathfrak{t}}
\newcommand\fg{\mathfrak{g}}
\newcommand\fh{\mathfrak{h}}
\newcommand\fl{\mathfrak{l}}
\newcommand\fc{\mathfrak{c}}
\newcommand\fd{\mathfrak{d}}

\newcommand\fs{\mathfrak{s}}
\newcommand\fe{\mathfrak{e}}

\newcommand\Carl{\psi^{\text{Carl}}}

\nc {\form}[1] {\enm{\mbox{\underline{for}}}_{#1}}
\nc {\prol}[1] {\enm{\mbox{\underline{prol}}_{{#1}^*}}}

\nc {\stk} {\stackrel}

\nc{\beqar}{\begin{eqnarray*}}
	\nc{\eeqar}{\end{eqnarray*}}            
\newcommand{\GL}{\mathrm{GL}}
\newcommand{\SL}{\mathrm{SL}}
\newcommand{\Pn}[2] {\ensuremath{ {\mathbb{P}}^{#1}_{#2}}}
\nc{\Quot}[3]{\enm{ {\mathfrak{Quot}_{ {#1}/{#2}/{#3}}}}}
\nc{\Hilb}[2]{\enm{ {\mathfrak{Hilb}_{ {#1}/{#2}}}}}

\newcommand{\mf}[1]{\mathfrak{#1}}
\newcommand{\bb}[1]{\mathbb{#1}}

\newcommand{\lcm}{\text{lcm}}

\nc {\Coh}[4] {\ensuremath{H^{#1}(\Pn{#2}{},{#3}({#4}))}}
\nc {\Ch}[3] {\enm{H^{#1}(X_t,{#2}_t({#3}))}}
\nc {\Qphi}[4]{\enm{ {\mathfrak{Quot}^{~#4}_{ {#1}/{#2}/{#3}}}}}
\nc {\Gra}[4]{\enm{ {\mathfrak{Grass}_{#2}({#3},{#4})}}}
\nc {\HomA}[2]{\enm{\mathrm{Hom}_A{#1}{#2}}}
\nc {\tr}{\mathrm{tr}}

\nc {\C}[2]{\binom{#1}{#2}}
\def \mb{\mbox}
\def \Z{{\mathbb Z}} 

 \def \Z{{\mathbb Z}}

   \def \h{\hat{\ }}

\def \d{\delta}

  \def \bX{{\bf X}} \def \bH{{\bf H}}

\def \Ga{\mathbb{G}_{\mathrm{a}}}

\newcommand{\Q}{{\mathbb Q}}

\newcommand{\F}{{\mathbb{F}}}
\newcommand{\fb}{{\mathfrak{b}}}
\def \R1{R((q))[q']\h}

\usepackage{xcolor}

\DeclareMathOperator{\Spec}{\mathrm{Spec}}

\newcommand{\End}{\mathrm{End}}

\newcommand{\oldmarginpar}[1]{}

\newcommand{\main}{\sideset{}{_{\text{main}}}\sum}
\newcommand{\error}{\sideset{}{_{\text{error}}}\sum}
\newcommand{\summonics}{\sideset{}{'}\sum}

\newcommand{\ff}{\mathfrak{f}}

\nc{\bx}{{\mathbf x}}
\nc{\by}{\bold{y}}
\nc{\bz}{\bold{z}}
\nc{\ba}{\bold{a}}
\nc{\Fp}{\tilde{F}}
\nc{\Rp}{\tilde{R}}
\nc{\mlow}{m_{\mathrm{l}}}
\nc{\mup}{m_{\mathrm{u}}}
\nc{\ord}{\mb{ord }}
\nc{\bXp}{\bX_{\mathrm{prim}}}
\nc{\bPsi}{\bold{\Psi}}
\nc{\mult}{\mathrm{mult}}
\nc{\mbB}{\mathbbm{B}}
\nc{\mfor}[1]{{#1}^{\mathrm{for}}}
\nc{\Hdr}{\bH^1_{\mathrm{dR}}(A)}

\nc{\ale}{\textcolor{red}}
\nc{\emma}{\textcolor{orange}}
\nc{\bluev}{\textcolor{blue}}
\nc{\mP}{\mathfrak{p}}

\nc{\prolarrow}[1]{\stk{(\uuu_{#1},\d_{#1})}{\longrightarrow}}
\nc{\tc}{\textcolor}
\nc{\nexp}[1]{\exp_{\d,{#1}}}
\nc{\gr}{\mathrm{Gr}}

\nc{\tJ}[1]{{{\mathcal{J}}^n{#1}}} %Buium algebraic functor

\nc{\Sn}{S^{(n)}}
\nc{\Snr}[1]{S^{(n)}_{#1}}
\nc{\Sp}[1]{{\mathbf{Sp}_{#1}}}
\nc{\Aff}[1]{{\mathbf{Aff}_{#1}}}
\nc{\Sch}[1]{{\mathbf{Sch}_{#1}}}
\nc{\ov}[1]{{\overline{#1}}}

\nc{\Qbar}{\overline{\bb{Q}}}
\nc{\Aut}{\mathrm{Aut}}

\title{Cyclicity of Drinfeld modules modulo primes in arithmetic progressions}
\author[Prassanna Nand Jha]{Prassanna Nand Jha}

\address{Department of Mathematics, 
	Indian Institute of Technology Gandhinagar, Gandhinagar, 
	Gujarat 382355, India}

\email{prassanna.jha@iitgn.ac.in}

\author[Arnab Saha]{Arnab Saha}

\address{Department of Mathematics, 
	Indian Institute of Technology Gandhinagar, Gandhinagar, 
	Gujarat 382355, India}
\email{arnab.saha@iitgn.ac.in}

\author[Akshaa Vatwani]{Akshaa Vatwani}

\address{Department of Mathematics, 
	Indian Institute of Technology Gandhinagar, 
	Gandhinagar, 
	Gujarat 382355, India}

\email{akshaa.vatwani@iitgn.ac.in}

\date{}

\begin{document}
	
	%\\

	\begin{abstract}
		Let $A=\mathbb{F}_q[T]$,  $K=\mathbb{F}_q(T)$, and let $\phi$ be
		a Drinfeld $A$-module over $K$ of rank $r\geq 2$. We study the reduction
		of $\phi$ at primes $\fp\in A$ lying in a prescribed arithmetic
		progression modulo $\ft$. We obtain asymptotic formulas for the number
		of such primes of degree $x$ for which $\phi(\mathbb{F}_{\fp})$ has at
		most $r-1$ cyclic components, with error terms depending on $r$
		and $\deg\ft$. 
		Under additional hypotheses, we obtain a range for $\deg \ft$
		in terms of $x$ and $r$,  allowing $\deg\ft$ to grow as $x \to \infty$,  in which, the
		error term admits a uniform power saving relative to the main term.

	We also investigate the positivity of the density for the primes counted above. 	
	More precisely, when $\End_{\overline K}(\phi)=A$ and the relevant
		composita of division fields are geometric, we establish a necessary
		and sufficient condition for the positivity of the resulting density.  Our analysis of positivity highlights some differences between the elliptic curve and Drinfeld module settings, particularly in the roles played by cyclotomic fields and  division fields of Carlitz modules, in these respective settings.
		\normalcolor
		
%		Let $A=\mathbb{F}_q[T]$ and $K=\mathbb{F}_q(T)$,		and let $\phi$ be a Drinfeld $A$-module over $K$		of rank $r\geq 2$. We study the reduction of		$\phi$ at primes $\fp\in A$ lying in a prescribed		arithmetic progression $\pmod \ft$, for which		$\phi(\mathbb{F}_{\fp})$ has at most $r-1$	cyclic components. We obtain asymptotic formulas		for the number of such primes of a given degree,	with  error terms depending on the rank $r$ and 	 the degree of the modulus $\deg \ft$.  Under certain additional hypotheses, we look at what happens when deg t grows with x.  As $\deg \ft$ grows with $x$, we determine  ranges of $\deg \ft$, in terms of $x$ and $r$, in which the asymptotic formula holds with uniformly bounded  error terms.  In particular, these error terms admit power-savings relative to the main term.  
%	Under additional hypotheses, we obtain an asymptotic formula uniform for moduli in a rank-dependent range, with an error term that is exponentially smaller than the main term.
	%	Furthermore, when $\End_{\overline K}(\phi)=A$,		under additional hypotheses, we establish a	necessary and sufficient condition for positivity	of the resulting density.
	\end{abstract}
	
	\maketitle

	\section{Introduction and statement of results}\label{sec:intro}
	
	A classical theme in arithmetic is to understand how the algebraic structure associated to an arithmetic object varies after reduction modulo primes. Let $E/\mathbb{Q}$ be an elliptic curve with conductor $N(E)$. For every prime $p\nmid N(E)$, the reduction $\overline{E}/\F_p$ is again an elliptic curve, and the set $\overline{E}(\F_p)$ of its rational points forms a finite abelian group. It is therefore natural to ask how the group structure of $\overline{E}(\F_p)$ varies as $p$ ranges over  primes of good reduction. In particular, one may count   how often $\overline{E}(\F_p)$ is cyclic using the  function
	\begin{align}\label{def:pi-E-x}
		\pi_E(x)\coloneqq\#\{p\leq x:(p,N(E))=1,\ \overline{E}(\F_p)\text{ is cyclic}\}.
	\end{align}
	The cyclicity problem has its origins in the work of Borosh, Moreno and Porta \cite{Borosh-Moreno-Porta}, Lang, Trotter \cite{Lang-Trotter} and Serre  \cite{Serre1978} (see  Section \ref{sec: lit-cyclicity} for a detailed description).

	For $m \in \mathbb N$, let
	$E[m]$
	denote the group of $m$-torsion points of $E$, and let $\mathbb{Q}(E[m])$ be the $m$th division field of $E$,  obtained by adjoining to $\mathbb{Q}$ the coordinates of all the points in $E[m]$. 	Under the generalized Riemann hypothesis (GRH), Serre showed that
	$
	%	\begin{align}\label{eq:Serre-cyclicity}
		\pi_E(x)=\delta_E\operatorname{li}(x)+o\left({x}/{\log x}\right), 
		%	\end{align}
	$
	where $\operatorname{li}(x)=\int_2^x\frac{dt}{\log t}$, and  $\delta_E$ is given by the density 
	\begin{align}\label{def:delta-E}
		\delta_E\coloneqq\sum_{m\geq 1}\frac{\mu(m)}{[\mathbb{Q}(E[m]):\mathbb{Q}]},
	\end{align}
	where  $\mu$ is the M\"{o}bius function. 
	Significant progress has been made on the cyclicity problem over the years - we discuss this in Section \ref{sec: lit-cyclicity}. 
	Recently, Akbal and G\"ulo\u{g}lu \cite{Akbal-Guloglu} imposed an additional arithmetic progression condition on the primes in \eqref{def:pi-E-x}. Letting $f\in \mathbb N$ and $a\in(\mathbb{Z}/f\mathbb{Z})^*$, they obtained  asymptotic estimates for  
	\begin{align}\label{def0:pi-E-fa}
		\pi_E(x;f,a)\coloneqq\#\{p\leq x:(p,N(E))=1,\ p\equiv a\pmod f,\ \overline{E}(\F_p)\text{ is cyclic}\}. 
	\end{align} 
	We give a detailed survey of the surrounding literature in Section \ref{sec:literature}.  
	
	Function fields provide a particularly attractive setting for the study of such  problems. The polynomial ring $\mathbb{F}_q[T]$ exhibits a close arithmetic analogy with the integers, while the distribution of prime polynomials is often more explicit and allows sharper counting estimates. It is therefore natural to ask for function-field analogues of classical questions concerning reductions of elliptic curves. This naturally leads to Drinfeld modules, whose arithmetic provides a framework for formulating and studying such questions over function fields.
	
	\subsection{Notation and setup}
	\label{sec:notation}
	Let $\F_q$ be a finite field of characteristic $p$,  $A=\F_q[T]$ and $K=\F_q(T)$.   Let $\tau:x\mapsto x^q$ be the Frobenius endomorphism of the additive group scheme  $\Ga$ over $\Spec K$.  Let $i:A\to K$ be the canonical inclusion.  Let $\End_K(\Ga)$ be the $K$-endomorphism ring $\Ga$ over $\Spec K$.  Define the subring $K\{\tau\}$ of $\End_K(\Ga)$ by
	$	
	%\begin{align}
	K\{\tau\}\coloneqq \left\{\sum_{j=0}^n \fa_j\tau^j:n\geq0,\ \fa_j\in K\right\},
	%	\end{align}
$
where $\fa_j\in K$ is identified with the endomorphism defined by multiplication by $\fa_j$.  For all $\fa\in K$, one has $\tau\fa=\fa^q\tau$, and hence $K\{\tau\}$ is a noncommutative ring, usually called the ring of twisted polynomials.  We have two homomorphisms: $	\epsilon: K \longrightarrow K\{\tau\},  \epsilon(\fa)=\fa,$ and 
$	%\begin{align}\label{def:epsilon-D}
D: K\{\tau\} \longrightarrow K$ with $ D\Big(\sum_{j=0}^n \fa_j\tau^j\Big)=\fa_0.
$%	\end{align}

A Drinfeld $A$-module $\phi$ over $K$ is defined as a ring  homomorphism 
\begin{align}
\phi: A\to K\{\tau\}, \:
\fa \mapsto\phi_\fa,
\end{align} 
such that $i=D\circ\phi$ and $\phi\neq\epsilon\circ i$.  Let $\deg_\tau\phi_\fa$ denote the degree of $\phi_\fa$ in the variable $\tau$.  The rank $r$ of a Drinfeld module $\phi$ is the unique positive integer \cite[Proposition 2.1]{Drinfeld} such that $\deg_\tau\phi_\fa=r\cdot\deg\fa$ for all $\fa\in A$ with $\fa\neq0$.  We write $\phi\simeq_K\psi$ if there exists some $u\in K^*$ such that
$\psi_\fa=u^{-1}\phi_\fa u$ for every $\fa\in A$. The canonical example of a Drinfeld $A$-module over $K$ of rank $1$ is the Carlitz module $\Carl$, defined by
$
\Carl_T=T+\tau.
$

Throughout the paper, by a prime $\fp\in A$ we mean a monic irreducible polynomial in $A$.  A Drinfeld module has good reduction at all but finitely many primes
\cite[Definition 4.10.1]{Goss}. We denote by $\mathcal P_\phi$ the set of
primes of good reduction for $\phi$. For each $\fp\in\mathcal P_\phi$, let
$\F_\fp=A/\fp A$ be the residue field at $\fp$. The reduction of $\phi$
modulo $\fp$ defines an algebra homomorphism
$
A\longrightarrow\F_\fp\{\tau\},
$
and hence an $A$-module structure on $\F_\fp$, which we denote by
$\phi(\F_\fp)$. Then, by the theory of torsion points for Drinfeld modules and that of finitely generated modules over a PID, there exist uniquely determined monic polynomials $\fd_{1,\fp},\ldots,\fd_{r,\fp}\in A$ satisfying $\fd_{i-1,\fp}A\mid \fd_{i,\fp}A$ for $2\leq i\leq r$ such that
\begin{align}\label{eq:phiFp}
\phi(\F_\fp)\simeq \frac{A}{\fd_{1,\fp}A}\times \frac{A}{\fd_{2,\fp}A}\times\cdots\times \frac{A}{\fd_{r,\fp}A}.
\end{align}

Let $K_\infty=\F_q((1/T))$ denote the completion of $K$ at the prime at infinity in $K$.
A sign function is a homomorphism
$\operatorname{sgn}:K_\infty^*\to\F_q^*$, whose
restriction to $\F_q^*$ is the identity (see \cite[Definition 7.2.1(1)]{Goss}).
A Drinfeld $A$-module $\psi$ over $K$ of rank $1$ is called
sign-normalized if, for some sign function $\operatorname{sgn}$,
the leading coefficient of $\psi_\fa$ equals
$\operatorname{sgn}(\fa)$ for every nonzero $\fa\in A$ (see \cite[Definition 7.2.14]{Goss}).
In our setting, this is equivalent to saying that
$\psi_T=T+c\tau$ for some $c\in\F_q^*$.

\subsection{Summary of results}		
We follow the notation set up in  Section \ref{sec:notation}. 
From \eqref{eq:phiFp}, one can  see that the analogue of  Serre's cyclicity problem for a Drinfeld $A$-module $\phi$ over $K$ of rank $r$,   is concerned with  the quantity
\begin{align}\label{def:pi-phi-x}
\pi_\phi(x)\coloneqq
\#\{\fp\in\mathcal{P}_\phi:\deg\fp=x,\ 
\phi(\F_\fp)\text{ has at most }r-1\text{ cyclic components}\},
\end{align}
%	and obtained an asymptotic formula for $\pi_\phi(x)$. Note that 
which for $r=2$  is analogous to \eqref{def:pi-E-x}. We discuss existing results for this problem in Section \ref{sec:lit-drinfeld}. 	It is  natural to study  the density  of primes $\mf p$ lying in a given  arithmetic progression, for which $	\phi(\F_\fp)$ has at most $r-1$ cyclic components. 
%and explore whether we can adopt Akbal and G{\"u}lo{\u{g}}lu's work in the setting of Drinfeld modules. In particular, we establish asymptotic density estimates for
More precisely, 	for a monic polynomial  $\ft\in A$  and  $\fa\in(A/(\ft))^*$, consider the   quantity
\begin{align}\label{def:pi-phi-x,t,a}
 \! \! \! \!  \!\: \pi_\phi(x;\ft,\fa)  \! \coloneqq\#\{\fp\in\mathcal P_\phi  \! : \! \deg\fp=x,  \fp\equiv\fa\pmod{\ft},\  \! \phi(\F_\fp)\text{ has at most }r-1\text{ cyclic components}\}. 
\end{align}
We may let   $\deg \ft \leq x$ since   $\pi_\phi(x;\ft,\fa) \le 1$   if $\deg \ft>x$.   

%then the only monic $\fp \in A$ with $\deg \fp =x$ and $\fp \equiv \fa \pmod \ft$ is $\fp = \fa$. Hence, the count
\textbf{	Cyclicity results:} 
In this paper, we obtain an asymptotic formula for $	\pi_\phi(x;\ft,\fa)$, with error terms depending on $\deg \ft$ and the rank $r$ of $\phi$. This is Theorem \ref{thm:estimation} below. We use  inclusion-exclusion to express the cyclicity condition in terms of a Frobenius condition in  the division fields $K(\phi[\fm])$ of $\phi$. In addition, we encode  {the arithmetic progression   via a Frobenius condition in the division fields $K(\Carl[\ft])$ of the Carlitz module}, so that we finally work with Chebotarev   conditions in  composita of  the form  $L_{\fm,\ft}= K(\phi[\fm])K(\Carl[\ft])$.  We then apply the  Chebotarev density theorem to obtain the main term and  estimate remainder terms to obtain the error term.  

The density $\delta_{\phi,\ft,\fa}(x)$ in
Theorem~\ref{thm:estimation} may depend on the degree $x$ of the prime
polynomials being counted. This dependence arises from the constant
fields of the extensions $L_{\fm,\ft}/K$ (see Corollary \ref{for-main} and the discussion preceding it).  However, it is possible to choose an integer $B = 
B(\phi,K)$, such that for  given \(\mathfrak t,\mathfrak a\),  on each arithmetic progression $x \equiv j \pmod B$, the density is constant, although its value may vary between residue classes. This is expressed more precisely in Remark \ref{rmk:density in AP}. Under additional hypotheses, we determine in Theorem \ref{thm:uniform-AP} a
rank-dependent range for $\deg\ft$ in which the error term admits
a uniform power saving relative to the main term. \normalcolor

\begin{comment}
	In Theorem \ref{thm:estimation-non-CM}, we refine this result under   the condition $\End_{\overline K}(\phi)=A$.  This  is analogous to the non-CM condition for an elliptic curve $E/\mathbb Q$, namely $\End_{\overline{\mathbb Q}}(E)=\mathbb Z$. In this setting, the adelic openness theorem of Pink and Rütsche  \cite[Theorem~0.1]{PinkRutsche}  plays the role of Serre’s open image theorem for non-CM elliptic curves. In particular, the Galois groups of the division fields $K(\phi[\mathfrak m])$ are isomorphic to $\mathrm{GL}_r(A/\mathfrak m A)$ whenever $\mathfrak m$ is coprime to a fixed nonzero polynomial depending on $\phi$ (see Lemma \ref{lem:M-phi}). This stronger information allows us to obtain an asymptotic formula with  improved error term in suitable ranges of the modulus in Theorem \ref{thm:estimation-non-CM}. 
\end{comment}

\textbf{Positivity results:} 
We next  restrict to the case in which the field extension $L_{\mathfrak m,\mathfrak t}/K$ is geometric  (equivalently  $[L_{\fm,\ft}\cap\overline{\F_q}:\F_q]=1$), for every monic $\mathfrak m\in A$. The density $\delta_{\phi,\mathfrak t,\mathfrak a}$ is then independent of $x$, and we study conditions under which $\delta_{\phi,\mathfrak t,\mathfrak a}$  is positive. It is at this point that an important distinction between the elliptic curve and Drinfeld module settings arises due to the dual role played by cyclotomic fields in the positivity problem for  elliptic curves. First, the condition \(p\equiv a\pmod f\) is encoded by the Frobenius action on \(\mathbb Q(\zeta_f)\). Second, one has the containment
$
\mathbb Q(\zeta_f)\subseteq\mathbb Q(E[f])
$
for every \(f\). These inclusions provide the auxiliary subfields and   suitable disjointness results in the positivity arguments of  \cite{cojo-murty}, \cite{Akbal-Guloglu}.  

In the Drinfeld setting, these two roles need not be fulfilled by the same family of fields. The canonical  rank $1$ Drinfeld modules whose division fields give an explicit class field theory are precisely those that are sign-normalized. The  arithmetic progression condition is thus encoded by the division field  \(K(\psi^{\mathrm{Carl}}[\mathfrak t])\) of the  Carlitz  module $\psi^{\mathrm{Carl}}$ (one may replace $\psi^{\mathrm{Carl}}$ by any chosen rank $1$ sign-normalized Drinfeld module). However,   \(K(\psi^{\mathrm{Carl}}[\mathfrak t])\)  need not be contained in \(K(\phi[\mathfrak t])\). The analogue of the Weil pairing instead gives
$
K(\phi^1[\mathfrak m])\subseteq K(\phi[\mathfrak m]),
$
where \(\phi^1\) is the rank \(1\) Drinfeld module associated to \(\phi\) (see \eqref{def0:phi1}), which   need not be sign-normalized and need \textit{not} satisfy $\phi^1 \simeq_K \Carl$.  In particular, when $\phi^1 \not\simeq_K \Carl$, the  division fields encoding the arithmetic progression cannot be identified with the ones naturally contained in  the division fields of \(\phi\), and thus the two roles of the cyclotomic fields in the elliptic case are now played by  distinct entities.  Indeed, one may still obtain positivity of $\delta_{\phi,\mathfrak t,\mathfrak a}$ in such a scenario, as demonstrated by Example \ref{prop:non-Carlitz-positive-intro}. 

%$\phi^1 \simeq_K \Carl$ is not a necessary condition for positivity of $\delta_{\phi,\mathfrak t,\mathfrak a}$, as demonstrated by Example \ref{prop:non-Carlitz-positive-intro}. 

As a consequence,  methods of  \cite{cojo-murty}, \cite{Akbal-Guloglu},  \cite{Kuo-Liu2009} for positivity of density  -- which construct  auxiliary fields using containment and disjointness relations between division fields --   do not transfer  to our setting. 	In this paper, we give a Galois-theoretic argument inspired by Milner and Shotton \cite{MilnerShotton}, which  characterizes  the positivity problem in terms of the existence of a compatible automorphism in a finite Galois extension. This is Theorem \ref{thm:non-CM-density} below. As Example~\ref{prop:non-Carlitz-positive-intro} illustrates,
the criterion given by Theorem \ref{thm:non-CM-density} can establish positivity even when $\phi^1\not\simeq_K\Carl$.

As noted in \cite{MilnerShotton}, the argument of \cite{cojo-murty} requires a correction, since pairwise linear disjointness of the auxiliary fields does not, in general, imply their joint linear disjointness. 
%A similar issue arises in the proof of \cite[Theorem 3]{Kuo-Liu2009}, which therefore also requires correction. 
Our proof of Theorem~\ref{thm:non-CM-density} explicitly accounts
for entanglements among the relevant division fields and gives a necessary and sufficient criterion for positivity of the density. 
When there is no arithmetic progression condition, it yields
Corollary~\ref{cor:positivity for t=1}. This  addresses the
positivity problem considered in \cite{Kuo-Liu2009} and gives  an alternative  criterion for positivity, as compared to \cite[Theorem~3]{Kuo-Liu2009}.

%As noted in \cite{MilnerShotton}, the argument of \cite{cojo-murty} requires a correction, since pairwise linear disjointness of the auxiliary fields does not, in general, imply their joint linear disjointness. A similar issue arises in the proof of \cite[Theorem 3]{Kuo-Liu2009}, which therefore also requires correction. Our proof of Theorem~\ref{thm:non-CM-density} explicitly accounts for entanglements among the relevant division fields. When there is no arithmetic progression condition, it yields Corollary~\ref{cor:positivity for t=1}. This  addresses the positivity problem considered in \cite[Theorem~3]{Kuo-Liu2009} by giving  an alternative necessary and sufficient criterion for positivity. 

Finally, if $	\phi^1\simeq_K\Carl$, the positivity criterion simplifies to yield Corollary \ref{cor:positivity-final-explicit}.

\subsection{Statement of results}		
In order to state our results, we need the following notation. 
Let $\pi(x)$ denote the number of prime polynomials in $A$ of degree $x$, $\mu_q$ denote the M\"obius function on $A$, and $\varphi(\fm)=|(A/\fm A)^*|$.
Let $\ft\in A$ be monic and let $\fa\in (A/(\ft))^*$. For a monic $\fm\in A$, put
\begin{align}\label{def:Lmt-intro}
	L_{\fm,\ft}\coloneqq K(\phi[\fm])K(\Carl[\ft]),\quad
	r_{\fm,\ft}\coloneqq [L_{\fm,\ft}\cap\overline{\F_q}:\F_q].
\end{align}
For Galois extensions $L/K$, define
\begin{align}\label{def:gamma-intro}
	\gamma_{\fa,\ft}(L)\coloneqq
	\begin{cases}
		1,&\Carl_\fa\text{ fixes }L\cap K(\Carl[\ft]),\\
		0,&\text{otherwise}.
	\end{cases}
\end{align}
The condition $\gamma_{\fa,\ft}(K(\phi[\fm]))=1$ detects precisely when there is a unique automorphism
$
\nu_{\fm,\ft}\in\Gal(L_{\fm,\ft}/K),
$
which fixes $K(\phi[\fm])$ and which acts as $\Carl_\fa$ on $\lambda_\ft$, where $\lambda_\ft$ generates the $A/\ft A$-module $\Carl[\ft]$.   Moreover, there exists an integer $0\leq \epsilon_{\fm,\ft}<r_{\fm,\ft}$ such that  the restriction of $\nu_{\fm,\ft}$ to the constant  field $L_{\fm,\ft}\cap\overline{\F_q}$ is the $ \epsilon_{\fm,\ft}$-th power of the Frobenius  
$\tau$, that is, 
\begin{align} 
	\label{def:epsilon}
	\nu_{\fm,\ft}\big|_{L_{\fm,\ft}\cap\overline{\F_q}}=\tau^{\epsilon_{\fm,\ft}}.
\end{align}
Whenever $\gamma_{\fa,\ft}(K(\phi[\fm]))=0$, we set $\epsilon_{\fm,\ft}=0$. Also, for 
$
\fd=\fp_1^{e_1}\cdots\fp_k^{e_k}
$,
we define
$
\fd_0(\fd)\coloneqq\fp_1^{\lceil e_1/r\rceil}\cdots\fp_k^{\lceil e_k/r\rceil}.
$
Then, for a monic $\ft\in A$, set
\begin{align}\label{def:J-alpha-intro}
	J_\phi(\ft)\coloneqq\summonics_{\fd\mid\ft}\varphi(\fd)q^{-\deg\fd_0},\quad
	\alpha_\phi(\ft)\coloneqq\deg\ft-\log_qJ_\phi(\ft).
\end{align}
Let  $\End_{\overline K}(\phi)$ denote the endomorphism ring of the Drinfeld module $\phi$ over the algebraic closure $\overline K$ of $K$. Let $\rho$ denote the rank of $\End_{\overline K}(\phi)$ as an $A$-module. 

For functions $f$ and $g$, with $g\geq 0$, the notation $f\ll_{\phi,K}g$ or $f=O_{\phi,K}(g)$ means that $|f|\leq Cg$ for some constant $C>0$ depending only on $\phi$ and $K$. 
Throughout the paper, let $\sum'$ denote a sum over monics in $A$.

Our first main result is as follows. 

\begin{theorem}\label{thm:estimation}
Let $\ft\in A$ be monic and let $\fa\in(A/(\ft))^*$ be arbitrary.  Let $\phi$ be a Drinfeld $A$-module over $K$ of rank $r\geq2$.
For monics $\fm \in A$, define
\begin{align}\label{def:r-m,t-x}
	r_{\fm,\ft}(x)\coloneqq\begin{cases}
		r_{\fm,\ft}, &\text{if } x\equiv\epsilon_{\fm,\ft}\pmod{r_{\fm,\ft}},\\
		0, &\text{otherwise},
	\end{cases}
\end{align}
where $r_{\fm,\ft}$, $\epsilon_{\fm,\ft}$  are as defined in \eqref{def:Lmt-intro}, \eqref{def:epsilon} respectively. Let $\gamma_{\fa,\ft}$ be as defined in \eqref{def:gamma-intro} and define the density 
\begin{align}\label{def:delta-intro}
	\delta_{\phi,\ft,\fa}(x)
	\coloneqq
	\summonics_{\substack{\fm\in A}}
	\frac{\mu_q(\fm)r_{\fm,\ft}(x)\gamma_{\fa,\ft}(K(\phi[\fm]))}
	{[K(\phi[\fm])K(\Carl[\ft]):K]},
\end{align}
where $\mu_q$ denotes the M\"obius function on $A$. Then, for integers $x\geq\max\{2,\deg\ft\}$,
\begin{align}\label{eq:main-asymptotic}
	\pi_\phi(x;\ft,\fa)=\delta_{\phi,\ft,\fa}(x)\pi(x)+E_{\phi,K}(x,\ft),
\end{align}
where 
\begin{align}\label{eq:error-main}
	E_{\phi,K}(x,\ft)=
	\begin{cases}
		O_K\!\left(\frac{x}{r}q^{\frac{r-1}{r}x}\right)
		+O_{\phi,K}\!\left(
		q^{\frac{r^2+4r-2}{2r(r+1)}x-\frac{\alpha_\phi(\ft)}{r+1}}
		\left(1+\frac{\deg\ft}{x}\right)^{\frac{r}{r+1}}
		\right),
		& \text{if } r=2,3, \\ 
		O_{\phi,K}\!\left(q^{\frac{x}{2}+\frac{x}{r}}\left(1+\frac{\deg\ft}{x}\right)\right),
		& \text{if }r\geq 4,
	\end{cases}		
\end{align}
where $\alpha_\phi(\ft)$ is as defined in \eqref{def:J-alpha-intro}.
\end{theorem}

\begin{remark}\label{rem:density-bound}
	We note that, for every monic $\ft\in A$,
	$\fa\in(A/(\ft))^*$ and $x\in\mathbb N$, we have $0\leq\delta_{\phi,\ft,\fa}(x)
		\leq1/\varphi(\ft)$.
	This is shown at the end of Subsection \ref{subsec:CM-main}.
\end{remark}

\begin{remark}\label{rmk:density in AP}
	By Lemma \ref{E-phi}, $r_{\fm,\ft}<E(\phi,K)$ for all monic
	$\fm,\ft\in A$. Put
	$B=B(\phi,K)\coloneqq\operatorname{lcm}\{1,\ldots,E(\phi,K)\}$.
	Then, for each $0\leq j<B$, define the constant
	\begin{align}
		\delta_{\phi,\ft,\fa,j}
		\coloneqq
		\summonics_{\substack{\fm : 
				\epsilon_{\fm,\ft}			\equiv 	j \pmod{r_{\fm,\ft}}}}
		\frac{\mu_q(\fm)r_{\fm,\ft}
			\gamma_{\fa,\ft}(K(\phi[\fm]))}
		{[K(\phi[\fm])K(\Carl[\ft]):K]}.
	\end{align}
	Let $x \in \mathbb N$ be such that $x\equiv j\pmod{B}$. Since $r_{\fm,\ft}\mid B$, we have that $x\equiv\epsilon_{\fm,\ft}\pmod{r_{\fm,\ft}}$ if and only if $j\equiv\epsilon_{\fm,\ft}\pmod{r_{\fm,\ft}}$. Hence, for $x\geq\max\{2,\deg\ft\}$, satisfying $x\equiv j\pmod B$, we have  $\pi_\phi(x;\ft,\fa)
	=\delta_{\phi,\ft,\fa,j}\,\pi(x)
	+E_{\phi,K}(x,\ft)$. 
\end{remark}

We  mention here that Theorem \ref{thm:estimation} recovers
\cite[Theorem 1]{Kuo-Liu2009} for $K=\F_q(T)$ upon taking $\ft=1$ and $\rho =1$.
Indeed, then $\fa=1$, $\alpha_\phi(1)=0=\epsilon_{\fm,1}$ and
$\gamma_{1,1}(K(\phi[\fm]))=1$, so
$\delta_{\phi,1,1}(x)=\sum'_{\fm\in A}
\mu_q(\fm)r_{\fm}(x)/[K(\phi[\fm]):K]$,
where $r_{\fm}(x)=r_{\fm,1}\mathbbm{1}_{r_{\fm,1}\mid x}$ and
$r_{\fm,1}=[K(\phi[\fm])\cap\overline{\F_q}:\F_q]$.
Moreover, $E_{\phi,K}(x,1)\ll_{\phi,K}q^{\Delta_r x}$, with
$\Delta_r$ matching the error exponent in \cite{Kuo-Liu2009} (see \eqref{exponent kuo-liu}). 

\begin{remark}\label{constant-delta}
	If the field extension $L_{\mathfrak m,\mathfrak t}/K$ is geometric  (equivalently  $r_{\fm,\ft}=1$), for every monic $\mathfrak m\in A$, then the congruence condition in \eqref{def:r-m,t-x} is vacuously true and $\delta_{\phi,\ft,\fa}(x)$ is independent of $x$.  In this case, we denote it as
	\begin{align}\label{def0:constant-delta}
		\delta_{\phi,\ft,\fa}
		\coloneqq
		\summonics_{\fm\in A}
		\frac{\mu_q(\fm)\gamma_{\fa,\ft}(K(\phi[\fm]))}
		{[K(\phi[\fm])K(\Carl[\ft]):K]}.
	\end{align}
	% not needed ---
	%Then by Theorem \ref{thm:estimation} and Theorem \ref{thm:estimation-non-CM},\begin{align}	\pi_\phi(x;\ft,\fa)=\delta_{\phi,\ft,\fa}\pi(x)+o(\pi(x)),\end{align} and hence by the nonnegativity of $\pi_\phi(x;\ft,\fa)$, we have that for any fixed $\ft$ and $\fa$, the density $\delta_{\phi,\ft,\fa}\geq0$. 
\end{remark}

Note that there are exactly $q^{x-\deg\ft}$ monic polynomials of degree $x$ in a prescribed residue class modulo $\ft$, since $\deg\ft\leq x$. Thus, 
$ q^{x-\deg\ft}$ is a trivial bound for $\pi_\phi(x;\ft,\fa)$. 
Indeed, the main term  $\delta_{\phi,\ft,\fa}(x)\pi(x)$ in \eqref{eq:main-asymptotic} aligns with  this --  using the prime number theorem for polynomials and the bound
$\varphi(\ft)\gg_K q^{\deg \ft}/(1+\deg \ft)$, one can see that  
$\delta_{\phi,\ft,\fa}(x)\pi(x) \le \pi(x) / \varphi(\ft) \ll_{\phi,K}q^{x-\deg \ft}$. 
%In particular, the error bound in Theorem \ref{thm:estimation} is better than the trivial bound $q^{x-\deg t}$ whenever $3\deg t-\alpha_\phi(t)<x/2$ for $r=2$, $4\deg t-\alpha_\phi(t)<5x/6$ for $r=3$, and $\deg t<\frac{r-2}{2r}x$ for $r\geq4$.
Due to the varying  size of the main term in \eqref{eq:main-asymptotic}, it can be a delicate task to check whether \eqref{eq:main-asymptotic} is truly an asymptotic formula, especially when the degree of the modulus $\ft$ grows with $x$.  
We hence aim to give an explicit range for $\deg\ft$,
in terms of $x$ and the rank $r$, throughout which the error term
admits a uniform power saving relative to the main term.  To do this, we assume additional hypotheses under which we have $\delta_{\phi,\ft,\fa}(x) >0$ (see Corollary 
\ref{cor:positivity-final-explicit}).

For a Drinfeld $A$-module $\phi$ over $K$ of rank $r$, given by
$
\phi_T \!=\! T + \fa_1 \tau + \cdots + \fa_r \tau^r,
$
the rank $1$ Drinfeld $A$-module over $K$ associated to $\phi$ via the Weil pairing of Tate modules (see \cite{Heiden-Jan} , \cite[Section 3.7]{Papikian}) is given  by
\begin{align}\label{def0:phi1}
	\phi^1_T \coloneqq T + (-1)^{r-1} \fa_r \tau.
\end{align}
Under the assumption	
$\End_{\overline K}(\phi)=A$, 
by a result of Pink and R\"utsche \cite{PinkRutsche}, there exists a nonzero monic polynomial $M(\phi)\in A$ (see Lemma \ref{lem:M-phi}) such that the division fields $K(\phi[\fm])$, for $\fm\in A$ coprime to $M(\phi)$, have full Galois image. Define
\begin{align}
\label{def:N-phi}
N(\phi)\coloneqq\text{lcm}\big(M(\phi),\prod_{\deg \fp = 1}\fp\big). 
\end{align}
%Put $\lambda_2=1/5$, $\lambda_3=5/22$ and $\lambda_r=\frac{r-2}{2r}$ for $r\geq4$. Fix $0<\varepsilon<\lambda_r$. The next corollary establishes an asymptotic formula for $\pi_\phi(x;\ft,\fa)$, under the assumptions of Corollary \ref{cor:positivity-final-explicit}, uniformly for monics $\ft\in A$ with $\deg\ft\leq(\lambda_r-\varepsilon)x$.

\begin{theorem}\label{thm:uniform-AP}
	Let $\phi$  be as in Theorem \ref{thm:estimation}.
	Assume further that $q>2$,
	$\End_{\overline K}(\phi)=A$,  and 
	$\phi^1\simeq_K\Carl$.  
	Put $\lambda_2=1/5$, $\lambda_3=5/22$
	and $\lambda_r=\frac{r-2}{2r}$ for $r\geq4$.
	For every fixed $0<\varepsilon<\lambda_r$,
	we have
	\begin{align}
		\label{asymp}
		\pi_\phi(x;\ft,\fa)
		=\delta_{\phi,\ft,\fa}\pi(x)
		+
		\begin{cases}
			O_{\phi,K}\left(
			q^{5x/6-\alpha_\phi(\ft)/3}
			\right), & r=2,\\[4pt]
			O_{\phi,K}\left(
			q^{19x/24-\alpha_\phi(\ft)/4}
			\right), & r=3,\\[4pt]
			O_{\phi,K}\left(
			q^{(1/2+1/r)x}
			\right), & r\geq4,
		\end{cases}
	\end{align}
	uniformly for all monic $\ft\in A$ and
	$\fa\in(A/(\ft))^*$ satisfying
	$\gcd(\ft,N(\phi))=1$,
	$r_{\fm,\ft}=1$ for every monic $\fm\in A$,
	and $\deg\ft\leq(\lambda_r-\varepsilon)x$. 
	In particular, for $ \fa$ and $\ft$ as above,  we have uniformly 
		\begin{align} 
			\label{asymp r=4}
	\pi_\phi(x;\ft,\fa) = \delta_{\phi,\ft,\fa}\pi(x) \big(1+ O_{\phi, K}(x  q^{-\varepsilon x} \big)\big), 	\end{align}
 for $r\ge 4$ and 
	\begin{align} 
	\label{asymp r=2,3}
\pi_\phi(x;\ft,\fa) = \delta_{\phi,\ft,\fa}\pi(x) \big(1+ O_{\phi, K}(x \tau(\ft)^{1/(r+1)}  q^{-c_r \varepsilon x} \big) \big), 
\end{align}
 with $c_r =\frac{5}{6}, \frac{11}{12}$, for $r=  2, 3$ respectively.  
\end{theorem}

We next study when the density $\delta_{\phi,\ft,\fa}$ is  positive. Our next result  is a corollary to Theorem \ref{thm:estimation}.

\begin{corollary}\label{cor:density-zero}
Let $\ft\in A$ be monic and let $\fa\in(A/(\ft))^*$ be arbitrary. Let $\phi$ be a Drinfeld $A$-module over $K$ of rank $r\geq2$ such that
$r_{\fm,\ft}=1$ for every monic $\fm\in A$. If
$K(\phi[\fp])=K
$
for some prime $\fp\in A$ with $\deg\fp=1$, then
$\delta_{\phi,\ft,\fa}=0$.
\end{corollary}
Note that when we put $\ft =1$, Corollary \ref{cor:density-zero} yields one direction of Theorem 3 of Kuo-Liu \cite{Kuo-Liu2009}.

%Similar results have also been established in \cite{cojo-murty}, \cite{Kuo-Liu2009} and \cite{Akbal-Guloglu} by introducing suitable auxiliary families of division fields. Milner and Shotton \cite{MilnerShotton} recently pointed out that the linear-disjointness argument underlying this approach requires additional care, since pairwise linear disjointness of the auxiliary fields does not in general imply their joint linear disjointness. They repaired the elliptic-curve argument by explicitly accounting for the resulting cyclotomic entanglement. Their repair, however, is specific to the structure of the auxiliary family arising in that setting and does not directly carry over to our situation, where the arithmetic-progression condition introduces further intersections with the Carlitz division field. 

%We therefore follow an alternative Galois-theoretic argument as developed in \cite{MilnerShotton}, which involves isolating all the entanglement relevant to the density. This reduces the positivity problem to the existence of a compatible automorphism in a finite Galois extension.

We now show that assuming  $
\End_{\overline K}(\phi)=A
$, the positivity of $\delta_{\phi,\ft,\fa}$ is equivalent to a condition on division fields. 
Let $\mathcal B$ be the set of all primes in $A$ dividing $N(\phi)$ and $\mathcal G$ be the set of all primes in $A$ coprime to $N(\phi)$. Also define, for a finite set $\mathcal S \subseteq \mathcal G$,
\begin{align}\label{def:DBUS}
D_{\mathcal B \cup \mathcal S} \coloneqq \prod_{\fp \in \mathcal B \cup \mathcal S}K(\phi[\fp]), \quad \text{and} \quad L_{\phi,\ft,\mathcal S}
\coloneqq
D_{\mathcal B\cup\mathcal S}
K(\Carl[\ft]),
\end{align}
where the product denotes compositum of fields. We call a finite set $\mathcal S\subseteq\mathcal G$ an
\emph{entanglement set} for $(\phi,\ft)$ if, for every square-free monic
$\fm$ supported on $\mathcal G\setminus\mathcal S$,
$K(\phi[\fm])\cap
D_{\mathcal B\cup\mathcal S}K(\Carl[\ft])
=K$.
By Lemma \ref{finite-entanglement}, such a set $\mathcal S$ exists. We prove the following positivity criterion. 

\begin{theorem}\label{thm:non-CM-density}
Let $\ft\in A$ be monic and let
$\fa\in(A/(\ft))^*$ be arbitrary. Let $\lambda_\ft$ be a generator of $\Carl[\ft]$, so that
$K(\Carl[\ft])=K(\lambda_\ft)$. Let $\phi$ be a Drinfeld
$A$-module over $K$ of rank $r\geq2$ such that
$r_{\fm,\ft}=1$ for every monic $\fm\in A$, and
$\End_{\overline K}(\phi)=A$. Let $\mathcal S$ be an entanglement set for $(\phi,\ft)$, and define
\begin{align}\label{def:Omega-intro}
	\Omega_{\phi,\ft,\fa,\mathcal S}
	\coloneqq
	\left\{
	\sigma\in
	\Gal(L_{\phi,\ft,\mathcal S}/K):
	\sigma|_{K(\Carl[\ft])}(\lambda_\ft)=\Carl_\fa(\lambda_\ft),
	\sigma|_{K(\phi[\fp])}\neq\mathbbm{1}
	\text{ for every }\fp\in\mathcal B\cup\mathcal S
	\right\}.
\end{align}
Then
\begin{align}\label{density-finite-factor-intro}
	\delta_{\phi,\ft,\fa}
	=
	\frac{|\Omega_{\phi,\ft,\fa,\mathcal S}|}
	{[L_{\phi,\ft,\mathcal S}:K]}
	\prod_{\fp\notin\mathcal B\cup\mathcal S}
	\left(
	1-\frac{1}{[K(\phi[\fp]):K]}
	\right).
\end{align}
In particular, $\delta_{\phi,\ft,\fa}>0$ if and only if $\Omega_{\phi,\ft,\fa,\mathcal S}$ is nonempty.
\end{theorem}

In the case $\ft =1$, let us denote the density \eqref{def0:constant-delta} by $\delta_{\phi}$.  We also  have $
K(\psi^{\mathrm{Carl}}[1])=K,
$ and $L_{\phi,\ft,\mathcal S}= L_{\phi,\mathcal S}=\prod_{\fp \in \mathcal B \cup \mathcal S}K(\phi[\fp])$,    so that the set $\Omega_{\phi,\ft,\fa,\mathcal S}$ reduces to 
\begin{align}
\Omega_{\phi, \mathcal S} = \{ \sigma 
\in
\Gal(L_{\phi,\mathcal S}/K):
\sigma|_{K(\phi[\fp])}\neq\mathbbm{1}
\text{ for every }\fp\in\mathcal B\cup\mathcal S
\}.
\end{align}
We thus obtain the following corollaries. 
\begin{corollary} \label{cor:positivity for t=1}
If $\phi$ is  a Drinfeld
$A$-module over $K$ of rank $r\geq2$ with  $\End_{\overline K}(\phi)=A$ and all  division fields of $\phi$ are geometric, then 
 $\delta_{\phi}>0 \iff  \Omega_{\phi, \mathcal S}  \neq \varnothing$. 
\end{corollary}

%As a corollary to Theorem \ref{thm:non-CM-density}, we have the following result. 
%\begin{corollary}\label{cor:ositivity-final}
%Let $q>2$. Let $\ft\in A$ be monic and let $\fa\in(A/(\ft))^*$ be arbitrary. Let $\phi$ be a Drinfeld $A$-module over $K$ of rank $r\geq2$ such that $	\End_{\overline K}(\phi)=A$ and $r_{\fm,\ft}=1$ for every monic $\fm\in A$. Assume further that  $	\gcd(\ft,N(\phi))=1$ and $	\phi^1\simeq_K\Carl.$ Then the density $\delta_{\phi,\ft,\fa}>0$.
%\end{corollary}

\begin{corollary}\label{cor:positivity-final-explicit}
	Let $q>2$. Let $\ft\in A$ be monic and let
	$\fa\in(A/(\ft))^*$ be arbitrary. Let $\phi$ be a
	Drinfeld $A$-module over $K$ of rank $r\geq2$
	such that $\End_{\overline K}(\phi)=A$ and
	$r_{\fm,\ft}=1$ for every monic $\fm\in A$.
	Assume further that
	$
	\gcd(\ft,N(\phi))=1
	$
	and
	$
	\phi^1\simeq_K\Carl.
	$
	Then the density $\delta_{\phi,\ft,\fa}>0$.
	In particular, the density is given by
	\begin{align}
		\delta_{\phi,\ft,\fa}
		=
		\frac{C_\phi}{\varphi(\ft)}
		\prod_{\fp\nmid N(\phi)\ft}
		\left(
		1-\frac{1}{|\GL_r(A/\fp A)|}
		\right)
		\prod_{\substack{\fp\mid\ft \\ \fp \mid \fa-1}}
		\left(
		1-
		\frac{1}
		{|\SL_r(A/\fp A)|}
		\right),
		\label{density-euler-product}
	\end{align}
	where $C_\phi
	\coloneqq
	\sum_{\fd\mid\fb}
	\frac{\mu_q(\fd)}{[K(\phi[\fd]):K]}
	>0$ and $\fb$ is the product of all primes dividing $N(\phi)$.  
\end{corollary}

\begin{remark}
	In the case $\ft =1$, by Corollary \ref{cor:positivity-final-explicit}, we have that
	\begin{align}
		\delta_{\phi}
		=
		C_\phi
		\prod_{\fp\nmid N(\phi)}
		\left(
		1-\frac{1}{|\GL_r(A/\fp A)|}
		\right)=\sum_{\fd\mid\fb}
		\frac{\mu_q(\fd)}{[K(\phi[\fd]):K]}\prod_{\fp\nmid N(\phi)}
		\left(
		1-\frac{1}{|\GL_r(A/\fp A)|}
		\right),
	\end{align}
	and hence 
	\begin{align}
		\delta_{\phi,\ft,\fa}
		=
		\frac{\delta_{\phi}}{\varphi(\ft)}
		\prod_{\fp\mid\ft}
		\left(
		1-\frac{1}{|\GL_r(A/\fp A)|}
		\right)^{-1}
		\prod_{\substack{\fp\mid\ft\\\fp\mid\fa-1}}
		\left(
		1-\frac{1}{|\SL_r(A/\fp A)|}
		\right).
	\end{align}
	Thus, under the hypothesis of Corollary \ref{cor:positivity-final-explicit}, $\delta_{\phi,\ft,\fa}>0$ if and only if $\delta_\phi>0$.
\end{remark}

\normalcolor

We apply the positivity criterion in
Theorem~\ref{thm:non-CM-density} to  give the following example,  demonstrating that $\delta_{\phi,\ft,\fa}$ can be positive
even when $\phi^1\not\simeq_K\Carl$.  
\begin{example}\label{prop:non-Carlitz-positive-intro}
Let $q=5$, and let $\phi$ be the rank-$2$ Drinfeld $A$-module over
$K$ defined by
\begin{align}
	\phi_T=T-T\tau-T^{10}\tau^2.
\end{align}
Then, the associated rank $1$ Drinfeld module 	$\phi^1_T=T+T^{10}\tau$ is not  sign-normalized and $\phi^1 \not\simeq_K\Carl$. It satisfies   $\End_{\overline K}(\phi)=A$ and $r_{\fm,T}=1$ for every monic $\fm \in A$. Moreover, for every $\fa\in(A/(T))^*$,
$\delta_{\phi,T,\fa}>0$.
\end{example}

The paper is organized as follows. Section \ref{sec:literature} gives an account of the surrounding literature on the cyclicity problem for elliptic curves and Drinfeld modules. In Section \ref{sec:prelims1}, we recall the relation between the number of cyclic components of $\phi(\F_\fp)$ and complete splitting in the division fields of $\phi$. We also translate the arithmetic progression condition into a Frobenius condition in a Carlitz division field, recall the Chebotarev density theorem for function fields, and establish the uniform estimates required in the counting argument. Sections \ref{sec:estimation-proof} and \ref{sec:uniform} are devoted to the proofs of Theorem \ref{thm:estimation} and Theorem \ref{thm:uniform-AP} respectively. In Section \ref{sec:non-CM}, we consider the case when $\End_{\overline{K}}(\phi)=A$ and obtain an improved estimate for the tail of the density. Finally in Section \ref{sec:density}, we study the positivity of the density and give  proofs of Corollary \ref{cor:density-zero}, Theorem \ref{thm:non-CM-density}, Corollary \ref{cor:positivity-final-explicit} and Example \ref{prop:non-Carlitz-positive-intro}.

	\section{Surrounding literature}\label{sec:literature}
	
	\subsection{Cyclicity of elliptic curves}
	\label{sec: lit-cyclicity}
	
	In 1975, Borosh, Moreno and Porta \cite{Borosh-Moreno-Porta} conjectured
	that, for 'many' elliptic curves $E/\Q$, the group $\overline{E}(\F_p)$
	is cyclic for infinitely many primes $p$. Subsequently, in 1977, Lang and Trotter
	\cite{Lang-Trotter} formulated an elliptic curve analogue of Artin's
	primitive root conjecture. 
	%: for a fixed point $P\in E(\Q)$ of infinite	order, the set of primes $p$ of good reduction	for which $P\bmod p$ generates $\overline{E}(\F_p)$ has a density. 
	The cyclicity conjecture is a subproblem of this, and concerns finding an asymptotic density estimate for $\pi_E(x)$ (see \eqref{def:pi-E-x}). 
	In 1978,   assuming GRH for the Dedekind zeta functions of the division fields of $E/\Q$, Serre \cite{Serre1978} established the asymptotic formula \begin{align}
		\pi_E(x)=\delta_E\operatorname{li}(x)+o({x}/{\log x}),
		\label{error serre}
	\end{align}
	where $\delta_E$ is as  in \eqref{def:delta-E}.
	In 1983, Murty \cite{Murty1983} refined Serre's argument, obtaining an improved error term under GRH. He also proved the asymptotic formula unconditionally for elliptic curves with complex multiplication (CM). Subsequently, in 1987, Murty \cite{Murty1987} proved unconditionally
	the infinitude of cyclic reductions for certain non-CM
	elliptic curves, using the infinitude of supersingular primes. In 1990, Gupta and Murty \cite{Gupta-Murty1990} established unconditionally that for an arbitrary elliptic curve $E/\Q$ satisfying $\Q(E[2])\neq\Q$, we have
	$\pi_E(x)\gg_E x/(\log x)^2$.

	Subsequent improvements by relaxing the assumption of GRH for non-CM elliptic curve and obtaining better error terms for  CM elliptic curves were obtained by Cojocaru in \cite{Cojocaru2002} and \cite{Cojocaru2003} respectively. 
%	\red{	In 2002, for non-CM elliptic curves, Cojocaru \cite{Cojocaru2002} relaxed the assumption of GRH to the condition that the Dedekind zeta functions of the relevant division fields have no zeros in $\Re(s)>3/4$. For CM elliptic curves, she gave a new unconditional proof of the asymptotic formula with explicit error estimates in \cite{Cojocaru2003}.}
 Cojocaru and Murty \cite{cojo-murty} improved the error term in \eqref{error serre} to  $O_{N(E)}(x^{5/6}(\log x)^{2/3})$ under GRH for a non-CM elliptic curve $E/\Q$;  and to  $O\left(x^{3/4}(\log(N(E)x))^{1/2}\right)$ if 
%	\begin{align}\label{eq:CM-survey-nonCM}		\pi_E(x)		=\delta_E\operatorname{li}(x)		+O_{N(E)}\Big(x^{5/6}(\log x)^{2/3}\Big),	\end{align}
	$E$ has complex multiplication by the full ring of integers of an imaginary quadratic field. 
%	\begin{align}\label{eq:CM-survey-CM}		\pi_E(x)	=\delta_E\operatorname{li}(x)		+O\left(x^{3/4}(\log(N(E)x))^{1/2}\right).	\end{align}
	Further assuming Artin's Holomorphy Conjecture and a Pair Correlation Conjecture for the associated Artin $L$-functions, they further improved the exponent $5/6$  to $7/10$, up to a logarithmic factor. Their work also gives bounds for the least prime of cyclic reduction in terms of $N(E)$. Later results include average cyclicity estimates by Banks and Shparlinski \cite{Banks-Shparlinski2009} in 2009 and a Siegel--Walfisz-type theorem for CM elliptic curves by Akbary and V.~K.~Murty \cite{Akbary-Murty2010} in 2010.

	The positivity of the density $\delta_E$ was first studied by Cojocaru and Murty in \cite{cojo-murty}, following an observation of Serre. 
	%Recently, in their work on densities associated with the prescribed smallest invariant factor for non-CM elliptic curves,   Milner and Shotton \cite{MilnerShotton}    repaired a gap in their argument. 
	In 2017, Brau \cite{Brau2017} obtained nonvanishing criteria for $\delta_E$ using character sums for entanglement correction factors. More generally, in 2023, Campagna and Stevenhagen \cite{Campagna-Stevenhagen2023} studied the density of cyclic reductions for elliptic curves over number fields.

	\subsection{Cyclicity of elliptic curves in arithmetic progressions}
	\label{sec:lit-cyclicity in AP}
	
	Let $f\geq 1$ and $a\in(\Z/f\Z)^*$. In 2017, Brau \cite{Brau2017}  studied cyclicity density constants in arithmetic progressions using character sums, obtaining explicit formulas and nonvanishing results. Subsequently, Akbal and G\"ulo\u{g}lu \cite{Akbal-Guloglu} investigated the counting function $\pi_E(x;f,a)$ defined in \eqref{def0:pi-E-fa}.

	%Let $\zeta_f$ be a primitive $f$-th root of unity, and let $\sigma_a\in\Gal(\Q(\zeta_f)/\Q)$ be defined by $\sigma_a(\zeta_f)=\zeta_f^a$. For each integer $m\geq1$, put $\gamma_{a,f}(\Q(E[m]))=1$ whenever $\sigma_a$ fixes $\Q(E[m])\cap\Q(\zeta_f)$, and zero otherwise. The corresponding density constant is given by
%	\begin{align}\label{def:delta-E-fa-literature}		\delta_E(f,a)		\coloneqq		\sum_{m\geq1}	\frac{\mu(m)\gamma_{a,f}(\Q(E[m]))}		{[\Q(E[m])\Q(\zeta_f):\Q]}.	\end{align}
	Under GRH for the relevant Dedekind zeta functions, Akbal and G\"ulo\u{g}lu \cite{Akbal-Guloglu} established
	\begin{align}
		\pi_E(x;f,a)
		=\delta_E(f,a)\operatorname{li}(x)+E_{E,f,a}(x),
	\end{align}
with the density $\delta_E(f,a)$   and  bounds for errors $E_{E,f,a}(x)$ explicitly given, for both non-CM elliptic curves and CM elliptic curves with multiplication by the full ring of integers of an imaginary quadratic field. They also obtained sufficient conditions for the positivity of $\delta_E(f,a)$.
	
	In 2024, Wong \cite{Wong} obtained an unconditional asymptotic formula for CM elliptic curves in arithmetic progressions. Under GRH, he also refined the estimates of Akbal and G\"ulo\u{g}lu, obtaining logarithmic savings for suitable ranges of $f$.	More recently, Lee, Mayle and Wang \cite{Lee-Mayle-Wang2025} and Lee \cite{Lee2026} investigated average aspects of cyclicity in arithmetic progressions, while Jones and Lee \cite{Jones-Lee2026} studied the vanishing of the associated densities.
	
	\subsection{Cyclicity of Drinfeld modules}
	\label{sec:lit-drinfeld}
	In 2009, Kuo and Liu \cite{Kuo-Liu2009} studied the analogue
	of Serre's cyclicity problem for Drinfeld $A$-modules over finite
	extensions of $K$. Let $F/K$ be a finite extension, and let
	$\phi$ be a Drinfeld $A$-module over $F$ of rank $r\geq2$.
	Let $\F_F$ denote the constant field of $F$, and for monic
	$\fm\in A$, define
	$
	r_\fm\coloneqq
	[F(\phi[\fm])\cap\overline{\F}_q:\F_F].
	$
	Kuo and Liu defined the density
$\label{def:delta-KL-literature}
		\delta_{\phi,F}(x)
		\coloneqq
		\sum'_{\substack{\fm\in A}}
		{\mu_q(\fm)r_\fm(x)}/ {[F(\phi[\fm]):F]},
$
	where $r_\fm(x)=r_\fm$ whenever $r_\fm \mid x$, and zero otherwise. Let $\pi_F(x)$ denote the number of primes of $F$ of degree $x$,
	and let $\pi_{\phi,F}(x)$ denote the number of primes of good
	reduction of degree $x$ for which $\phi(\F_\fp)$ has at most
	$r-1$ cyclic components. Assuming $\End_{\overline F}(\phi)=A$,
	they established
$
		\pi_{\phi,F}(x)
		=
		\delta_{\phi,F}(x)\pi_F(x)
		+
		O_{\phi,F}(q_F^{\Delta_r x}),
$
	where $q_F=|\F_F|$ and
	\begin{align} \label{exponent kuo-liu}
		\Delta_r
		\coloneqq
		\begin{cases}
			\frac{r^2+4r-2}{2r^2+2r}, & r=2,3,\\
			\frac{r+2}{2r}, & r\geq4.
		\end{cases}
	\end{align}
	They also investigated the positivity of $\delta_{\phi,F}(x)$, under the hypotheses that $r_\fm =1 $ for every $\fm \in A$.
	
	In 2013, Cojocaru and Shulman \cite{cojo-shulman-I} established
	an unconditional average effective Chebotarev density
	theorem for the division fields 
	of   rank-$2$ Drinfeld modules $\phi$  over $K=\F_q(T)$, with $q$ odd,  such that $\End_K(\phi)$ is the full ring of integers of an
	imaginary quadratic extension of $K$.   As an
	application, they obtained an asymptotic formula
	for cyclic reductions, with an error term of
	$O_\phi(q^{3x/4}\log x)$. Subsequently, in 2015, Cojocaru and Shulman
	\cite{cojo-shulman-II} studied the distribution of the first
	elementary divisor $\fd_{1,\fp}$ (see \eqref{eq:phiFp}) of the reductions of generic
	Drinfeld modules of arbitrary rank. For rank $2$, they further studied the
	average growth of the second elementary divisor $\fd_{2,\fp}$.
	
	In 2017, Virdol \cite{Virdol}  generalized  the result in \cite{Kuo-Liu2009}
  to Drinfeld modules with more general
	endomorphism rings, obtained better  error terms and  established an analogue of the Titchmarsh
	divisor problem for Drinfeld modules.
	In 2021, Kuo and Tweedle \cite{Kuo-Tweedle} studied
	the full cyclicity problem for Drinfeld modules of arbitrary
	rank. They established the existence of a density for the primes
	$\fp$ for which $\phi(\F_\fp)$ is cyclic.	They further proved that this density is positive if and
	only if $\phi(\F_\fp)$ is cyclic for infinitely many
	primes $\fp$.

	\section{Notation and Preliminaries}\label{sec:prelims1}
	
	\subsection{Basic notation}
	
	We have the following notation
	\begin{itemize}
		\item $\pi(x)$: the number of primes in $A$ of degree $x$;
		\item $\F_q$: the field of characteristic $p$ having $q$ elements;
		\item $A=\F_q[T]$: the polynomial ring over $\F_q$ in the variable $T$;
		\item $\mu_q$: the M\"obius function on $\F_q[T]$;
		\item $K=\F_q(T)$: the field of fractions of $A$;
		\item $\F_L=L\cap\overline{\F_q}$: the constant field of a finite extension $L/K$;
		\item $\infty$: the prime at infinity in $K$;
		\item $K_\infty=\F_q((1/T))$: the completion of $K$ at $\infty$;
		\item $C_\infty$: the completion of $\overline{K_\infty}$ at $\infty$.
	\end{itemize}
	
	\subsection{Characteristic polynomials of Frobenius}
	
	We will require some properties of the characteristic polynomial of the Frobenius map on the Tate module of $\phi$. For a prime $\fl \in A$, consider the $\fl$-adic completion of $A$, denoted by
	$A_\fl \coloneqq \varprojlim A/(\fl^n)$.
	Let $K_\fl \coloneqq \mathrm{frac}(A_\fl)$ be the fraction field of $A_\fl$. For any Drinfeld module $\phi$ of rank $r$, the $\fl$-adic \textit{Tate module} $T_\fl(\phi)$ is defined as 
	$T_\fl(\phi)\coloneqq \varprojlim \phi[\fl^n]$,
	where the inverse limit is taken with respect to the morphisms $\phi[\fl^{n+1}]\to \phi[\fl^n]$, given by $\alpha \mapsto \phi_\fl(\alpha)$.
	Then $T_\fl(\phi)\simeq A_\fl^r$.
	
	Let $\fp \in \mathcal P_\phi$ be coprime to $\fl$ and let $\overline{\phi}$ denote the reduction mod $\fp$ of our Drinfeld module $\phi$. Since $(\fp,\fl)=1$, we have that the natural reduction mod $\fp$ induces the isomorphism of free $A_\fl$-modules $A_\fl^r \simeq T_\fl(\phi) \simeq T_\fl(\overline{\phi})$.
	Therefore, the Frobenius map $\eta_\fp$ at $\fp$ given by
	$x \mapsto x^{|A/\fp|}$,
	induces an endomorphism on $T_\fl(\overline{\phi})$ and hence on $T_\fl(\phi)$ via the above isomorphism. The resulting induced action of $\eta_\fp$ on the $K_\fl$-vector space $T_\fl(\phi)\otimes_{A_\fl}K_\fl$ is considered to be the action of the Frobenius at $\fp$ on $T_\fl(\phi)\otimes K_\fl$.
	
	Let $P_{\phi,\fp}^\fl(X)$ denote the  characteristic polynomial of $\eta_\fp$, defined as
	\begin{align}
		P_{\phi,\fp}^\fl(X) \coloneqq \det(X\text{ Id}-{\eta}_\fp)
		= X^r+\fa_{\fp,r-1}X^{r-1}+\cdots+\fa_{\fp,1}X+\fa_{\fp,0}.
	\end{align}
	We state some of the properties of this polynomial, which will be beneficial in obtaining our estimation. We will recall the following properties as stated in \cite[Proposition 8]{cojo-shulman-I}.
	
	\begin{lemma}(see \cite[Corollary 3.4, p.193; Theorem 5.1, p.199]{Gekeler})\label{lem:char_pol}
		Let $\phi$ be a Drinfeld $A$-module over $K$ of rank $r\geq 1$. Let $\fp \in \mathcal{P}_\phi$ and $\fl$ be primes in $K$ such that $\fl \neq \fp$. Then,
		\begin{enumerate}
			\item[(i)] $P_{\phi,\fp}^\fl(X) \in A[X]$; in particular, $P_{\phi,\fp}^\fl(X)$ is independent of $\fl$, and as such, we will drop the superscript $\fl$ from notation and simply write $P_{\phi,\fp}(X)$.
			\item[(ii)] $\fa_{\fp,0}=u_\fp \fp$ for some $u_\fp \in \F_q^*$.
			\item[(iii)] $\deg \fa_{\fp, i} \leq \frac{(r-i)\deg \fp}{r}$, for $i=0,\cdots,r-1$. 
			\item[(iv)] $P_{\phi,\fp}(1)A=d_{1,\fp}d_{2,\fp}\cdots d_{r,\fp} A$, where $d_{i,\fp}$ are as defined in \eqref{eq:phiFp}.
		\end{enumerate}
	\end{lemma}

	For a positive integer $k$, the cyclotomic field $\Q(\zeta_k)$ is a subset of the $k$-division field $\Q(E[k])$. Recall that for a Drinfeld $A$-module $\phi$ over $K$ of rank $r$, defined by
$
		\phi_T = T + \fa_1 \tau + \cdots + \fa_r \tau^r,
$
	the rank $1$ Drinfeld $A$-module $\phi^1$  over $K$, associated via Weil pairing, is given by \eqref{def0:phi1}. 
	Then, we have that $K(\phi^1[\fm])\subseteq K(\phi[\fm])$, for every $\fm \in A$ (see \cite[Theorem 3.7.1]{Papikian}). With this in place, we state the following result of \cite{cojo-shulman-I}.
	\begin{proposition}\cite[Proposition 10]{cojo-shulman-I}\label{char_pol_2}
		Let $\phi$ be a Drinfeld $A$-module over $K$ of rank $r\geq 2$. There exists a rank $1$ Drinfeld $A$-module $\phi^1$ over $K$, uniquely determined up to $\overline{K}$-isomorphism, such that:
		\begin{enumerate}
			\item[(i)] $\mathcal{P}_\phi \subseteq \mathcal{P}_{\phi^1}$;
			\item[(ii)] for any $\fp \in \mathcal{P}_\phi$, the characteristic polynomials of $\phi$ and $\phi^1$ at $\fp$ satisfy the relation:
			\begin{align}
				P_{\phi,\fp}(X) = X^r+\fa_{\fp,r-1}X^{r-1}+\cdots+\fa_{\fp,1}X+u_\fp \fp,\quad
				P_{\phi^1,\fp}(X) = X+(-1)^{r-1}u_\fp \fp,\label{eq:char-pols}
			\end{align}
			where, we recall, $u_\fp \in \F_q^*$;
			\item[(iii)] for any $\fp \in \mathcal{P}_\phi$ and any nonzero $\fm\in A$ coprime to $\fp$, we have that, if $\fp$ splits completely in $K(\phi[\fm])$, then
			\begin{enumerate}
				\item $\fp$ also splits completely in $K(\phi^1[\fm])$;
				\item $\fm^r \mid P_{\phi,\fp}(1)$;
				\item $\fm \mid P_{\phi^1,\fp}(1)$.
				
			\end{enumerate}
		\end{enumerate}
	\end{proposition}
	
	\begin{corollary}\label{cor:for-error}
		Let $\phi$ be a Drinfeld $A$-module over $K$ of rank $r\geq 2$ and let $\fm \in A\setminus\{0\}$. If $\fp\in \mathcal{P}_\phi$ is coprime to $\fm$, and $\fp$ splits completely in $K(\phi[\fm])$, then there exist $u_\fp \in \F_q^*$ and $\fb_\fp\in A$ with $\deg \fb_\fp\leq \frac{(r-1)\deg\fp}{r}$,  satisfying $\fm^r \mid 1+\fb_\fp+u_\fp\fp$ and \: $\fm\mid 1+(-1)^{r-1}u_\fp\fp$.
	\end{corollary}
	\begin{proof}
		Consider the expressions \eqref{eq:char-pols} for the characteristic polynomials $P_{\phi,\fp}(X)$ and $P_{\phi^1,\fp}(X)$.
		We define $\fb_\fp \coloneqq \fa_{\fp,1} + \cdots + \fa_{\fp,r-1}$. Then, by Lemma \ref{lem:char_pol}(iii), $\deg \fb_\fp \leq \frac{(r-1)\deg \fp}{r}$. By Proposition \ref{char_pol_2}(iii), if $\fp \in \mathcal{P}_\phi$ splits completely in $K(\phi[\fm])$, we have
$
			\fm^r \mid P_{\phi,\fp}(1)= 1+\fb_\fp+u_\fp \fp$ and 
		$	\fm \mid P_{\phi^1,\fp}(1)= 1+(-1)^{r-1}u_\fp \fp,
$
		as required.
	\end{proof}

	\subsection{Artin symbols and the Chebotarev Density Theorem}
	
	The following lemma of \cite{Kuo-Liu2009} relates  the $A$-module structure of $\phi(\F_\fp)$  to the splitting of primes in certain division fields.
	
	\begin{lemma}\cite[Proposition 11]{Kuo-Liu2009}\label{lem}
		Let $\phi$ be a Drinfeld $A$-module over $K$ of rank $r\geq 2$, let $\fp\in\mathcal{P}_\phi$, and let $\fl\in A$ be a prime such that $\fl\neq\fp$. Then, the following are equivalent.
		\begin{enumerate}
			\item $\phi(\F_\fp)$ contains an $A$-submodule isomorphic to $(A/\fl A)^r$.
			\item The prime $\fp$ splits completely in $K(\phi[\fl])$.
		\end{enumerate}
		Moreover, the $A$-module $\phi(\F_\fp)$ contains at least $r$ cyclic components if and only if the equivalent conditions above hold for some prime $\fl \neq \fp$.
	\end{lemma}
	
	By the Chinese Remainder Theorem, we have an $A$-module isomorphism
	$\frac{A}{\fm A} \cong \bigoplus_{\fl \mid \fm} \frac{A}{\fl A}$,
	where $\fm$ is a square-free monic in $A$. Since for every nonzero $\fm \in A$, $\phi[\fm]\cong (A/\fm A)^r$, we obtain $\phi[\fm] \cong \bigoplus_{\fl \mid \fm} \phi[\fl]$.
	This implies that we have the compositum
	$K(\phi[\fm]) = \prod_{\fl \mid \fm} K(\phi[\fl])$.
	Therefore, a prime $\fp$ splits completely in $K(\phi[\fm])$ if and only if it splits completely in $K(\phi[\fl])$ for every prime $\fl\mid\fm$. Thus, if $(\fm,\fp)=1$, then by Lemma \ref{lem} and the Chinese Remainder Theorem, $\fp$ splits completely in $K(\phi[\fm])$ if and only if $\phi(\F_\fp)$ contains an $A$-submodule isomorphic to $(A/\fm A)^r$.
	
	Using this, one obtains the following expression for $\pi_\phi(x;\ft,\fa)$.
	
	\begin{lemma}\label{lem:mobius_sum}
		Let $\ft \in A$ be a monic polynomial and $\fa \in (A/(\ft))^*$ be arbitrary. For $x\geq 2$, we have that
		\begin{align}\label{eq:finite-mobius-sum}
			\pi_\phi(x;\ft,\fa)
			=
			\summonics_{\deg \fm\leq x/r}
			\mu_q(\fm)\pi_{\phi,\fm}(x;\ft,\fa),
		\end{align}
		where
		\begin{align}
			\pi_{\phi,\fm}(x;\ft,\fa)
			\coloneqq
			\#\{\fp\in\mathcal{P}_\phi:
			\deg\fp=x,\,
			\fp\equiv\fa \pmod \ft,\,
			\fp \text{ splits completely in }K(\phi[\fm])\}.\label{def0:pi-phi-m}
		\end{align}
	\end{lemma}
	
	\begin{proof}
		By \eqref{def:pi-phi-x,t,a} and Lemma \ref{lem},
		\begin{align}
			\pi_{\phi}(x;\ft,\fa)
			=
			\#\{&\fp\in\mathcal{P}_\phi:
			\deg\fp=x,\,
			\fp\equiv \fa \pmod \ft,\\ &\fp \text{ does not split completely in }
			K(\phi[\fl])
			\text{ for all primes }\fl\neq\fp\}.
		\end{align}
		If $\fp$ splits completely in $K(\phi[\fl])$, then by Lemma \ref{lem}, $\phi(\F_\fp)$ contains an $A$-submodule isomorphic to $(A/\fl A)^r$. Since $\#\phi(\F_\fp)=\#\F_\fp=q^x$ and $\#(A/\fl A)^r=q^{r\deg\fl}$, it then follows that
		$q^{r\deg\fl}\leq q^x$,
		and hence $\deg\fl\leq x/r$. In particular, since $r\geq 2$ and $\deg\fp=x$, we must have $\fl\neq\fp$. Hence,
		\begin{align}
			\pi_{\phi}(x;\ft,\fa)
			=
			\#\{&\fp\in\mathcal{P}_\phi:
			\deg\fp=x,\,
			\fp\equiv \fa \pmod \ft,\\&\fp \text{ does not split completely in }
			K(\phi[\fl])\text{ for all primes }\fl
			\text{ with }\deg\fl\leq x/r\}.
		\end{align}
		For monics $\fm \in A$, define the set
		\begin{align}
			S_{\phi,\fm}(x;\ft,\fa)
			\coloneqq
			\{\fp\in\mathcal{P}_\phi:
			\deg\fp=x,\,
			\fp\equiv\fa \pmod \ft,\,
			\fp \text{ splits completely in }K(\phi[\fm])\}.
		\end{align}
		Also define
$
			X\coloneqq
			\{\fp\in\mathcal{P}_\phi:
			\deg\fp=x,\,
			\fp\equiv\fa \pmod \ft\}
			=
			S_{\phi,1}(x;\ft,\fa).
$
		As a result, we have 	
		\begin{align}
			\pi_{\phi}(x;\ft,\fa)
			=
			\#\bigcap_{\substack{\text{primes }\fl\\
					\deg\fl\leq x/r}}
			\left(
			X\setminus S_{\phi,\fl}(x;\ft,\fa)
			\right)
			=\#
			\Big(
			X\setminus
			\bigcup_{\substack{\text{primes }\fl\\
					\deg\fl\leq x/r}}
			S_{\phi,\fl}(x;\ft,\fa)
			\Big).
		\end{align}
		The set of primes $\fl$ with $\deg\fl\leq x/r$ is finite. Fix a total ordering of the prime polynomials of $A$. Hence, by the inclusion-exclusion principle, 
		\begin{align}
			\pi_{\phi}(x;\ft,\fa)
			&=
			\#X
			-
			\sum_{k\geq 1}(-1)^{k+1}
			\sum_{\substack{
					\fl_1<\fl_2<\cdots<\fl_k\\
					\deg\fl_i\leq x/r}}
			\#
			\Big(
			\bigcap_{i=1}^k
			S_{\phi,\fl_i}(x;\ft,\fa)
			\Big).
		\end{align}
		By the discussion preceding this lemma, this simplifies to
		\begin{align}\label{eq:inclusion-exclusion}
			\pi_{\phi}(x;\ft,\fa)
			=
			\# X
			-
			\sum_{k=1}^\infty (-1)^{k+1}
			\sum_{\substack{
					\fl_1<\fl_2<\cdots<\fl_k\\
					\deg\fl_i\leq x/r}}
			\#
			S_{\phi,\fl_1\fl_2\cdots \fl_k}(x;\ft,\fa).
		\end{align}
		Let $\fm=\fl_1\cdots\fl_k$ be a product occurring in \eqref{eq:inclusion-exclusion}. If $S_{\phi,\fm}(x;\ft,\fa)$ is nonempty and $\fp\in S_{\phi,\fm}(x;\ft,\fa)$, then
		$
		\deg\fl_i\leq x/r<x=\deg\fp
		$
		for every $i$. Hence $(\fm,\fp)=1$. Then, by the discussion preceding the lemma, $\phi(\F_\fp)$ contains an $A$-submodule isomorphic to $(A/\fm A)^r$. Thus
$
			q^{r\deg\fm}
			=
			\#(A/\fm A)^r
			\leq
			\#\phi(\F_\fp)
			=
			q^x,
$
		and hence $\deg\fm\leq x/r$. Consequently, \eqref{eq:inclusion-exclusion} is equivalent to
		\begin{align}
			\pi_{\phi}(x;\ft,\fa)
			&=
			\#X
			-
			\sum_{k\geq 1}(-1)^{k+1}
			\sum_{\substack{
					\fm\text{ square-free }\\\fm\text{ has }k\text{ distinct prime divisors}\\
					\deg\fm\leq x/r}}
			\#
			S_{\phi,\fm}(x;\ft,\fa),
		\end{align}
		or
$
			\pi_{\phi}(x;\ft,\fa)
			=
			\sum'_{\deg\fm\leq x/r}
			\mu_q(\fm)
			\#S_{\phi,\fm}(x;\ft,\fa).
$
		Since $\#S_{\phi,\fm}(x;\ft,\fa)=\pi_{\phi,\fm}(x;\ft,\fa)$, this completes the proof.
	\end{proof}
	
	We now recall the definition of the Artin symbol [cf. \cite{Rosen}, pp. 121-122]. Let $L/K$ be a Galois extension and $O_L$ denote the ring of integers of $L$. For a prime $\fp \in K$ unramified in $L$ and a prime $\mathfrak{P}$ in $L$ with $\mathfrak{P}\mid\fp$, there exists a unique element $\sigma_\mathfrak{P} \in \Gal(L/K)$ such that for every $\alpha \in \mathcal{O}_L$,
	\begin{align}\label{def:artin-symbol}
		\sigma_\mathfrak{P}(\alpha)=\alpha^{q^{\deg \fp}} \pmod {\mathfrak{P}}.
	\end{align}
	This map is called the Frobenius automorphism at $\mathfrak{P}$. As $\mathfrak{P}$ runs over all the primes dividing $\fp$ in $L$, the corresponding Frobenius automorphisms form a conjugacy class in $\Gal(L/K)$, called the Artin symbol of $\fp$ in $L/K$ and denoted by $(\fp,L/K)$. When $\Gal(L/K)$ is abelian, this conjugacy class has only one element, which we identify with that element.
	
	Let $\mathbbm{1}$ denote the conjugacy class of the identity element in the Galois group $\Gal(L/K)$. A prime $\fp\in K$ splits completely in $L/K$ if and only if $\fp$ is unramified in $L/K$ and the Artin symbol
	$(\fp,L/K)=\mathbbm{1}$.
	Therefore, the definition \eqref{def0:pi-phi-m} of $\pi_{\phi,\fm}(x;\ft,\fa)$ can be rewritten as
	\begin{align}\label{def:pi-phi-m}
		\pi_{\phi,\fm}(x;\ft,\fa)
		=
		\#\{\fp\in\mathcal{P}_\phi:
		\deg\fp=x,\,
		\fp\equiv\fa \pmod \ft,\,
		(\fp,K(\phi[\fm])/K)=\mathbbm{1}\}.
	\end{align}
	
	We now translate the condition $\fp\equiv\fa \pmod \ft$ in the language of Artin symbols, using the Carlitz module $\Carl$. Let $\lambda_\ft$ generate the $A/\ft A$-module $\Carl[\ft]$. From Theorem 12.10 of \cite{Rosen}, we see that for primes $\fp$ coprime to $\ft$, the Artin symbol $(\fp,K(\Carl[\ft])/K)$ is the automorphism which takes $\lambda_\ft$ to $\Carl_\fp(\lambda_\ft)$. For primes counted by $\pi_{\phi,\fm}(x;\ft,\fa)$, the condition
	$\fp\equiv\fa\pmod\ft$, with $\fa\in(A/\ft A)^*$, implies
	$\fp\nmid\ft$.
	\begin{lemma}\label{lem:AP-psi}
		Let $\ft \in A$ be a monic polynomial and $\fa \in A$ be such that $\fa \pmod\ft$ is invertible in  $A/(\ft)$. Let $\fp\in A$ be a prime so that $\fp \nmid \ft$. Then the prime $\fp \equiv \fa \pmod\ft$ if and only if Artin symbol $(\fp,K(\Carl[\ft])/K)$ acts as $\Carl_\fa$ on $\lambda_\ft$.
	\end{lemma}
	\begin{proof}
		Since $\fa \pmod \ft \in (A/(\ft))^*$, we know that $(\fp,K(\Carl[\ft])/K)$ acts as $\Carl_\fp$ on $\lambda_\ft$. Now, by Theorem 2.3 of \cite{Hayes}, there exists an isomorphism 
		\begin{align}
			\Phi : (A/\ft A)^* &\to \Gal(K(\Carl[\ft])/K)\\
			\fm &\mapsto \xi_\fm,
		\end{align}
		where $\xi_\fm(\lambda_\ft)=\Carl_\fm(\lambda_\ft)$. Hence $\fp \equiv \fa \pmod \ft$ is equivalent to $\Carl_\fp(\lambda_\ft) = \Carl_\fa(\lambda_\ft)$.
	\end{proof}
	
	As a result, we can modify and write \eqref{def:pi-phi-m} as
	\begin{align}
		\pi_{\phi,\fm}(x;\ft,\fa)=\#\{\fp\in\mathcal{P}_\phi:\deg\fp=x,(\fp,K(\Carl[\ft])/K) (\lambda_\ft)=\Carl_\fa (\lambda_\ft) \text{ and } (\fp,K(\phi[\fm])/K)=\mathbbm{1}\}.
	\end{align}
	From Proposition 9.11 of \cite{Rosen}, we see that the Artin symbol of the compositum field 
	\begin{align}\label{def:sigma_p}
		\sigma_\fp = (\fp,K(\phi[\fm])K(\Carl[\ft])/K)
	\end{align}
	captures the behaviour of the individual Artin maps when restricted to subfields as follows
	\begin{align}
		\sigma_\fp|_{K(\phi[\fm])} &= (\fp,K(\phi[\fm])/K), \text{ and } \; 
		\sigma_\fp|_{K(\Carl[\ft])} = (\fp,K(\Carl[\ft])/K).
	\end{align}
	Let $\nu_{\fm,\ft}$ be the unique automorphism on $K(\phi[\fm])K(\Carl[\ft])$ which fixes $K(\phi[\fm])$ and acts as $\Carl_\fa$ on $K(\Carl[\ft])$, provided it exists. Let $\mathcal{C}_{\fm,\ft}$ denote the conjugacy class of $\nu_{\fm,\ft}$ in $G=\Gal(L_{\fm,\ft}/K)$.
	This implies that $\pi_{\phi,\fm}(x;\ft,\fa)$ can be rewritten as
	\begin{align}
		\pi_{\phi,\fm}(x;\ft,\fa)=
		\begin{cases}
			\#\{\fp\in\mathcal{P}_\phi:\deg\fp=x,\sigma_\fp = \mathcal{C}_{\fm,\ft}\} \text{ if } \nu_{\fm,\ft} \text{ exists},\\
			0 \text{ otherwise},
		\end{cases}
	\end{align}
	Such a $\nu_{\fm,\ft}$ exists only if the intersection field $K(\phi[\fm]) \cap K(\Carl[\ft])$ is fixed by $\Carl_\fa$. Thus, by the definition of $\gamma_{\fa,\ft}$
	in \eqref{def:gamma-intro}, we may write
	\begin{align}\label{eq:for-CDT}
		\pi_{\phi,\fm}(x;\ft,\fa)=
		\begin{cases}
			\#\{\fp\in\mathcal{P}_\phi:\deg\fp=x,\sigma_\fp = \mathcal{C}_{\fm,\ft}\} &\text{if }\gamma_{\fa,\ft}(K(\phi[\fm]))=1,\\
			0 &\text{otherwise}.
		\end{cases}
	\end{align}
	\begin{lemma}\label{conjugacy-singleton}
		Let $\mathcal{C}_{\fm,\ft}$ denote the conjugacy class of $\nu_{\fm,\ft}$ in $G=\Gal(L_{\fm,\ft}/K)$. Then $\mathcal{C}_{\fm,\ft}=\{\nu_{\fm,\ft}\}$.
	\end{lemma}
	\begin{proof}
		Let $\nu_{\fm,\ft}'$ be a conjugate of $\nu_{\fm,\ft}$ in $G$. Since $\nu_{\fm,\ft}$ acts as identity on $K(\phi[\fm])$, so will $\nu_{\fm,\ft}'$. Moreover, as $K(\Carl[\ft])/K$ is an abelian extension, $\nu_{\fm,\ft}$ restricted to $K(\Carl[\ft])$ has only one conjugate automorphism. Thus,  $\nu_{\fm,\ft}'(\lambda_\ft)=\nu_{\fm,\ft}(\lambda_\ft)=\Carl_\fa(\lambda_\ft)$. Hence, $\nu_{\fm,\ft}'=\nu_{\fm,\ft}$ on $L_{\fm,\ft}$.
	\end{proof}	
	This motivates invoking the Chebotarev density theorem over function fields.
	For a finite Galois extension $L/K$, recall that $\F_L=L\cap\overline{\F_q}$, and put $r_L=[\F_L:\F_q]$. Also, recall that  $\tau:c\mapsto c^q$ denotes the Frobenius automorphism. For $x \in \mathbb{N}$, we define
	\begin{equation}
		\pi_\mathcal{C}(x,L/K)\coloneqq\#\{\fp\in A:\deg\fp=x, \fp \text{ unramified in } L/K \text{ and }(\fp,L/K) \subseteq \mathcal{C}\}.
	\end{equation}
	\begin{theorem} [\cite{Fried-Jarden}, Prop. 6.4.8, \cite{Kuo-Liu2009}, Theorem 7] \label{CDT}
		Recall that $K=\F_q(T)$.  Let $L/K$ be a finite Galois extension with Galois group $G$ and $\mathcal{C}\subset G$ be a conjugacy class whose restriction to $\F_L$ is equal to $\tau^n$.
		Then, for $x \not\equiv n (\text{mod }r_L)$, we have $\pi_\mathcal{C}(x,L/K) = 0$ and for $x \equiv n (\text{mod }r_L)$, we have
		\begin{equation}
			\left|\pi_\mathcal{C}(x,L/K)-r_L\frac{|\mathcal{C}|}{|G|}\frac{q^x}{x}\right| \leq  \frac{2|\mathcal{C}|}{x|G|}\big((|G|+g_Lr_L)q^{x/2}+|G|(2g_K+1)q^{x/4}+g_Lr_L+|G|\big),
		\end{equation}
		where $g_L$ and $g_K$ are the genus of $L$ and $K$, respectively.
	\end{theorem}
	
	Let $\pi(x)$ denote the number of primes in $K$ of degree $x$. Substituting $L=K$ in Theorem \ref{CDT} yields the prime number theorem for polynomials:
	\begin{align}\label{eq:PNT}
		\pi(x)={q^x}/{x}+O\big({q^{x/2}}/{x}\big).
	\end{align}
	
	We will apply the Chebotarev Density Theorem on the compositum fields
	$L_{\fm,\ft}$,
	for monics $\fm\in A$ (see \eqref{def:Lmt-intro}). Let $\mathcal{C}_{\fm,\ft}$ be as defined in Lemma \ref{conjugacy-singleton} and $\epsilon_{\fm,\ft}$ be a non-negative integer such that the restriction of $\nu_{\fm,\ft}$ to the field $\F_{L_{\fm,\ft}}$ is equal to the map
	$\tau^{\epsilon_{\fm,\ft}}$.
	Using Lemma \ref{conjugacy-singleton}, \eqref{def:gamma-intro}, \eqref{eq:PNT}, and the fact that there are only finitely many primes $\fp$ of bad reduction, we have the following corollary of Theorem \ref{CDT}. Henceforth, we let 
	\begin{align}\label{def:rm-nm-gm}
		r_{\fm,\ft} = [L_{\fm,\ft}\cap\overline{\F_q}:\F_q], \quad
		n_{\fm,\ft} = [L_{\fm,\ft} : K], \; \text{ and}
		\;\;
		g_{\fm,\ft} = \text{ the genus of }L_{\fm,\ft}.
	\end{align}
	\begin{corollary}\label{for-main}
		For $x \in \mathbb{N}$, if $\gamma_{\fa,\ft}(K(\phi[\fm])) = 1$ and $x \equiv \epsilon_{\fm,\ft} \pmod {r_{\fm,\ft}}$, then we have
		\begin{align}
			\pi_{\phi,\fm}(x;\ft,\fa) = 
			\frac{r_{\fm,\ft}}{n_{\fm,\ft}}\pi(x) + O\left(\left(1+\frac{r_{\fm,\ft} g_{\fm,\ft}}{n_{\fm,\ft}}\right)\frac{q^{x/2}}{x}\right).
		\end{align}
		Otherwise $\pi_{\phi,\fm}(x;\ft,\fa) = 0$.
	\end{corollary}

	\subsection{Some useful bounds}
	
	We now give some useful bounds below. 
	
	\begin{theorem}\cite[Theorem 3.14]{cojo-shulman-II}\label{lower-bound-K(phi[m]):K}
		Let $\phi$ be a Drinfeld $A$-module over $K$ of rank $r\geq 2$ and define $\rho\coloneqq\text{rank}_A\End_{\overline{K}}(\phi)$. Then, for any $\fm \in A \setminus \F_q$, we have
		\begin{align}
			[K(\phi[\fm]):K]\gg_{\phi,K} \frac{q^{\frac{r^2\deg\fm}{\rho}}}{\log \rho + \log \deg \fm + \log \log q}.
		\end{align} 
	\end{theorem}
	
	Let $K_\infty=\F_q\big(\big(\frac{1}{T}\big)\big)$ be the completion of $K$ at the prime at infinity. Then, $K_\infty$ is the quotient field of the power series ring
	\begin{align}
		\F_q\Big[\Big[\frac{1}{T}\Big]\Big] = \Big\{\sum_{i=0}^{\infty}b_iT^{-i} : b_i \in \F_q\Big\}.
	\end{align}
	
	\begin{lemma}\label{geometric-finite-extension}
		If $H/K_\infty$ is a finite field extension, then $H\cap \overline{\F_q}$ is a finite extension of $\F_q$. In particular, we have
		$[H \cap \overline{\F_q}:\F_q]\leq [H:K_\infty]$.
	\end{lemma}
	\begin{proof}
		We proceed by contradiction. Let $H/K_\infty$ be a finite field extension of degree $d$. Suppose $H\cap \overline{\F_q}$ is an extension of $\F_q$ of degree $\geq d+1$. Then, there exists a set $\{\alpha_i\}_{i=1}^{d+1} \subseteq H\cap \overline{\F_q}$, which is linearly independent over $\F_q$. But since $H/K_\infty$ is a finite extension, there exist  $c_1, c_2, \cdots, c_{d+1} \in K_\infty$, not all zero, such that $
			c_1 \alpha_1 + c_2 \alpha_2 + \cdots + c_{d+1} \alpha_{d+1} = 0$.
		Clearing denominators, we obtain $f_i = \sum_{j=0}^\infty b_{ij}T^{-j} \in \F_q\big[\big[\frac{1}{T}\big]\big]$, for $i=1,2,\cdots,d+1$, not all zero, such that
		\begin{align}\label{not-all-zero}
			f_1 \alpha_1 + f_2 \alpha_2 + \cdots + f_{d+1} \alpha_{d+1} = 0.
		\end{align}
		Collecting coefficients of the same degree , we see that
		\begin{align}
			\sum_{j=0}^\infty (b_{1j}\alpha_1+b_{2j}\alpha_2+\cdots+b_{d+1j}\alpha_{d+1})T^{-j}=0.
		\end{align}
		Hence, for every $j \geq 0$, we have that $b_{1j}\alpha_1+b_{2j}\alpha_2+\cdots+b_{d+1j}\alpha_{d+1} = 0$, where each $b_{ij} \in \F_q$. By linear independence of $\alpha_i's$ over $\F_q$, we must then have that $b_{ij}=0$, for all $i,j$. Thus, each $f_i=0$, which is a contradiction since $f_i$'s were not all zero in \eqref{not-all-zero}.
		
		Hence, if $H/K_\infty$ is a finite field extension, then $H\cap \overline{\F_q}$ is a finite extension of $\F_q$.
	\end{proof}
	
	Let $C_\infty$ denote the completion of $\overline{K_\infty}$ at the prime at infinity. For any Drinfeld $A$-module $\phi$ over $K$, let $\phi_{\text{tors}}$ be defined as
	\begin{align}\label{phi-tors}
		\phi_{\text{tors}}\coloneqq\{\alpha\in C_\infty:\phi_\fa(\alpha)=0 \text{ for some non-zero }\fa\in A\}.
	\end{align}
	Note that $\phi_\text{tors}$ contains the torsion points $\phi[\fm]$, for every $\fm \in A$.
	
	\begin{lemma}\label{E-phi}
		Let $L_{\fm,\ft}$ be as in \eqref{def:Lmt-intro}. There exists a positive integer $E(\phi,K)$ (depending only on $\phi$ and $K$) such that for every $\fm,\ft \in A$, we have $
			r_{\fm,\ft}=[L_{\fm,\ft}\cap \overline{\F_q}:\F_q] < E(\phi,K)$.
	\end{lemma}
	\begin{proof}
		Let $\phi_{\text{tors}}$ and $\Carl_{\text{tors}}$ be as defined in \eqref{phi-tors}. We have $\phi:A\to K\{\tau\}\to K_\infty\{\tau\}$. Then viewing $\phi$ as a Drinfeld $A$-module over $K_\infty$, \cite[Theorem 4.6.9]{Goss} allows us to take $M=K_\infty$ in \cite[Remark 4.3.6.2]{Goss}. Hence, there exist finite field extensions $K_{\phi}/K_\infty$ and $K_C/K_\infty$ such that $\phi_{\text{tors}}\subseteq K_{\phi}$ and $\Carl_{\text{tors}}\subseteq K_C$.
		Hence, we have the inclusion of compositums
		\begin{align}\label{inclusion}
			\prod_{\fm, \ft \in A} K(\phi[\fm])K(\Carl[\ft]) \subseteq K_\phi K_C.
		\end{align}
		Since both $K_{\phi}/K_\infty$ and $K_C/K_\infty$ are finite extensions, we have that $K_{\phi}K_C/K_\infty$ is also a finite extension and hence by Lemma \ref{geometric-finite-extension}, $K_{\phi}K_C \cap \overline{\F_q}$ is a finite extension of $\F_q$.
		Thus, by \eqref{inclusion}, there exists a positive integer $E(\phi,K)$ such that for every $\fm,\ft \in A$, 
		$[L_{\fm,\ft}\cap \overline{\F_q}:\F_q] < E(\phi,K)$
	\end{proof}
	
	We observe that the genus $g_L$ of a field $L/K$, as denoted in Theorem \ref{CDT}, is essentially defined over the full constant field $\mathbb{F}_L$ of $L$. So, we have that $g_L$ is the genus of the field $L/\mathbb{F}_L$.
	To specify the role of the base constant field, we introduce the notation $g_{L/F}$ for any constant field $F \subseteq \mathbb{F}_L$ to denote the genus of the field $L/F$.
	Hence, we have that $g_L = g_{L/\mathbb{F}_L}$. 
	
	\begin{lemma}\label{genus-support}
		For every $\fm,\ft \in A$, we have that
		\begin{align}\label{eq:genus-support}
			g_{\fm,\ft} 
			&\leq [L_{\fm,\ft} : K(\phi[\fm])]\cdot g_{K(\phi[\fm])}+[L_{\fm,\ft} : K(\Carl[\ft])]\cdot g_{K(\Carl[\ft])}\\
			&+([L_{\fm,\ft} : K(\phi[\fm])] - 1)([L_{\fm,\ft} : K(\Carl[\ft])]-1).
		\end{align}
	\end{lemma}
	\begin{proof}
		Let $F = \F_{L_{\fm,\ft}}$. We define the compositum fields $L_1 \coloneqq K(\phi[\fm])F$ and $L_2 \coloneqq K(\Carl[\ft])F$.
		Then, we note that $L_{\fm,\ft}=L_1 L_2$ and $\F_{L_1}=\F_{L_2}=F$. Hence, Castelnuovo's inequality \cite[Theorem 3.11.3]{Stichtenoth}
		implies that 
		\begin{align}\label{castel-appl}
			g_{\fm,\ft} \leq [L_{\fm,\ft} : L_1]\cdot g_{L_1}+[L_{\fm,\ft} : L_2]\cdot g_{L_2}+([L_{\fm,\ft} : L_1]-1)([L_{\fm,\ft} : L_2]-1).
		\end{align}
		By \cite[Theorem 3.6.3(b)]{Stichtenoth}, we have
		\begin{align}\label{genus-equality}
			g_{L_1}=g_{L_1/F}=g_{K(\phi[\fm])} \quad \text{and} \quad g_{L_2}=g_{L_2/F}=g_{K(\Carl[\ft])}.
		\end{align}
		We also have that
		\begin{align}\label{degree-inequalities}
			[L_{\fm,\ft} : L_1] \leq [L_{\fm,\ft} : K(\phi[\fm])] \quad \text{and} \quad [L_{\fm,\ft} : L_2] \leq [L_{\fm,\ft} : K(\Carl[\ft])].
		\end{align}
		Hence, using \eqref{genus-equality} and \eqref{degree-inequalities} in \eqref{castel-appl}, we obtain
		\eqref{eq:genus-support}.
	\end{proof}

	We will need the following inequality as stated in \cite[Proposition 3.9]{cojo-shulman-II}.
	
	\begin{proposition}\cite[Corollary 7]{Gardeyn}\label{pre-D-phi}
		Let $\psi$ be a Drinfeld $A$-module over $K$. Let $g_{K(\psi[\fm])}$ denote the genus of the field $K(\psi[\fm])$. There exists a constant $G(\psi,K)\in\mathbb{N}$, depending on $K$ and $\psi$, such that for any $\fm \in A\setminus\F_q$, we have $g_{K(\psi[\fm])} \leq G(\psi,K)\cdot[K(\psi[\fm]):K]\cdot\deg\fm$.
	\end{proposition}

	This allows us to deduce
	
	\begin{proposition}\label{D-phi}
		There exist constants $D(\Carl,K), D(\phi,K) \in \mathbb{N}$ such that for every $\fm \in A\setminus\F_q$, 
		\begin{align}
			\frac{g_{\fm,\ft}}{n_{\fm,\ft}} \ll D(\phi,K)\deg\fm + D(\Carl,K)\deg \ft.
		\end{align}
	\end{proposition}
	\begin{proof}
		By Proposition \ref{pre-D-phi}, there exist constants $D(\Carl,K) \in \mathbb{N}$ and $D(\phi,K) \in \mathbb{N}$ satisfying
		\begin{align}
			g_{K(\Carl[\ft])} &\leq D(\Carl,K)\cdot[K(\Carl[\ft]):K]\cdot\deg\ft,\\
			g_{K(\phi[\fm])} &\leq D(\phi,K)\cdot[K(\phi[\fm]):K]\cdot\deg\fm.\label{bounds-genus}
		\end{align}
		By Lemma \ref{genus-support}, we have the inequality	
		\begin{align}\label{genus-before}
			g_{\fm,\ft} 
			&< [L_{\fm,\ft} : K(\phi[\fm])]\cdot g_{K(\phi[\fm])}+[L_{\fm,\ft} : K(\Carl[\ft])]\cdot g_{K(\Carl[\ft])}+[L_{\fm,\ft} : K(\phi[\fm])][L_{\fm,\ft} : K(\Carl[\ft])]\\
			&= [L_{\fm,\ft} : K]\left(\frac{g_{K(\phi[\fm])}}{[K(\phi[\fm]):K]}+\frac{g_{K(\Carl[\ft])}}{[K(\Carl[\ft]):K]}+\frac{[L_{\fm,\ft}:K]}{[K(\phi[\fm]):K][K(\Carl[\ft]):K]}\right).
		\end{align}
		Hence, by \eqref{bounds-genus}, we obtain
		\begin{align}
			\frac{g_{\fm,\ft}}{n_{\fm,\ft}}&< \frac{1}{[K(\phi[\fm]) : K]}g_{K(\phi[\fm])} + \frac{1}{[ K(\Carl[\ft]):K]}g_{K(\Carl[\ft])} + \frac{1}{[K(\phi[\fm])\cap K(\Carl[\ft]):K]}\\
			&\leq D(\phi,K) \deg \fm + D(\Carl,K) \deg \ft + 1\ll D(\phi,K) \deg \fm + D(\Carl,K) \deg \ft.
		\end{align}
		This proves the proposition.
	\end{proof}

	We note that the estimate in Proposition \ref{D-phi} also holds when $\fm=1$. 	Indeed, if $\ft=1$, then $L_{1,1}=K$, and hence $g_{1,1}=0$. 	If $\deg \ft\geq 1$, then $L_{1,\ft}=K(\Carl[\ft])$, and by Proposition \ref{pre-D-phi}, 	
	\begin{align} 		\frac{g_{1,\ft}}{n_{1,\ft}} 		= 		\frac{g_{K(\Carl[\ft])}} 		{[K(\Carl[\ft]):K]} \ll D(\Carl,K)\deg \ft. 	\end{align} 	
	Thus, the estimate in Proposition \ref{D-phi} holds for every monic $\fm\in A$.

	We apply Lemma \ref{E-phi} and Proposition \ref{D-phi} to Corollary \ref{for-main}. When either $\fm\neq1$ or $\ft\neq 1$, the initial $1$ in the error term of Corollary \ref{for-main} is absorbed by the bounds from Lemma \ref{E-phi} and Proposition \ref{D-phi}. Then, for $x \in \mathbb{N}$, if $\gamma_{\fa,\ft}(K(\phi[\fm])) = 1$ (see \eqref{def:gamma-intro}) and $x \equiv \epsilon_{\fm,\ft} \pmod {r_{\fm,\ft}}$, we have
	\begin{align}\label{eq:pi-phi-m}
		\pi_{\phi,\fm}(x;\ft,\fa)
		=
		\frac{r_{\fm,\ft}}{n_{\fm,\ft}}\pi(x)
		+ O\left(E(\phi,K)D(\phi,K) \deg \fm\frac{q^{x/2}}{x}\right)
		+ O\left(E(\phi,K)D(\Carl,K) \deg \ft\frac{q^{x/2}}{x}\right).
	\end{align}
	Otherwise $\pi_{\phi,\fm}(x;\ft,\fa) = 0$. When $\fm=\ft=1$, we have $L_{1,1}=K$, and hence by Corollary \ref{for-main}, we have
	an additional $O_{\phi,K}(q^{x/2}/x)$, which is absorbed into the first error term below. Hence, the M\"{o}bius sum in \eqref{eq:finite-mobius-sum} yields 
	\begin{align}
		\pi_\phi(x;\ft,\fa)
		& =\summonics_{\substack{\deg \fm \leq x/r}} \frac{\mu_q(\fm) \gamma_{\fa,\ft}(K(\phi[\fm])) r_{\fm,\ft}(x)}{n_{\fm,\ft}}\pi(x)+ O\bigg(\frac{q^{x/2}}{x} E(\phi,K)D(\phi,K)
		\summonics_{\deg \fm \leq x/r} \deg \fm \bigg)
		\\
		&  + O\left(\frac{q^{x/2+x/r}}{x} E(\phi,K)D(\Carl,K)\deg \ft\right),
		\label{eq:naive}
	\end{align}
	where $r_{\fm,\ft}(x)$ is equal to $r_{\fm,\ft}$ whenever $x\equiv\epsilon_{\fm,\ft}\pmod{r_{\fm,\ft}}$, and $0$ otherwise.

	\section{Proof of Theorem \ref{thm:estimation}}\label{sec:estimation-proof}

	We wish to estimate the quantity \eqref{def:pi-phi-x,t,a}.
	By \eqref{eq:finite-mobius-sum}, we have that 
	\begin{align} \label{eq:finite-mobius-sum-2}
		\pi_\phi(x;\ft,\fa)=\summonics_{\deg\fm\leq x/r} \mu_q(\fm)\pi_{\phi,\fm}(x;\ft,\fa).
	\end{align}
	One can observe that in the expression \eqref{eq:naive} for $\pi_\phi(x;\ft,\fa)$, when $r=2$, the error terms are not $o(q^x/x)$, and are thus as large as $\pi(x)$ (see \eqref{eq:PNT}). Hence we will split the sum \eqref{eq:finite-mobius-sum-2} into two parts and use \eqref{eq:pi-phi-m} only for one of them. We will let $y=y(x)$ be a parameter to be chosen later, satisfying $y(x) \to \infty$ as $x \to \infty$. Then
	\begin{align}\label{eq:main+error}
		\pi_\phi(x;\ft,\fa)=\main+\error,
	\end{align}
	where, 
	\begin{align}\label{def:main-and-error}
		\main \coloneqq\summonics_{\deg\fm\leq y} \mu_q(\fm)\pi_{\phi,\fm}(x;\ft,\fa),
		\quad
		\error \coloneqq\summonics_{y<\deg\fm\leq x/r} \mu_q(\fm)\pi_{\phi,\fm}(x;\ft,\fa).
	\end{align}

	\subsection{Main term}\label{subsec:CM-main}
	Using the asymptotic formula \eqref{eq:pi-phi-m} for $\pi_{\phi,\fm}(x;\ft,\fa)$ in \eqref{def:main-and-error}, we obtain 
	\begin{align}
		\main
		=&\summonics_{\substack{\deg \fm \leq y}} \frac{\mu_q(\fm) \gamma_{\fa,\ft}(K(\phi[\fm])) r_{\fm,\ft}(x)}{n_{\fm,\ft}}\pi(x)
		+ O\Bigg(\frac{q^{x/2}}{x} E(\phi,K)D(\phi,K)\summonics_{\deg \fm \leq y} \deg \fm \Bigg) \\&+ O\Bigg(\frac{q^{\frac{x}{2}+y}}{x}  E(\phi,K)D(\Carl,K)\deg \ft\Bigg),\label{eq:main-without-tail}
	\end{align}
	where $r_{\fm,\ft}(x)$ is as defined in \eqref{def:r-m,t-x}.
	Since $\sum_{\deg\fm\leq y} \deg \fm \ll y q^y $, the first error term above is  $
		\ll_{\phi,K}q^{\frac{x}{2}+y}\frac{y}{x}$.
	The first and second error terms thus are
	\begin{align}\label{eq:first-main-error}
		\ll_{\phi,K} q^{\frac{x}{2}+y}\Big(\frac{y}{x}+\frac{\deg \ft}{x}\Big).
	\end{align} 
	The first sum in \eqref{eq:main-without-tail} can be extended to run over all monics $\fm \in A$ satisfying $x \equiv \epsilon_{\fm,\ft}\pmod {r_{\fm,\ft}}$, up to an error $E_\text{tail}$ which is
	\begin{align}\label{Etail-pre}
		\ll \pi(x)\summonics_{\substack{\deg \fm> y }} \frac{|\mu_q(\fm)|r_{\fm,\ft}(x)\gamma_{\fa,\ft}(K(\phi[\fm]))}{n_{\fm,\ft}}.
	\end{align}
	Recalling that $n_{\fm,\ft} = [K(\phi[\fm])K(\Carl[\ft]):K]$, we see that
	\begin{align}\label{eq:updated-main-without-tail}
		\main
		= \delta_{\phi,\ft,\fa}(x)\cdot\pi(x)
		&+ E_{\text{tail}} + O_{\phi,K}\Big(q^{\frac{x}{2}+y}\Big(\frac{y}{x}+\frac{\deg \ft}{x}\Big)\Big),
	\end{align}
	where $\delta_{\phi,\ft,\fa}(x)$ is as defined in \eqref{def:delta-intro}, and is shown to be an absolutely convergent sum below.
	
	Using the bound $n_{\fm,\ft} \geq [K(\phi[\fm]):K]$, and using Lemma \ref{E-phi} to bound $r_{\fm,\ft}(x)$, we see that
	\begin{align}\label{eq:E1-bound}
		E_{\text{tail}} \ll E(\phi,K)\pi(x)\summonics_{\substack{\deg \fm > y }}\frac{1}{[K(\phi[\fm]):K]}.
	\end{align}
	Invoking the lower bound for $[K(\phi[\fm]):K]$ given in Theorem \ref{lower-bound-K(phi[m]):K}, we see that \eqref{eq:E1-bound} is 
	\begin{align}
		\ll_{\phi,K} \pi(x) \summonics_{\substack{\deg \fm > y }}\frac{\log \rho + \log \deg \fm + \log \log q}{q^{r^2\deg\fm/\rho}}\ll_{K} \pi(x)\sum_{n>y}q^n\cdot \frac{\log \rho + \log n + \log \log q}{q^{r^2n/\rho}},\label{eq:E1-updated-bound}
	\end{align}
	where $\rho = \text{rank}_A \End_{\overline{K}}(\phi)$. 
	
	Note that $y\to\infty$ as $x\to\infty$, $q\geq 2$, $r\geq 2$ and $\rho \leq r$ (see \cite[Theorem 3.2(ii)]{cojo-shulman-II}). For $n$ sufficiently large, we have 
	$\log \rho + \log n + \log \log q \ll_{\phi, K} \log n$.
	Denote the exponent of $1/q$ in the summand of \eqref{eq:E1-updated-bound} by $na$, where $a=a(r,\rho)=(r^2/\rho) -1$. It is easy to see that $a\geq 1$.
	
	The sum over $n$ in \eqref{eq:E1-updated-bound} is then
	\begin{align}
		\ll_{\phi,K} \sum_{n>y}\frac{\log n}{q^{na/2}}\cdot\frac{1}{q^{na/2}}\ll_{\phi,K} \frac{\log y}{q^{ya/2}}\sum_{n\geq y}\frac{1}{q^{na/2}},
	\end{align}
	since $\frac{\log n}{q^{na/2}}$ is a decreasing function of $n$ as $n \to \infty$. The final sum above is a geometric progression summing to $\ll \frac{1}{q^{ya/2}}$. 
	
	We have thus obtained
	\begin{align}\label{eq:tail-final-bound}
		E_{\text{tail}} \ll_{\phi,K} \pi(x)\cdot \frac{\log y}{q^{ya}} \ll_{\phi,K} \frac{q^{x-\left(\frac{r^2}{\rho}-1\right)y}\log y}{x},
	\end{align}
	using the prime number theorem \eqref{eq:PNT}.

	Moreover, the above discussion also allows us to conclude that the density $\delta_{\phi,\ft,\fa}(x)$ is finite. Indeed, from the bound  $n_{\fm,\ft} \geq [K(\phi[\fm]):K]$, and Lemma \ref{E-phi} for $r_{\fm,\ft}$, we see that the $\delta_{\phi,\ft,\fa}(x)$
	is dominated in absolute value by the series 
	\begin{align}\label{bound-for-density}
		E(\phi,K)\summonics_{\fm\in A}\frac{1}{[K(\phi[\fm]):K]}.
	\end{align}
	However the tail of \eqref{bound-for-density} is simply the sum in  \eqref{eq:E1-bound}, which we showed to be $\ll_{\phi,K}\frac{\log y}{q^{ya}}$
	in the discussion above \eqref{eq:tail-final-bound}, where $a=\frac{r^2}{\rho} -1$. As $q\geq 2$ and $a\geq 1$, the tail converges to $0$ as $y\to\infty$. By Cauchy's criterion, the series \eqref{def:delta-intro} converges absolutely, as required.

\subsubsection{Proving Remark 	\ref{rem:density-bound}}	
	Recall the notation  used in \eqref{def:Lmt-intro}. 
	Let $\fd,\fm\in A$ be monics with $\fd\mid\fm$.
	
	Suppose that
	$\gamma_{\fa,\ft}(K(\phi[\fd]))=1$ and
	$x\equiv\epsilon_{\fd,\ft}\pmod{r_{\fd,\ft}}$.
	Then, there exists a unique automorphism
	$\nu_{\fd,\ft}\in\Gal(L_{\fd,\ft}/K)$
	which fixes $K(\phi[\fd])$, acts as
	$\psi_{\fa}^{\mathrm{Carl}}$ on
	$K(\psi^{\mathrm{Carl}}[\ft])$, and whose
	restriction to the constant field $\F_{L_{\fd,\ft}}$
	of $L_{\fd,\ft}$ is $\tau^x$.
	Since $L_{\fd,\ft}\subseteq L_{\fm,\ft}$ and
	$L_{\fd,\ft}\cap \F_{L_{\fm,\ft}}
	=\F_{L_{\fd,\ft}}$,
	the restrictions of $\nu_{\fd,\ft}$ and
	$\tau^x$ are compatible. Hence, the number
	of automorphisms of $L_{\fm,\ft}/K$ extending
	$\nu_{\fd,\ft}$ and acting as $\tau^x$ on
	$\F_{L_{\fm,\ft}}$ is
	$[L_{\fm,\ft}:L_{\fd,\ft}\F_{L_{\fm,\ft}}]$,
	which is equal to
	\begin{align}
		\frac{[L_{\fm,\ft}:K]
			[\F_{L_{\fd,\ft}}:\mathbb F_q]}
		{[L_{\fd,\ft}:K]
			[\F_{L_{\fm,\ft}}:\mathbb F_q]}
		=
		\frac{n_{\fm,\ft}r_{\fd,\ft}}
		{n_{\fd,\ft}r_{\fm,\ft}}.
	\end{align}

	If either of these conditions fails, there is no
	such automorphism. Therefore, 
	\begin{align}
		\frac{\gamma_{\fa,\ft}(K(\phi[\fd]))
			r_{\fd,\ft}(x)}{n_{\fd,\ft}}
		=
		\frac{r_{\fm,\ft}}{n_{\fm,\ft}}
		\#\Big\{
		\sigma\in\operatorname{Gal}(L_{\fm,\ft}/K):
		\sigma|_{\F_{L_{\fm,\ft}}}=\tau^x,\,
		\sigma|_{K(\psi^{\mathrm{Carl}}[\ft])}
		=\psi_{\fa}^{\mathrm{Carl}},\,
		\sigma|_{K(\phi[\fd])}=1
		\Big\}.
	\end{align}
	Thus, by the inclusion-exclusion principle,
	\begin{align}
		\sum_{\fd\mid\fm}
		\frac{\mu_q(\fd)\gamma_{\fa,\ft}(K(\phi[\fd]))
			r_{\fd,\ft}(x)}{n_{\fd,\ft}}
		=
		\frac{r_{\fm,\ft}}{n_{\fm,\ft}}
		\#\Big\{&
		\sigma\in\operatorname{Gal}(L_{\fm,\ft}/K):
		\sigma|_{\F_{L_{\fm,\ft}}}=\tau^x,
		\sigma|_{K(\psi^{\mathrm{Carl}}[\ft])}
		=\psi_{\fa}^{\mathrm{Carl}},\,\\
		&\sigma|_{K(\phi[\fl])}\neq1
		\text{ for every prime }\fl\mid\fm
		\Big\}\label{IE-upperbound}.
	\end{align}
	Note that the right-hand side in \eqref{IE-upperbound}
	is nonnegative and is bounded above by
	\begin{align}
		\frac{r_{\fm,\ft}}{n_{\fm,\ft}}
		\#\Big\{
		\sigma\in\operatorname{Gal}(L_{\fm,\ft}/K):
		\sigma|_{\F_{L_{\fm,\ft}}}=\tau^x,\,
		\sigma|_{K(\psi^{\mathrm{Carl}}[\ft])}
		=\psi_{\fa}^{\mathrm{Carl}}
		\Big\}
		=
		\frac{r_{\fm,\ft}}{n_{\fm,\ft}}
		\frac{n_{\fm,\ft}}
		{r_{\fm,\ft}\varphi(\ft)}
		=\frac{1}{\varphi(\ft)},
	\end{align}
	since $K(\psi^{\mathrm{Carl}}[\ft]) \cap \overline{\F_q} = \F_q$,
	and hence $K(\psi^{\mathrm{Carl}}[\ft])\cap
	\F_{L_{\fm,\ft}}=\F_q$.
	Consequently, by \eqref{IE-upperbound},
	\begin{align}
		0\leq
		\sum_{\fd\mid\fm}
		\frac{\mu_q(\fd)\gamma_{\fa,\ft}(K(\phi[\fd]))
			r_{\fd,\ft}(x)}{n_{\fd,\ft}}
		\leq\frac{1}{\varphi(\ft)}.
	\end{align}
	Letting $\fm$ run over products of an
	increasing collection of primes in $A$
	whose union contains every prime in $A$,
	the remark follows from the absolute
	convergence of $\delta_{\phi,\ft,\fa}(x)$.

	\normalcolor

	\subsection{Error term}\label{subsec:CM-error}
	
	In order to estimate the term $\sum_{\text{error}}$ in \eqref{eq:main+error}, recall the definition of $\pi_{\phi,\fm}(x;\ft,\fa)$ from \eqref{def:pi-phi-m}. Using Corollary \ref{cor:for-error}, we see that each prime $\fp$ being counted by $\pi_{\phi,\fm}(x;\ft,\fa)$ satisfies
	\begin{align}\label{two-congruences}
		1+\fb_\fp+u_\fp \fp \equiv 0 \pmod {\fm^r} \quad \text{and} \quad 1+(-1)^{r-1}u_\fp \fp \equiv 0 \pmod \fm,
	\end{align}
	for some $u_\fp \in \F_q^*$ and some $\fb_\fp \in A$ with $\deg \fb_\fp \leq \frac{(r-1)x}{r}$. Since $1+\fb_\fp+u_\fp \fp$ and $1+(-1)^{r-1}u_\fp \fp$ are both divisible by $\fm$, we may eliminate $\fp$ from the second congruence in \eqref{two-congruences} to obtain the equivalent system of congruences 
	\begin{align}\label{equivalent-two-congruences}
		1+\fb_\fp+u_\fp \fp \equiv 0 \pmod {\fm^r} \quad \text{and} \quad 1+\fb_\fp+(-1)^r \equiv 0 \pmod \fm.
	\end{align}
	Summing over all possibilities $u_\fp = u$ and $\fb_\fp = \fb$, where $u\in\F_q^*$ and $\fb \in A$ with $\deg \fb \leq \frac{(r-1)x}{r}$, we see that 
	\begin{align}\label{eq:triple-sum}
		\pi_{\phi,\fm}(x;\ft,\fa)\leq \sum_{u\in\F_q^*}
		\sum_{\substack{\fb\in A\\ \deg \fb \leq \frac{(r-1)x}{r}\\ 1+\fb+(-1)^r \equiv 0 \pmod \fm}} \quad
		\sum_{\substack{\deg \fp = x\\ \fp \equiv \fa \pmod \ft\\ 1+\fb+u\fp \equiv 0 \pmod {\fm^r}}} 1.
	\end{align}
	We observe that if the congruences $ \fp \equiv \fa \pmod \ft, \fp \equiv -(1+\fb)u^{-1} \pmod {\fm^r}$ are incompatible modulo $\gcd(\ft,\fm^r)$, then the innermost sum is empty. Otherwise, by the Chinese remainder theorem, they can be written as a single congruence $\fp \equiv \fc \pmod {[\ft,\fm^r]}$ for some nonzero $\fc \in A$ with $\deg \fc < \deg [\ft,\fm^r]$ and where $[\ft,\fm^r]$ is the least common multiple of the monics $\ft$ and $\fm^r$ in $A$.
	
	Since the prime $\fp \in A$ is monic, we may write $\fp$ as $\fc +[\ft,\fm^r]\cdot\fh_1$ for some $\fh_1\in A$ where $\fh_1 =0$ or $\fh_1$ is monic with $\deg \fh_1 = x - \deg [\ft,\fm^r]$. Hence, the innermost sum in \eqref{eq:triple-sum} is at most $ 1+ q^{x-\deg[\ft,\fm^r]}$, so that \eqref{eq:triple-sum} is 
	\begin{align}\label{eq:double-sum}
		\ll \left(1+q^{x-\deg [\ft,\fm^r]}\right)
		\sum_{u\in\F_q^*}
		\sum_{\substack{\fb\in A\\ \deg \fb \leq \frac{(r-1)x}{r}\\ \fb \equiv -1-(-1)^r \pmod \fm}} 1.
	\end{align}
	To estimate the innermost sum above, we may again write $\fb$ as $\fb = -1-(-1)^r+\fm\cdot\fh_2$ for some $\fh_2 \in A$. Since $\deg \fb \leq (r-1)x/r$, we must have $\deg \fh_2 \leq \frac{(r-1)x}{r} - \deg \fm$. This means that the innermost sum in \eqref{eq:double-sum} is $\leq q^{\frac{(r-1)x}{r}-\deg \fm + 1} $. Putting these bounds back into \eqref{eq:triple-sum}, we obtain 
	\begin{align}
		\pi_{\phi,\fm}(x;\ft,\fa)
		&\leq \left(q^{\left(\frac{r-1}{r}\right)x-\deg \fm + 1}+q^{\frac{(2r-1)}{r}x - \deg \fm - \deg [\ft,\fm^r]+1}\right)\sum_{u\in\F_q^*}1\\
		&\ll_K q^{\left(\frac{r-1}{r}\right)x-\deg \fm }\left(1+q^{x - \deg [\ft,\fm^r]}\right).
	\end{align}
	We now come back to the expression \eqref{def:main-and-error} for $\sum_{\text{error}}$ to see that
	\begin{align}
		\error
		&\ll_K q^{\left(\frac{r-1}{r}\right)x} \summonics_{y<\deg \fm \leq x/r}  q^{- \deg \fm} + q^{\left(\frac{2r-1}{r}\right)x} \summonics_{y<\deg \fm \leq x/r} q^{- \deg \fm - \deg [\ft,\fm^r]}\\
		&= q^{\left(\frac{r-1}{r}\right)x} \summonics_{y<\deg \fm \leq x/r}  q^{- \deg \fm} + q^{\left(\frac{2r-1}{r}\right)x - \deg \ft} \summonics_{y<\deg \fm \leq x/r} q^{- (r+1)\deg \fm + \deg (\ft,\fm^r)},\label{eq:penultimate-bound-error}
	\end{align}
	where we have used 
	$\deg [\ft,\fm^r]=\deg \ft + r\deg \fm - \deg (\ft,\fm^r)$,
	and $(\ft,\fm^r)\in A$ denotes the monic polynomial which is the gcd of $\ft$, $\fm^r$. We rewrite the first term of \eqref{eq:penultimate-bound-error} as  
	\begin{align}\label{first-term-bound}
		q^{\left(\frac{r-1}{r}\right)x}\sum_{y < \ell \leq x/r}q^{-\ell} \summonics_{\fm : \deg \fm = \ell} 1 = q^{\left(\frac{r-1}{r}\right)x}\sum_{y < \ell \leq x/r} 1 \leq \frac{x}{r} q^{\left(\frac{r-1}{r}\right)x}.
	\end{align}
	We  use the identity
$\label{eq:sum-totient}
		\sum'_{\substack{\fd \mid \ft,  \fd \mid \fm^r}}\varphi(\fd)=q^{\deg (\ft,\fm^r)},
$
	to rewrite the sum in the second term of \eqref{eq:penultimate-bound-error} as
	\begin{align}\label{second-term-double-sum}
		\summonics_{y<\deg \fm \leq x/r} q^{- (r+1)\deg \fm} \summonics_{\substack{\fd \mid \ft \\ \fd \mid \fm^r}}\varphi(\fd).
	\end{align}
	Interchanging the order of summation, we get that \eqref{second-term-double-sum} is equal to
	\begin{align}\label{second-term-interchanged-sum}
		\summonics_{\fd \mid \ft} \varphi(\fd) \summonics_{\substack{y<\deg \fm \leq x/r \\ \fd \mid \fm^r}}  q^{- (r+1)\deg \fm}.
	\end{align}
	Let $\fd_0\coloneqq \fd_0(\fd)$ be defined such that for $\fd = \fp_1^{\alpha_1} \fp_2^{\alpha_2} \cdots \fp_k^{\alpha_k}$, we have that
	\begin{align}\label{d_0}
		\fd_0(\fp_1^{\alpha_1} \fp_2^{\alpha_2} \cdots \fp_k^{\alpha_k})=\fp_1^{\lceil \alpha_1/r \rceil}\fp_2^{\lceil \alpha_2/r \rceil}\cdots\fp_k^{\lceil \alpha_k/r \rceil}.
	\end{align}
	Then, $\fd \mid \fm^r$ if and only if $\fd_0 \mid \fm$. Writing $\fm = \fd_0 \fh_3$ for some monic $\fh_3 \in A$, we rewrite the innermost sum in \eqref{second-term-interchanged-sum} as 
	\begin{align}\label{second-term-interchanged-sum-bound}
		q^{-(r+1)\deg \fd_0} \summonics_{\substack{\fh_3 : \\ y<\deg \fd_0 + \deg \fh_3 \leq x/r}}q^{-(r+1)\deg \fh_3}.
	\end{align}
	Then \eqref{second-term-interchanged-sum-bound} is bounded above by
	\begin{align}
		q^{-(r+1)\deg \fd_0} \sum_{\ell>y-\deg \fd_0} q^{-(r+1)\ell}\summonics_{\fh_3 : \deg \fh_3=\ell} 1
		&= q^{-(r+1)\deg \fd_0} \sum_{\ell>y-\deg \fd_0} q^{-r\ell}\\
		&\ll_K q^{-\deg \fd_0-ry} \label{second-term-interchanged-sum-final-bound},
	\end{align} 
	since the sum $\sum_{\ell>y-\deg \fd_0} q^{-r\ell}$ is the tail sum of a geometric progression.
	
	Hence, by \eqref{second-term-double-sum}, \eqref{second-term-interchanged-sum},\eqref{second-term-interchanged-sum-bound} and \eqref{second-term-interchanged-sum-final-bound}, the sum in the second term of \eqref{eq:penultimate-bound-error} is bounded above by
	\begin{align}\label{second-term-bound}
		q^{-ry}\summonics_{\fd \mid \ft} \varphi(\fd) q^{-\deg \fd_0}.
	\end{align}
	Recall that
	$J_\phi(\ft)=\summonics_{\fd\mid\ft}
	\varphi(\fd)q^{-\deg\fd_0}$,
	as defined in \eqref{def:J-alpha-intro}.
	Thus, by \eqref{eq:penultimate-bound-error}, \eqref{first-term-bound} and \eqref{second-term-bound}, we have obtained
	\begin{align}\label{eq:error-final-bound}
		\error \ll_K \frac{x}{r}q^{\left(\frac{r-1}{r}\right)x} + q^{\left(\frac{2r-1}{r}\right)x-\deg \ft -ry}J_\phi(\ft).
	\end{align}

	\subsection{Estimation of $J_\phi(\ft)$}\label{subsec:J-phi}
	 Observe that the summands $\varphi(\fd) q^{-\deg \fd_0}$ are multiplicative in $\fd$ and hence $J_\phi(\ft)$ is multiplicative in $\ft$. 
	 
	 Since $\fd\mid\fd_0^r$ by definition, we have $\deg \fd_0 \geq \frac{1}{r}\deg\fd$. Hence each summand of $J_\phi(\ft)$ is bounded above by $\varphi(\fd) q^{\frac{-\deg \fd}{r}}\leq q^{\left(1-\frac{1}{r}\right)\deg \fd}$. Thus one has 
	\begin{align}\label{J-trivial}
		J_\phi(\ft) \leq \tau(\ft)q^{\left(\frac{r-1}{r}\right)\deg \ft},
	\end{align} 
	where $\tau(\ft)$ is the number of monic polynomials in $A$ dividing $\ft$.
	
	We now obtain a nontrivial bound for $J_\phi(\ft)$. Let $\ft=\fp^\alpha$, a prime power. Then
	\begin{align}
		J_\phi(\fp^\alpha)
		=1+\sum_{i=1}^\alpha
		(q^{i \deg \fp}-q^{(i-1)\deg \fp})
		q^{-\lceil i/r\rceil \deg \fp}
		=\left(1-\frac{1}{q^{\deg \fp}}\right)
		H_\phi(\fp^\alpha)
		+\frac{1}{q^{\deg \fp}},
		\label{eq:J in terms of H}
	\end{align}
	where
	\begin{align}\label{def:H}
		H_\phi(\fp^\alpha)
		\coloneqq
		\sum_{i=0}^\alpha
		q^{(i-\lceil i/r\rceil)\deg \fp}.
	\end{align}
	This definition of $H_\phi$ on prime powers yields the multiplicative function
	\begin{align}
		H_\phi(\ff)
		\coloneqq
		\summonics_{\fg:\,\fg\mid\ff}
		\sum_{\substack{\fh\in A:\,\fg\mid\fh^r\\
				\deg\fh<\deg\fg}}1,
	\end{align}
	where the inner sum is over all polynomials $\fh\in A$, including $\fh=0$. Indeed, for $\fg=\fp^i$, the condition
	$\fp^i\mid\fh^r$ is equivalent to
	$\fp^{\lceil i/r\rceil}\mid\fh$, and hence the inner sum has exactly
	$q^{(i-\lceil i/r\rceil)\deg\fp}$ elements. This function $H_\phi(\ff)$ is
	the function field analogue of the function $H(f)$ defined in
	\cite[Theorem 3]{Akbal-Guloglu}.
	
	From \eqref{def:H}, it is easy to see that
	$H_\phi(\fp^\alpha)\geq 2$ whenever $\alpha\geq1$. Hence, for
	$\alpha\geq1$, one has
	$\frac{1}{q^{\deg\fp}}
		\leq
		\frac{H_\phi(\fp^\alpha)}{2q^{\deg\fp}}$,
	and substituting this inequality into \eqref{eq:J in terms of H} yields
	\begin{align}
		\left(1-\frac{1}{q^{\deg \fp}}\right)H_\phi(\fp^\alpha)
		\leq J_\phi(\fp^\alpha)
		\leq
		\left(1-\frac{1}{2q^{\deg \fp}}\right)H_\phi(\fp^\alpha),
	\end{align}
	for prime powers $\fp^\alpha$ with $\alpha\geq1$. Moreover,
	$J_\phi(1)=H_\phi(1)$. Hence, by multiplicativity of $H_\phi$ and
	$J_\phi$, for any monic $\ft\in A$,
	\begin{align}\label{penul-nontrivial-bound-J}
		J_\phi(\ft)
		\leq
		H_\phi(\ft)
		\prod_{\text{primes }\fp\mid\ft}
		\left(1-\frac{1}{2q^{\deg\fp}}\right).
	\end{align}
	Following \cite[Remark 1, p.1284]{Akbal-Guloglu}, one can show analogously that
	\begin{align}\label{bound-H}
		H_\phi\bigg(\prod_{i\leq k}\fp_i^{\alpha_i}\bigg)
		\leq
		2^k
		\sigma\bigg(
		\prod_{i\leq k}
		\fp_i^{\lfloor \alpha_i(r-1)/r\rfloor}
		\bigg), \quad \text{and} \quad
		H_\phi(\ft)
		\leq
		\frac94\,2^{\omega(\ft)}
		q^{\left(\frac{r-1}{r}\right)\deg\ft},
	\end{align}
	where $
		\sigma(\ff)\coloneqq
		\sum'_{\fd\mid\ff}q^{\deg\fd}$,
	and $\omega(\ft)$ denotes the number of distinct prime divisors of $\ft$.
	Thus, by \eqref{penul-nontrivial-bound-J} and \eqref{bound-H}, we have
	\begin{align}\label{J-nontrivial-bound}
		J_\phi(\ft)
		\leq
		\frac94\,2^{\omega(\ft)}
		q^{\left(\frac{r-1}{r}\right)\deg \ft}
		\prod_{\text{primes }\fp\mid \ft}
		\left(1-\frac{1}{2q^{\deg \fp}}\right).
	\end{align}

	\subsection{Final expression for $\pi_\phi(x;\ft,\fa)$ }\label{subsec:CM-final}
	Combining \eqref{eq:main+error}, \eqref{eq:updated-main-without-tail}, \eqref{eq:tail-final-bound}, \eqref{eq:error-final-bound}, we have obtained
	
	\begin{align}\label{optimize_first}
		\pi_\phi(x;\ft,\fa)&=\pi(x)\cdot\delta_{\phi,\ft,\fa}(x)+ E_1 + E_2 + E_{3a} + E_{3b},
	\end{align}
	where
	\begin{align}
		E_1 &= O_{\phi,K}\Big(q^{\frac{x}{2}+y}\Big(\frac{y}{x}+\frac{\deg \ft}{x}\Big)\Big), \quad
		E_2 = O_{\phi,K}\bigg(\frac{q^{x-\left(\frac{r^2}{\rho}-1\right)y}\log y}{x}\bigg),
		\\
		E_{3a} &= O_{K}\Big(\frac{x}{r}q^{\left(\frac{r-1}{r}\right)x}\Big), \qquad
		E_{3b} = O_{K}\Big(q^{\left(\frac{2r-1}{r}\right)x-\alpha_\phi(\ft)- ry}\Big),
	\end{align}
	where, for monics $\ft \in A$, we define
	\begin{align}\label{def:alpha-t}
		\alpha_\phi(\ft)\coloneqq\deg \ft - \log_q J_{\phi}(\ft),
	\end{align}
	with $\log_q$ denoting logarithm to the base $q$ and $J_\phi(\ft) : A \to \mathbb{R}^+$ as in \eqref{def:J-alpha-intro}. In order to optimise our error terms, we make some observations about the function $\alpha_\phi(\ft)$.  
	\begin{lemma}\label{lem:alpha-non-neg}
		Let $\alpha_\phi(\ft)$ be as defined in \eqref{def:alpha-t}. Then $\alpha_\phi(\ft)\geq 0$ for every monic $\ft \in A$.
	\end{lemma}
	\begin{proof}
		We first note that $\alpha(1)=0$. By the multiplicativity of $J_\phi(\ft)$, it follows that when $\fa,\fb\in A$ are coprime, then
		$\alpha(\fa\fb)=\alpha(\fa)+\alpha(\fb)$.
		Hence, it suffices to show that $\alpha(\fp^k)\geq 0$ for every prime $\fp$ and $k\in \mathbb{N}$. Then, by \eqref{def:J-alpha-intro} and \eqref{d_0}, 
		\begin{align}\label{J-phi-pk}
			J_\phi(\fp^k)&=1+\sum_{i=1}^k (q^{i \deg \fp}-q^{(i-1)\deg \fp})q^{-\lceil i/r\rceil \deg \fp}.
		\end{align}
		Since $\lceil i/r\rceil \geq 1$ for every $1\leq i \leq k$, \eqref{J-phi-pk} is bounded above by
		\begin{align}\label{J-phi-pk-b1}
			1+q^{-\deg \fp}\sum_{i=1}^k (q^{i \deg \fp}-q^{(i-1)\deg \fp}) = 1 + q^{-\deg \fp}(q^{k \deg \fp}-1).
		\end{align}
		We consider the difference $q^{k \deg \fp} - J_\phi(\fp^k)$. By \eqref{J-phi-pk} and \eqref{J-phi-pk-b1}, we have that 
		\begin{align}\label{J-phi-penul}
			q^{k \deg \fp} - J_\phi(\fp^k) &\geq q^{k\deg \fp} - q^{(k-1)\deg \fp} +q^{-\deg \fp} - 1\\
			&=(q^{(k-1)\deg \fp}-q^{-\deg \fp})(q^{\deg \fp}-1)\geq 0.
		\end{align}
		As a result, we obtain $q^{k \deg \fp}\geq J_{\phi}(\fp^k)$ and so
		$\alpha(\fp^k)= \log_q \left( \frac{q^{k \deg \fp}}{J_\phi(\fp^k)} \right)\geq 0$.
	\end{proof}

	\subsection{Choice of parameter $y$}\label{subsec:CM-optimization}
	The displayed error bounds are compared up to bounded factors, and we work with sufficiently large $x$. We observe that $E_{3a}$ is independent of $y$, whereas $E_1$ increases and $E_2$ and $E_{3b}$ decrease as $y$ increases. We begin by balancing $E_1$ and $E_{3b}$. At $y=x/r$, the sum $\sum_{\text{error}}$ in \eqref{def:main-and-error} vanishes and hence the errors $E_{3a}$ and $E_{3b}$ vanish. Furthermore, the error $E_1$ is worse than $E_2$ and hence \eqref{optimize_first} reduces to
	\begin{align}\label{eq:endpoint-y}
		\pi_\phi(x;\ft,\fa)&=\pi(x)\cdot\delta_{\phi,\ft,\fa}(x)+ O_{\phi,K}
		q^{\frac{x}{2}+\frac{x}{r}}
		\left(1+\frac{\deg\ft}{x}\right).
	\end{align}
	Let $r \geq 4$. Then, if $y<x/r$, the bound for $E_{3a}$ in \eqref{optimize_first} is worse than \eqref{eq:endpoint-y}. Hence we choose $y=x/r$ for $r\geq 4$ and the required error bound follows from \eqref{eq:endpoint-y}.
	
	Now let $r=2$ or $r=3$. Note that for $y \leq \frac{\rho}{2r^2}x$, the error term $E_{3b}$ can be seen to be worse than the trivial bound $q^{x-\deg \ft}$, using $\alpha_\phi(\ft)\leq \deg \ft$ and $\rho \leq r$. Hence we will restrict $y$ to be $>\frac{\rho}{2r^2}x$. In this range of $y$, note that $E_1$ dominates over $E_2$. So, we balance $E_1$ and $E_{3b}$. Since $y\leq x/r$, we have
	$\frac{y}{x}+\frac{\deg\ft}{x}\leq 1+\frac{\deg\ft}{x}$.
	This suggests the choice
	\begin{align}\label{def:y-2}
		y_2(x)\coloneqq
		\frac{1}{r+1}\left(
		\frac{3r-2}{2r}x-\alpha_\phi(\ft)
		-\log_q\left(1+\frac{\deg\ft}{x}\right)
		\right),
	\end{align}
	which yields the bound
	\begin{align}\label{eq:balanced-y2}
		E_1 = E_{3b} = O_{\phi, K}\left(q^{\frac{r^2+4r-2}{2r(r+1)}x
			-\frac{\alpha_\phi(\ft)}{r+1}}
		\left(1+\frac{\deg\ft}{x}\right)^{\frac{r}{r+1}}\right).
	\end{align}
	Adding $E_{3a}$ to the error in \eqref{eq:balanced-y2} completes the proof of Theorem \ref{thm:estimation}.

	\section{Proof of Theorem \ref{thm:uniform-AP}}\label{sec:uniform}

		Put $d=\deg\ft$. By Theorem
		\ref{thm:estimation}, we have
		\begin{align}
			|\pi_\phi(x;\ft,\fa)
			-\delta_{\phi,\ft,\fa}\pi(x)|
			\ll_{\phi,K}
			\begin{cases}
				xq^{x/2}
				+q^{5x/6-\alpha_\phi(\ft)/3}
				(1+d/x)^{2/3}, & r=2,\\[3pt]
				xq^{2x/3}
				+q^{19x/24-\alpha_\phi(\ft)/4}
				(1+d/x)^{3/4}, & r=3,\\[3pt]
				q^{(1/2+1/r)x}(1+d/x),
				& r\geq4.
			\end{cases}
		\end{align}
		Since $\alpha_\phi(\ft)\leq d\leq(\lambda_r-\varepsilon)x $,
		this proves the asserted error estimates in \eqref{asymp}.
		
		It remains to show that these errors are
		$o(\delta_{\phi,\ft,\fa}\pi(x))$
		uniformly in the stated ranges.
		By Corollary \ref{cor:positivity-final-explicit},
		there exists a constant $K_\phi>0$,
		independent of $\ft$ and $\fa$, such that
		$\delta_{\phi,\ft,\fa}
		\geq\frac{K_\phi}{\varphi(\ft)}$. Note also that the bound \eqref{J-trivial} on $J_\phi(\ft)$ yields  $\alpha_\phi(\ft)
		\geq d/r-\log_q\tau(\ft)$.
		Consequently, we obtain
		\begin{align}
			\frac{
				|\pi_\phi(x;\ft,\fa)
				-\delta_{\phi,\ft,\fa}\pi(x)|}
			{\delta_{\phi,\ft,\fa}\pi(x)}
			\ll_{\phi,K}
			\begin{cases}
				x\tau(\ft)^{1/3}
				q^{5d/6-x/6}, &r=2,\\[3pt]
				x\tau(\ft)^{1/4}
				q^{11d/12-5x/24}, &r=3,\\[3pt]
				xq^{d-(r-2)x/(2r)}, &r\geq4.\label{uniform}
			\end{cases}
		\end{align}
		Since $d\leq(\lambda_r-\varepsilon)x$ and $\tau(\ft)=q^{o(x)}$ uniformly for $d\leq x$, each of the
		bounds in \eqref{uniform} tends to zero, as 
		$x\to\infty$ and this yields \eqref{asymp r=4} and \eqref{asymp r=2,3} as needed.

	\section{The Case when $\End_{\overline{K}}(\phi)=A$}\label{sec:non-CM}
	
	\subsection{Preliminaries}
	
	Let $\phi$ be a Drinfeld $A$-module over $K$. Throughout this subsection, we assume that the endomorphism ring $\text{End}_{\overline{K}}(\phi)=A$. 
	We say that Galois extensions $L_1/K$ and $L_2/K$ are linearly disjoint over $K$ when 
	$\Gal(L_1L_2/K)\cong\Gal(L_1/K)\times \Gal(L_2/K)$.
	The following lemma is a corollary to the work of Pink and R\"{u}tsche \cite[Theorem 0.1]{PinkRutsche}. 
	\begin{lemma}{\cite[Corollary 9]{Kuo-Tweedle}}\label{lem:M-phi}
		There exists $M(\phi)\in A$ such that
		\begin{enumerate}
			\item if $\gcd(\fm,M(\phi))=1$, then $\Gal(K(\phi[\fm])/K)\cong\GL_r(A/\fm A)$;
			\item if $\fm=\fm_1\fm_2$, where $\fm_1$ is composed of primes that divide $M(\phi)$ and $\fm_2$ is coprime to $M(\phi)$, then $K(\phi[\fm_1])$ and $K(\phi[\fm_2])$ are linearly disjoint over $K$;
			\item if $\gcd(\fm_1,M(\phi))=\gcd(\fm_2,M(\phi))=\gcd(\fm_1,\fm_2)=1$, then $K(\phi[\fm_1])$ and $K(\phi[\fm_2])$ are linearly disjoint over $K$.
		\end{enumerate}
	\end{lemma}
	For the proof of Theorem \ref{thm:non-CM-density}, we will need to append primes of degree $1$ to $M(\phi)$, by defining  $	N(\phi)$ as in \eqref{def:N-phi}. 
	Then, note that all the points of Lemma \ref{lem:M-phi} still hold if we replace $M(\phi)$ by $N(\phi)$.
	\begin{lemma}\label{lem:N-phi}
		Let $N(\phi)$ be as defined in \eqref{def:N-phi}. Then, we have that
		\begin{enumerate}
			\item if $\gcd(\fm,N(\phi))=1$, then $\Gal(K(\phi[\fm])/K)\cong\GL_r(A/\fm A)$;
			\item if $\fm=\fm_1\fm_2$, where $\fm_1$ is composed of primes that divide $N(\phi)$ and $\fm_2$ is coprime to $N(\phi)$, then $K(\phi[\fm_1])$ and $K(\phi[\fm_2])$ are linearly disjoint over $K$;
			\item if $\gcd(\fm_1,N(\phi))=\gcd(\fm_2,N(\phi))=\gcd(\fm_1,\fm_2)=1$, then $K(\phi[\fm_1])$ and $K(\phi[\fm_2])$ are linearly disjoint over $K$.
		\end{enumerate}
	\end{lemma}
	\begin{proof}
		We note that for any $\fa \in A$ with $\gcd(\fa,N(\phi))=1$, we also have $\gcd(\fa,M(\phi))=1$. As a result, parts $(1)$ and $(3)$ follow from parts $(1)$ and $(3)$ of Lemma \ref{lem:M-phi}. 
		
		We now proceed to prove $(2)$. Note that for coprime $\fa,\fb \in A$, we have that
		$K(\phi[\fa\fb])=K(\phi[\fa])K(\phi[\fb])$.
		Let $\fm=\fm_1\fm_2$, where $\fm_1$ is made up of prime divisors of $N(\phi)$ and $\gcd(\fm_2,N(\phi))=1$. As a result, $\gcd(\fm_2,M(\phi))=1$ as well. Also, let $\fm_1=\fm_{11}\fm_{12}$, where $\fm_{11}$ is made up of prime divisors of $M(\phi)$ and $\gcd(\fm_{12},M(\phi))=1$. Note that $(\fm_{11},\fm_{12})=1$. Then,
		\begin{enumerate}
			\item[(i)] by part $(2)$ of Lemma \ref{lem:M-phi}, $K(\phi[\fm_{11}])$ and $K(\phi[\fm_{12}\fm_{2}])$ are linearly disjoint over $K$;
			\item[(ii)] by part $(3)$ of Lemma \ref{lem:M-phi}, $K(\phi[\fm_{12}])$ and $K(\phi[\fm_{2}])$ are linearly disjoint over $K$;
			\item[(iii)] by part $(2)$ of Lemma \ref{lem:M-phi}, 
			\begin{align}
				\Gal(K(\phi[\fm_{11}])/K)\times \Gal(K(\phi[\fm_{12}])/K)\cong\Gal(K(\phi[\fm_{1}])/K).
			\end{align}
		\end{enumerate} 
		Hence, we have that
		$K(\phi[\fm_1])$ and $K(\phi[\fm_2])$ are linearly disjoint over $K$.
	\end{proof}
	In the course of the proof, we come across intersections of the form
	$K(\phi[\fm])\cap K(\Carl[\ft])$.
	These can be explicitly characterized for $\fm$ coprime to $N(\phi)$. We  first compute the commutator subgroup of $\GL_r(A/\fm A)$, which will be useful later.
	\begin{lemma}\label{lem:commutator-1}
		Let $\fm \in A$. Then for $r\geq 3$, the commutator subgroup
		\begin{align}
			[\GL_r(A/\fm A),\GL_r(A/\fm A)]=\SL_r(A/\fm A).
		\end{align}
	\end{lemma}
	\begin{proof}
		Since $A$ is a Principal Ideal Domain, by the Chinese Remainder Theorem it suffices to show
		\begin{align}
			H_r \coloneqq [\GL_r(A/\fp^kA),\GL_r(A/\fp^kA)]=\SL_r(A/\fp^kA),
		\end{align}
		for all primes $\fp \in A$ and $k\in \mathbb{N}$. It is clear that $H_r \subseteq \SL_r(A/\fp^kA)$ since every element in the commutator $H_r$ has determinant $1$. We will now show that $H_r \supseteq \SL_r(A/\fp^kA)$. 
		
		For $1 \leq i, j \leq r$ with $i\neq j$ and $\lambda\in A/\fp^kA$, let $E_{ij}(\lambda)$ denote the elementary $r \times r$ matrix which has $1$'s along the diagonal, $\lambda$ at the $(i, j)$-th entry and $0$'s elsewhere. Since $R = A/\fp^kA$ is a local ring, these matrices generate $\SL_r(R)$ (see \cite[Lemma 2.1(iii)]{Eardley-Manoharmayum}). By \cite[Lemma 2.1(i)(b)]{Eardley-Manoharmayum}, for $r> 2$, and $1\leq i,j,k \leq r$ all distinct, we have the identity
		$[E_{ij}(\alpha),E_{jk}(\beta)]=E_{ik}(\alpha\beta)$.
		As a result, all elementary matrices $E_{ik}(\lambda)$ belong to the commutator subgroup, and hence 
		$\SL_r(A/\fp^kA) \subseteq H_r$, 
		as needed.
	\end{proof}

	\begin{lemma}\label{lem:commutator-2}
		Let $\fm \in A$ be such that $\deg \fp \geq 2$ for all primes $\fp \mid \fm$. Then, the commutator subgroup
		$$[\GL_2(A/\fm A),\GL_2(A/\fm A)]=\SL_2(A/\fm A).$$
	\end{lemma}
	\begin{proof}
		As in the beginning of the previous proof, we only need to show that 
$
			\SL_2(A/\fp^k A)\subseteq [\GL_2(A/\fp^k A),\GL_2(A/\fp^k A)].
$
		Consider the quotient map
		\begin{align}
			\xi : A/\fp^kA \to \frac{A/\fp^kA}{\fp A/\fp^kA} \cong \frac{A}{\fp A}.
		\end{align}
		Since $\deg(\fp)\geq 2$, the field $A/\fp A$ has size $q^{\deg \fp} \geq 4$. Thus, there exists a $\beta' \in A/\fp A$ such that $\beta' \neq 0$ and $\beta' - 1\neq 0$. Let $\beta \in \xi^{-1}(\beta')$ be a preimage of $\beta'$. Then $\beta \in A/\fp^kA$ is such that $\beta \not \in \fp A/\fp^kA$ and $\beta-1 \not \in \fp A/\fp^kA$. 
		
		Since $\fp A/\fp^kA$ is the unique maximal ideal of $A/\fp^kA$, we obtain that $\beta$ and $\beta-1$ are units in $A/\fp^kA$.
		Then for any $\lambda \in A/\fp^kA$, we have
		\begin{align}\label{commutator-r=2}
			E_{12}(\lambda) = [\text{Diag}(\beta,1), E_{12}((\beta-1)^{-1}\lambda)],
		\end{align}
		where $\text{Diag}(\beta,1)$ is the matrix $ \begin{pmatrix}
			\beta & 0\\
			0 & 1
		\end{pmatrix}.$
		By taking transpose on \eqref{commutator-r=2}, we see that $E_{21}(\lambda)$ is also a commutator. By \cite[Lemma 2.1(iii)]{Eardley-Manoharmayum}, $\SL_2(A/\fp^k A)$ is generated by the matrices $E_{12}(\lambda)$, $E_{21}(\lambda)$, $\lambda \in A/\fp^k A$. Hence,
		$\SL_2(A/\fp^kA) \subseteq[\GL_2(A/\fp^kA),\GL_2(A/\fp^kA)]$.
	\end{proof}

	Recall that $\phi^1$, as defined in \eqref{def0:phi1}, denotes the rank $1$ Drinfeld $A$-module over $K$ associated to $\phi$ via the Weil pairing of Tate modules. 
	
	\begin{lemma}\label{lem : max-abelian}
		Let $\End_{\overline{K}}(\phi)=A$. Let $N(\phi)$ be as defined in \eqref{def:N-phi}. Let $\fm \in A$ with $\gcd(\fm,N(\phi))=1$. Then the maximal abelian Galois extension of $K$ contained in $K(\phi[\fm])$ is $K(\phi^1[\fm])$. 
	\end{lemma}
	
	\begin{proof}
		Let $L$ be the maximal abelian Galois extension of $K$ in $K(\phi[\fm])$. We have the group isomorphism
		\begin{align}
			\Gal(L/K) \cong \frac{\Gal(K(\phi[\fm])/K)}{[\Gal(K(\phi[\fm])/K),\Gal(K(\phi[\fm])/K)]}.
		\end{align}
		Since $\gcd(\fm,N(\phi))=1$, we have that $\Gal(K(\phi[\fm])/K)\cong\GL_r(A/\fm A)$ by Lemma \ref{lem:N-phi}(1). By the construction of $N(\phi)$ (see \eqref{def:N-phi}) and Lemmas \ref{lem:commutator-1}, \ref{lem:commutator-2}, the commutator subgroup of $\GL_r(A/\fm A)$ is $\SL_r(A/\fm A)$. Thus,
$
			\Gal(L/K) \cong ( {\GL_r(A/\fm A)}/ {\SL_r(A/\fm A)} ) \cong (A/\fm A)^*.
$
		By \cite[Theorem 3.1.3]{Thakur}, $K(\phi^1[\fm])\subseteq K(\phi[\fm])$ is an abelian extension of $K$. Hence we must have 
		$K(\phi^1[\fm])\subseteq L$. Moreover since $\Gal(K(\phi^1[\fm])/K)\cong(A/\fm A)^*$ (again by \cite[Theorem 3.1.3]{Thakur}), we see that $L = K(\phi^1[\fm])$.
	\end{proof}

	\begin{lemma}\label{lcm}
		Let $\phi$ be an arbitrary Drinfeld $A$-module over $K$. For any $\fa, \fb \in A$, we have that the compositum 
		$K(\phi[\fa])K(\phi[\fb])=K(\phi[\lcm(\fa,\fb)])$,
		where $\lcm(\fa,\fb)$ is the   least common multiple of $\fa, \fb$ in $A$. 
	\end{lemma}
	
	\begin{proof}
		Recall that for a monic $\fa \in A$, the torsion group
		$\phi[\fa]=\{\alpha \in \overline{K}:\phi_\fa(\alpha)=0\}$.
		Hence, we have the set inclusions $\phi[\fa]\subseteq\phi[\lcm(\fa,\fb)]$ and $\phi[\fb]\subseteq\phi[\lcm(\fa,\fb)]$. Thus, 
		\begin{align}
			K(\phi[\fa])K(\phi[\fb])\subseteq K(\phi[\lcm(\fa,\fb)]).
		\end{align}
		To prove the reverse inclusion, let $\fa'=\lcm(\fa,\fb)/\fa$ and $\fb'=\lcm(\fa,\fb)/\fb$. Then, as $\gcd(\fa',\fb')=1$, there exist $\fe,\ff\in A$ such that $\fe \fa' + \ff \fb' =1$. So, we have that for every $\alpha \in \overline{K}$,
		\begin{align}\label{eq:bezout}
			\phi_{\fe}(\phi_{\fa'}(\alpha))+\phi_{\ff}(\phi_{\fb'}(\alpha))=\phi_{\fe\fa' + \ff\fb'}(\alpha)=\alpha.
		\end{align}
		Now, let $\alpha \in \phi[\lcm(\fa,\fb)]$ be arbitrary. Then, we have 
		\begin{align}
			\phi_\fa(\phi_{\fe}(\phi_{\fa'}(\alpha))) &= \phi_{\fe}(\phi_{\lcm(\fa,\fb)}(\alpha))=\phi_{\fe}(0)=0\\
			\phi_\fb(\phi_{\ff}(\phi_{\fb'}(\alpha)))  &= \phi_{\ff}(\phi_{\lcm(\fa,\fb)}(\alpha))=\phi_{\ff}(0)=0
		\end{align}
		Hence, we obtain $\phi_{\fe}(\phi_{\fa'}(\alpha)) \in \phi[\fa]$ and $\phi_{\ff}(\phi_{\fb'}(\alpha)) \in \phi[\fb]$,
		and as a result, by \eqref{eq:bezout},
		$\alpha \in \phi[\fa]+\phi[\fb]$.
		Thus, 
		$K(\phi[\lcm(\fa,\fb)]) \subseteq K(\phi[\fa])K(\phi[\fb])$.	\end{proof}
	
	\begin{lemma}\label{carl-gcd}
		For any nonzero monic $\fa,\fb\in A$, we have $K(\Carl[\fa])\cap K(\Carl[\fb])
			=
			K(\Carl[\gcd(\fa,\fb)])$.
	\end{lemma}
	
	\begin{proof}
		Since $\gcd(\fa,\fb)$ divides both $\fa$ and $\fb$, we have
		\begin{align}\label{containment}
			K(\Carl[\gcd(\fa,\fb)])
			\subseteq
			K(\Carl[\fa])\cap K(\Carl[\fb]).
		\end{align}
		On the other hand, by Lemma \ref{lcm},
		$K(\Carl[\fa])K(\Carl[\fb])
			=
			K(\Carl[\lcm(\fa,\fb)])$.
		Hence,
		\begin{align}
			[K(\Carl[\fa])\cap K(\Carl[\fb]):K]
			&=
			\frac{[K(\Carl[\fa]):K][K(\Carl[\fb]):K]}
			{[K(\Carl[\lcm(\fa,\fb)]):K]}\\
			&=
			\frac{\varphi(\fa)\varphi(\fb)}
			{\varphi(\lcm(\fa,\fb))}=
			\varphi(\gcd(\fa,\fb)).
		\end{align}
		Since
	$
			[K(\Carl[\gcd(\fa,\fb)]):K]
			=
			\varphi(\gcd(\fa,\fb)),
$
		the result follows by \eqref{containment}.
	\end{proof}

	Before proving the next lemma, recall that for the rank $1$ Drinfeld $A$-modules $\Carl$ and $\phi^1$, by \cite[Example 2.2.5, p.38]{Thakur}, there exists an isomorphism of Drinfeld modules
	given by
	\begin{align}\label{eq:iso}
	\Carl \to \phi^1,	\quad 	x \mapsto u\cdot x
	\end{align}
	for some $u \in \overline{K}$ satisfying $u^{q-1} \in K$.
	
	\begin{lemma}\label{lem : phi_cap_C}
		Let $\End_{\overline{K}}(\phi)=A$. Let $N(\phi)$ be as defined in \eqref{def:N-phi}. Let $\fm \in A$ with $\gcd(\fm,N(\phi))=1$ and let $u$ be as in \eqref{eq:iso} above. Then, we have
		\begin{equation}
			\label{d/c}
			[K(\phi[\fm]) \cap K(\Carl[\ft]) : K] = \frac{d_{\fm,\ft}}{c_{\fm,\ft}}\varphi(\gcd(\fm,\ft)),
		\end{equation}
		where $c_{\fm,\ft}$ and $d_{\fm,\ft}$ are certain positive divisors of $q-1$, given by 
		\begin{align}
			c_{\fm,\ft}&=[K(\Carl[\lcm(\fm, \ft)])(u):K(\Carl[\lcm(\fm, \ft)])], \text{ and}\\
			d_{\fm,\ft}&=[K(u\Carl[\fm])K(\Carl[\ft])(u):K(u\Carl[\fm])K(\Carl[\ft])]\label{cmt-dmt}.
		\end{align}
	\end{lemma}
	
	\begin{proof}		
		Consider the diagram
		\normalcolor
		\begin{center} 
			\begin{tikzcd}[
				row sep=6em,
				column sep=1.5em,
				cells={nodes={text=black}},
				arrows={draw=black, text=black}
				]
				% First Row
				K(\phi[\fm]) \arrow[d, dash] \arrow[dr, dash] & & \\
				% Second Row
				L \coloneqq K(\phi^1[\fm]) \arrow[dr, dash, "\substack{\text{maximal abelian}\\ \text{Galois extension}}"'] & 
				\begin{matrix} K(\phi[\fm]) \cap K(\Carl[\ft]) \end{matrix} \arrow[l, hook'] \arrow[r, hook] \arrow[d, dash, "\substack{\text{abelian} \\ \text{Galois}}"'] & 
				H \coloneqq K(\Carl[\ft])  \arrow[dl, dash, "\substack{\text{abelian} \\ \text{Galois}}"] \\
				% Third Row
				& K &
			\end{tikzcd}
		\end{center}

		By \cite[Theorem 3.1.3]{Thakur}, $H\coloneqq K(\Carl[\ft])$ is an abelian extension of $K$ (in fact $\Gal(H/K)\cong(A/\ft A)^*$). The abelian Galois extension $(K(\phi[\fm])\cap H)/K$ must hence be contained inside $L$ by Lemma \ref{lem : max-abelian}. We thus have
		\begin{align}
			[K(\phi[\fm])\cap H : K] = [L \cap H : K] 
			&= \frac{[L:K][H:K]}{[LH:K]} = \frac{\varphi(\fm)\varphi(\ft)}{[LH:K]}\label{penultimate-intersection-bound}, 
		\end{align}
		where $LH=K(\phi^1[\fm])K(\Carl[\ft])$.
		
		The isomorphism \eqref{eq:iso} induces an isomorphism of the torsion points
	$\label{iso:phi1-C}
			\Carl[\fm]\simeq \phi^1[\fm],
$
		over $\overline{K}$ given by
		$
			x \mapsto ux.
	$
		As a consequence we have
		$
			K(\phi^1[\fm])=K(u\cdot \Carl[\fm]).
$
		Hence we have
	$
			LH=K(u \Carl[\fm])K(\Carl[\ft]).
$
		We claim that
		\begin{align}\label{LH-u}
			LH(u)=K( \Carl[\fm])K(\Carl[\ft])(u).
		\end{align}
		In order to show this, first observe that one has $K(u\Carl[\fm])\subseteq K(\Carl[\fm])(u)$, which gives the forward inclusion in \eqref{LH-u}. For the reverse inclusion, note that if $\lambda \in \Carl[\fm]$, then writing $\lambda$ as $u^{-1}(u\lambda)$, we see that $\lambda \in K(u\Carl[\fm])(u)$. This yields $K(\Carl[\fm])\subseteq K(u\Carl[\fm])(u)$, as needed.

		Applying Lemma \ref{lcm} to \eqref{LH-u}, we see that
$
			LH(u)=K(\Carl[\lcm(\fm,\ft)])(u)=K_1(u),
$
		denoting $K(\Carl[\lcm(\fm,\ft)])$ as $K_1$. Hence
		\begin{align}
			[LH:K]&=\frac{[LH(u):K]}{[LH(u):LH]}
			=\frac{[K_1(u):K_1][K_1:K]}{[LH(u):LH]}=\frac{c_{\fm,\ft}}{d_{\fm,\ft}}\varphi(\lcm(\fm,\ft)),\label{LH-u-penul}
		\end{align} 
		where $c_{\fm,\ft}$ and $d_{\fm,\ft}$ are as defined in \eqref{cmt-dmt}. Since $u^{q-1}\in K^*$, we see that $c_{\fm,\ft}$ and $d_{\fm,\ft}$ are divisors of $q-1$. Substituting \eqref{LH-u-penul} into  \eqref{penultimate-intersection-bound}, by multiplicativity of $\varphi$, we obtain
	\eqref{d/c}, as needed. 
	\end{proof}
	
	We will also use the following bound. 
	\begin{lemma}\label{lem:recip-tot}
		For every integer $n\geq 1$, we have
		\begin{align}
			\label{sum 1/phi}
			\summonics_{\substack{\fm \text{ square-free}\\ \deg \fm = n}}
			\frac{1}{\varphi(\fm)} \ll 1.
		\end{align}
	\end{lemma}
	
	\begin{proof}
		We have
		\begin{align}
			\summonics_{\substack{\fm \text{ square-free}\\ \deg \fm = n}}
			\frac{1}{\varphi(\fm)}
			&= q^{-n}
			\summonics_{\substack{\fm \text{ square-free}\\ \deg \fm = n}}
			\frac{q^{\deg \fm}}{\varphi(\fm)}
			= q^{-n}
			\summonics_{\substack{\fm \text{ square-free}\\ \deg \fm = n}}
			\prod_{\fp\mid\fm}
			(1-q^{-\deg\fp})^{-1}\\
			&= q^{-n}
			\summonics_{\substack{\fm \text{ square-free}\\ \deg \fm = n}}
			\prod_{\fp\mid\fm}
			\left(1+\frac{1}{q^{\deg\fp}-1}\right)
			= q^{-n}
			\summonics_{\substack{\fm \text{ square-free}\\ \deg \fm = n}}
			\prod_{\fp\mid\fm}
			\left(1+\frac{1}{\varphi(\fp)}\right).\label{recip-tot-1}
		\end{align}
		Writing the inner product in \eqref{recip-tot-1} as a sum then reveals that \eqref{recip-tot-1} is equal to
		\begin{align}
			q^{-n}
			\summonics_{\substack{\fm \text{ square-free}\\ \deg \fm = n}}
			\sum_{\fd\mid\fm}\frac{1}{\varphi(\fd)}
			\leq q^{-n}
			\summonics_{\substack{\fd\text{ square-free}\\ \deg\fd\leq n}}
			\frac{q^{n-\deg\fd}}{\varphi(\fd)}\leq
			\prod_{\fp}
			\left(
			1+\frac{1}{q^{\deg\fp}(q^{\deg\fp}-1)}
			\right).\label{recip-tot-2}
		\end{align}
		The last product converges, since
		\begin{align}
			\sum_{\fp}
			\frac{1}{q^{\deg\fp}(q^{\deg\fp}-1)}
			\ll
			\sum_{\fp}q^{-2\deg\fp}
			\leq
			\sum_{j\geq 1}q^{-j}
			<\infty.
		\end{align}
		Hence, by \eqref{recip-tot-1} and \eqref{recip-tot-2}, we obtain 			\eqref{sum 1/phi}. 
	\end{proof}

	\subsection{An improved estimate for the tail of the density $\delta_{\phi,\ft,\fa}$}
	
	Recall that the tail of the density $\delta_{\phi,\ft,\fa}$ is the error term $E_{\text{tail}}$, satisfying \eqref{Etail-pre}. We refine this estimate in the case when $\End_{\overline{K}}(\phi)=A$. The argument below is analogous to the refinement of Wong \cite[pp. 771-772]{Wong}. We have from \eqref{Etail-pre} that
	\begin{align}
		E_{\text{tail}} \ll_{\phi, K}\pi(x)\summonics_{\substack{\fm \text{ square-free}\\ \deg \fm> y}}\frac{1}{[K(\phi[\fm])K(\Carl[\ft]):K]}.\label{tail-non-CM}
	\end{align}
	Let $\fm=\fm_1\fm_2$ and $\ft=\ft_1\ft_2$, where $\fm_1, \ft_1$ are composed of prime divisors of $N(\phi)$ and $\fm_2,\ft_2$ are coprime to $N(\phi)$. Then, we note that
	\begin{align}
		[K(\phi[\fm])K(\Carl[\ft]):K]
		\geq [K(\phi[\fm_2])K(\Carl[\ft]):K]=
		\frac{[K(\phi[\fm_2]):K][K(\Carl[\ft]):K]}
		{[K(\phi[\fm_2])\cap K(\Carl[\ft]):K]}.\label{nm-bound-non-CM}
	\end{align}
	By \cite[Proposition 12]{Kuo-Liu2009}, we have the bound 
	\begin{align}\label{full-galois}
		[K(\phi[\fm_2]):K]\gg_{\phi,K} \varphi(\fm_2)q^{(r^2-1)\deg\fm_2}.
	\end{align}
	By Lemma \ref{lem : phi_cap_C}, we have that
$\label{eq:non-CM-intersection-tail}
		[K(\phi[\fm_2])\cap K(\Carl[\ft]):K]
		=
		{d_{\fm_2,\ft}}\varphi(\gcd(\fm_2,\ft)) /{c_{\fm_2,\ft}},
$
	where $c_{\fm_2,\ft}$ and $d_{\fm_2,\ft}$ are positive divisors of $q-1$. In particular,
	\begin{align}\label{eq:non-CM-intersection-upper-tail}
		[K(\phi[\fm_2])\cap K(\Carl[\ft]):K]
		\leq(q-1)\varphi(\gcd(\fm_2,\ft)).
	\end{align}
	Hence, plugging \eqref{full-galois} and \eqref{eq:non-CM-intersection-upper-tail} into \eqref{nm-bound-non-CM}, we have
	\begin{align}
		[K(\phi[\fm])K(\Carl[\ft]):K]\gg_{\phi,K}
		\frac{\varphi(\fm_2)q^{(r^2-1)\deg\fm_2}\varphi(\ft)}
		{\varphi(\gcd(\fm_2,\ft))}.
	\end{align}
	Thus, by \eqref{tail-non-CM}, we have
	\begin{align}
		E_{\text{tail}}
		&\ll_{\phi, K}\pi(x)\summonics_{\substack{\fm_1\mid N(\phi)}}\summonics_{\substack{\fm_2 \text{ square-free}\\ \gcd(\fm_2,N(\phi))=1\\ \deg \fm_2> y-\deg \fm_1}}\frac{\varphi(\gcd(\fm_2,\ft))}
		{\varphi(\fm_2)q^{(r^2-1)\deg\fm_2}\varphi(\ft)}\\
		&\ll_{\phi, K}\pi(x)\frac{1}{\varphi(\ft)}\summonics_{\substack{\fm_1\mid N(\phi)}}\summonics_{\substack{\fm_2 \text{ square-free}\\ \gcd(\fm_2,N(\phi))=1\\ \deg \fm_2> y-\deg \fm_1}}\frac{\varphi(\gcd(\fm_2,\ft))}
		{\varphi(\fm_2)q^{(r^2-1)\deg\fm_2}}.\label{summands}
	\end{align}
	Since $\gcd(\fm_2,N(\phi))=1$, we have $\gcd(\fm_2,\ft)=\gcd(\fm_2,\ft_2)$. Let $
	\fg=\gcd(\fm_2,\ft_2)$,
	and write $\fm_2=\fg\fh$. For squarefree $\fm_2$, $\gcd(\fg,\fh)=1$ and hence $\varphi(\fm_2)=\varphi(\fg)\varphi(\fh)$. Thus, the innermost sum in \eqref{summands} is
	\begin{align}
		\leq 
		\summonics_{\fg\mid\ft_2}
		q^{-(r^2-1)\deg\fg}
		\summonics_{\substack{\fh \text{ square-free}\\ \deg\fh>y-\deg\fm_1-\deg\fg}}
		\frac{1}{\varphi(\fh)q^{(r^2-1)\deg\fh}}.\label{summands1}
	\end{align}
	The inner sum in \eqref{summands1} is then bounded above by
	\begin{align}
		\sum_{d>y-\deg\fm_1-\deg\fg}
		\frac{1}{q^{(r^2-1)d}}\summonics_{\substack{\fh \text{ square-free}\\\deg \fh = d}} \frac{1}{\varphi(\fh)}.\label{summands2}
	\end{align}
	By Lemma \ref{lem:recip-tot} and the properties of geometric sums, we have that \eqref{summands2} is 
	\begin{align}\label{summands2-bound}
		\ll_K q^{-(r^2-1)(y-\deg\fm_1-\deg\fg)}.
	\end{align}
	Hence, by \eqref{summands1}-\eqref{summands2-bound}, \eqref{summands} simplifies to
	\begin{align}
		E_{\text{tail}}
		&\ll_{\phi,K}\pi(x)\frac{1}{\varphi(\ft)}
		\summonics_{\fm_1\mid N(\phi)}
		\summonics_{\fg\mid\ft_2}
		q^{-(r^2-1)\deg\fg}
		q^{-(r^2-1)(y-\deg\fm_1-\deg\fg)}\\
		&\ll_{\phi,K}\pi(x)\frac{\tau(\ft_2)}{\varphi(\ft)}q^{-(r^2-1)y}\summonics_{\fm_1\mid N(\phi)}q^{(r^2-1)\deg \fm_1}\ll_{\phi,K}\pi(x)\frac{\tau(\ft_2)}{\varphi(\ft)}q^{-(r^2-1)y}.
	\end{align}
	Using the prime number theorem, we thus obtain
	\begin{align}\label{eq:non-CM-tail-final-new}
		E_{\text{tail}}
		\ll_{\phi,K}
		\frac{q^{x-(r^2-1)y}}{x}\frac{\tau(\ft_2)}{\varphi(\ft)}.
	\end{align}
	Let $E_1$, $E_2$, $E_{3a}$ and $E_{3b}$ be as defined in \eqref{optimize_first}. The estimate \eqref{eq:non-CM-tail-final-new} refines $E_2$. It does not, however, improve the
	final  error bound in Theorem \ref{thm:estimation}, since the analysis in Section \ref{subsec:CM-optimization} for the choice of $y$ remains unchanged.  Indeed, for $r\geq4$, the choice $y=x/r$ makes
	the sum $\sum_{\text{error}}$ in \eqref{def:main-and-error} vanish, and hence the errors $E_{3a}$ and $E_{3b}$ vanish. Furthermore $E_1$ already dominates over
	$E_2$. For $r=2,3$, whenever $y\leq \frac{\rho}{2r^2}x$, then the error $E_{3b}$ is worse than the trivial bound $q^{x-\deg\ft}$ and whenever $y> \frac{\rho}{2r^2}x$, $E_1$ again dominates over $E_2$. So, the optimal choice of $y$ for rank $r=2,3$ is obtained by balancing $E_1$ and $E_{3b}$, as before.

	\section{Positivity of density}\label{sec:density}
	
	Recall from Remark \ref{constant-delta} that if $r_{\fm,\ft}=1$ for every monic $\fm\in A$, then the density $\delta_{\phi,\ft,\fa}(x)$ is independent of $x$. In this case, we denote it by $\delta_{\phi,\ft,\fa}$ (see \eqref{def0:constant-delta}).
	Then, since $\pi_\phi(x;\ft,\fa)\geq 0$ by definition (see \eqref{def:pi-phi-x,t,a}), Theorem \ref{thm:estimation} implies $\delta_{\phi,\ft,\fa}\geq 0$.

	\subsection{Proof of Corollary \ref{cor:density-zero}}
		 Define
		\begin{align}
			\phi(K)_{\text{tors}}\coloneqq\{\alpha\in K:\phi_\fa(\alpha)=0 \text{ for some non-zero }\fa\in A\}.
		\end{align}
		Let $\fp\in A$ be a prime with $\deg \fp=1$ such that $K(\phi[\fp])=K$. By definition, we have
		$\phi[\fp]\subseteq \phi(K)_{\text{tors}}$.
		Since $\phi[\fp]\simeq (A/\fp A)^r$, we get that
		$\phi(K)_{\text{tors}}$ has a submodule isomorphic to $(A/\fp A)^r$.
		Since reduction embeds $\phi(K)_{\text{tors}}$ into
		$\phi(\F_{\mathfrak q})$ for all but finitely many primes
		$\mathfrak q$ (see \cite[Proof of Theorem 3'(i)]{Kuo-Liu2009}),
		there are only finitely many prime reductions
		$\phi(\F_{\mathfrak q})$ that have at most $(r-1)$ cyclic components. Recall by \eqref{def:pi-phi-x,t,a} that
		\begin{align}
			\pi_\phi(x;\ft,\fa)=\#\{\fp\in\mathcal{P}_\phi:\deg\fp=x,\fp\equiv \fa (\ft) \text{ and } \phi(\F_\fp) \text{ has at most } (r-1) \text{ cyclic components}\}.
		\end{align}
		Hence, by Theorem \ref{thm:estimation}, we have that the density $\delta_{\phi,\ft,\fa}=0$. 

		We can also obtain the result by the following computation. For any prime $\fp \in A$, we can rewrite \eqref{def0:constant-delta} as 
		\begin{align}\label{rewrite-constant-delta}
			\delta_{\phi,\ft,\fa} = \summonics_{\substack{\fm\\ \gcd(\fm,\fp)=1}}
			\left(
			\frac{\mu_q(\fm)\gamma_{\fa,\ft}(K(\phi[\fm]))}
			{[K(\phi[\fm])K(\Carl[\ft]):K]}
			+
			\frac{\mu_q(\fp\fm)\gamma_{\fa,\ft}(K(\phi[\fp\fm]))}
			{[K(\phi[\fp\fm])K(\Carl[\ft]):K]}
			\right).
		\end{align}
		Let $\fp$ be such that $\deg \fp =1$ and $K(\phi[\fp])=K$. Then, for $\fm\in A$ with $\gcd(\fm,\fp)=1$, Lemma \ref{lcm} gives
$
			K(\phi[\fp\fm])=K(\phi[\fp])K(\phi[\fm])=K(\phi[\fm]).
$
		Furthermore, for $\fm\in A$ with $\gcd(\fm,\fp)=1$, $\mu_q(\fp\fm)=-\mu_q(\fm)$. Hence, \eqref{rewrite-constant-delta} is the same as 
		\begin{align}
			\summonics_{\substack{\fm\\ \gcd(\fm,\fp)=1}}
			\left(
			\frac{\mu_q(\fm)\gamma_{\fa,\ft}(K(\phi[\fm]))}
			{[K(\phi[\fm])K(\Carl[\ft]):K]}
			-
			\frac{\mu_q(\fm)\gamma_{\fa,\ft}(K(\phi[\fm]))}
			{[K(\phi[\fm])K(\Carl[\ft]):K]}
			\right)=0.
		\end{align}
		This proves the corollary.\\

	\subsection{Preliminaries for Theorem \ref{thm:non-CM-density}}
	
	Throughout this subsection, we assume that $\End_{\overline K}(\phi)=A$.
	Let $N(\phi)$ be as defined in \eqref{def:N-phi}, and define
	\begin{align}
		\mathcal B &= \{\fp \in A: \fp \text{ prime}, \fp \mid N(\phi)\},  \:
		\mathcal G = \{\fp \in A: \fp \text{ prime}, \gcd(\fp,N(\phi))=1\}.
	\end{align}
	For a set $\mathcal S$ of primes, put
$\label{def:D-S-intro}
		D_{\mathcal S}\coloneqq\prod_{\fp\in\mathcal S}K(\phi[\fp]),
$
	where the product denotes the compositum. Then, by Lemma \ref{lcm}, for a finite $\mathcal S$, 
$
		D_{\mathcal S} = K(\phi[\fs])$, $ \fs = \prod_{\fp \in \mathcal S}\fp.
$
	We first record a simple field-theoretic observation.
	
	\begin{lemma}\label{linearly-disjoint-1}
		Let $L_1,L_2,L_3$ be finite Galois extensions of $K$ with
		$L_1\cap L_2=K$.
		Then
		$L_1\cap L_2L_3=K$
		if and only if
		$L_1L_2\cap L_3=L_2\cap L_3$.
	\end{lemma}
	
	\begin{proof}
		We write the degree $[L_1L_2L_3:K]$ in two different ways as
		\begin{align}\label{eq:degree-two-ways}
			\frac{[L_1L_2:K][L_3:K]}
			{[L_1L_2\cap L_3:K]}
			=
			\frac{[L_1:K][L_2L_3:K]}
			{[L_1\cap L_2L_3:K]}.
		\end{align}
		Since $L_1\cap L_2=K$, the left hand side of
		\eqref{eq:degree-two-ways} is equal to
		\begin{align}\label{eq:degree-left}
			\frac{[L_1:K][L_2:K][L_3:K]}
			{[L_1L_2\cap L_3:K]}.
		\end{align}
		Simplifying the numerator on the right hand side of
		\eqref{eq:degree-two-ways} gives
		\begin{align}\label{eq:degree-right}
			\frac{[L_1:K][L_2:K][L_3:K]}
			{[L_1\cap L_2L_3:K][L_2\cap L_3:K]}.
		\end{align}
		Equating \eqref{eq:degree-left} and \eqref{eq:degree-right} yields
		\begin{align}\label{eq:intersection-degree-identity}
			[L_1L_2\cap L_3:K]
			=
			[L_1\cap L_2L_3:K][L_2\cap L_3:K].
		\end{align}
		Hence, if $L_1\cap L_2L_3=K$, then by
		\eqref{eq:intersection-degree-identity},
$
			[L_1L_2\cap L_3:K]=[L_2\cap L_3:K].
$
		But since $L_2\cap L_3\subseteq L_1L_2\cap L_3$, we thus have
	$
			L_1L_2\cap L_3=L_2\cap L_3.
$
		Conversely, if $L_1L_2\cap L_3=L_2\cap L_3$, then by
		\eqref{eq:intersection-degree-identity},
		$[L_1\cap L_2L_3:K]=1$, and hence $L_1\cap L_2L_3=K$.
	\end{proof}
	
	The following lemma isolates all the entanglement which is relevant to the positivity of the density into finitely many prime-division fields.
	
	\begin{lemma}\label{finite-entanglement}
		Let $\End_{\overline K}(\phi)=A$. Then, there exists a finite set
		$\mathcal S\coloneqq\mathcal S_{\phi,\ft}\subseteq\mathcal G$
		such that for every square-free monic $\fm\in A$ supported on
		$\mathcal G\setminus\mathcal S$, one has
		$K(\phi[\fm])\cap
			D_{\mathcal B\cup\mathcal S}K(\Carl[\ft])
			=K$.
	\end{lemma}
	
	\begin{proof}
		Let $\mathcal S \subseteq \mathcal G$ be an arbitrary finite set. Then, by part (3) of Lemma \ref{lem:N-phi},  for every square-free monic $\fm\in A$ supported on
		$\mathcal G\setminus\mathcal S$, we have that
$
			K(\phi[\fm])\cap D_{\mathcal S}= K.
$
		Hence, by Lemma \ref{linearly-disjoint-1}, it suffices to show that there exists a finite $\mathcal S \subseteq \mathcal G$ satisfying
		\begin{align} \label{enough}
			K(\phi[\fm])D_{\mathcal S}
			\cap D_{\mathcal B}K(\Carl[\ft])
			=
			D_{\mathcal S}\cap D_{\mathcal B}K(\Carl[\ft]).
		\end{align} 
		Since $K(\phi[\fm])D_{\mathcal S}
		\subseteq D_{\mathcal G}$, we consider the field
$
			F\coloneqq
			D_{\mathcal G}\cap D_{\mathcal B}K(\Carl[\ft]).
$
		Since $\mathcal B$ is finite, the extension $D_{\mathcal B}K(\Carl[\ft])/K$
		is finite, and hence $F/K$ is also finite. Moreover,
		$
		F\subseteq D_{\mathcal G}
		$.
		Therefore, there exists a finite set
		$\mathcal S\subseteq\mathcal G$ such that
		\begin{align}\label{DS-cap-bad-C}
			F=D_{\mathcal S}\cap D_{\mathcal B}K(\Carl[\ft]).
		\end{align}
		Therefore,
		\begin{align}
			K(\phi[\fm])D_{\mathcal S}
			\cap D_{\mathcal B}K(\Carl[\ft])
			\subseteq
			D_{\mathcal G}\cap D_{\mathcal B}K(\Carl[\ft]) =
			D_{\mathcal S}\cap D_{\mathcal B}K(\Carl[\ft]),
		\end{align}
		where the last equality follows from \eqref{DS-cap-bad-C}. The reverse inclusion is immediate, giving \eqref{enough} and completing the proof. 
	\end{proof}
	
	Let $\mathcal S$ be an entanglement set for $(\phi,\ft)$ as in
	Lemma \ref{finite-entanglement}, and put
$
		L_{\phi,\ft,\mathcal S}
		\coloneqq
		D_{\mathcal B\cup\mathcal S}K(\Carl[\ft]).
$
	Let $\lambda_\ft$ generate the $A/\ft A$-module $\Carl[\ft]$. Recall, by \eqref{def:Omega-intro}, that
	\begin{align}
		\Omega_{\phi,\ft,\fa,\mathcal S}
		\coloneqq
		\{
		\sigma\in
		\Gal(L_{\phi,\ft,\mathcal S}/K):
		\sigma|_{K(\Carl[\ft])}(\lambda_\ft)=\Carl_\fa(\lambda_\ft),
		\sigma|_{K(\phi[\fp])}\neq\mathbbm{1}
		\text{ for every }\fp\in\mathcal B\cup\mathcal S
		\}.
	\end{align}
	
	\begin{lemma}\label{finite-factor-Omega}
		With the notation above,
		\begin{align}\label{finite-factor-sum}
			\summonics_{\substack{\fm\\
					\fp\mid\fm\Rightarrow\fp\in\mathcal B\cup\mathcal S}}
			\frac{\mu_q(\fm)\gamma_{\fa,\ft}(K(\phi[\fm]))}
			{[K(\phi[\fm])K(\Carl[\ft]):K]}
			=
			\frac{|\Omega_{\phi,\ft,\fa,\mathcal S}|}
			{[L_{\phi,\ft,\mathcal S}:K]}.
		\end{align}
	\end{lemma}
	
	\begin{proof}
		Put $G\coloneqq\Gal(L_{\phi,\ft,\mathcal S}/K)$, and define
$
			H_\fp
			\coloneqq
			\{\sigma\in G:
			\sigma|_{K(\phi[\fp])}=\mathbbm{1}\},
		\:	( \fp\in\mathcal B\cup\mathcal S)$, { and } 		$ 	X_{\fa,\ft}
			\coloneqq
			\{\sigma\in G:
			\sigma(\lambda_\ft)=\Carl_\fa(\lambda_\ft)\}.
$
		Then
		\begin{align}\label{Omega-as-complement}
			\Omega_{\phi,\ft,\fa,\mathcal S}
			=
			X_{\fa,\ft}
			\setminus
			\bigcup_{\fp\in\mathcal B\cup\mathcal S}H_\fp.
		\end{align}
		For square-free monic $\fm$ supported on
		$\mathcal B\cup\mathcal S$, the identity
	$
			K(\phi[\fm])
			=
			\prod_{\fp\mid\fm}K(\phi[\fp])
$
		gives
		\begin{align}\label{X-cap-H-description}
			X_{\fa,\ft}\cap\bigcap_{\fp\mid\fm}H_\fp
			=\{
			\sigma\in G:
			\sigma|_{K(\phi[\fm])}=\mathbbm{1},
				\sigma(\lambda_\ft)=\Carl_\fa(\lambda_\ft)
			\}.
		\end{align}
		
		The prescribed automorphisms agree on
		$K(\phi[\fm])\cap K(\Carl[\ft])$ if and only if
		$\gamma_{\fa,\ft}(K(\phi[\fm]))=1$.
		In this case, they extend to the compositum
		$K(\phi[\fm])K(\Carl[\ft])$, and then to
		$L_{\phi,\ft,\mathcal S}$. Thus the set in
		\eqref{X-cap-H-description} is empty when
		$\gamma_{\fa,\ft}(K(\phi[\fm]))=0$, and otherwise is a coset of
$
			\Gal(L_{\phi,\ft,\mathcal S}/K(\phi[\fm])K(\Carl[\ft])).
$
		Consequently,
		\begin{align}
			\Big|
			X_{\fa,\ft}\cap\bigcap_{\fp\mid\fm}H_\fp
			\Big|
			=
			\gamma_{\fa,\ft}(K(\phi[\fm]))
			[L_{\phi,\ft,\mathcal S}:K(\phi[\fm])K(\Carl[\ft])].
		\end{align}
		Dividing by $|G|=[L_{\phi,\ft,\mathcal S}:K]$ yields
		\begin{align}\label{finite-cardinality-identity}
			\frac{1}{|G|}
			\Big|
			X_{\fa,\ft}\cap\bigcap_{\fp\mid\fm}H_\fp
			\Big|
			=
			\frac{\gamma_{\fa,\ft}(K(\phi[\fm]))}
			{[K(\phi[\fm])K(\Carl[\ft]):K]}.
		\end{align}
		Since $\mu_q(\fm)=0$ unless $\fm$ is square-free,
		inclusion-exclusion, \eqref{Omega-as-complement}, and
		\eqref{finite-cardinality-identity} give
		\begin{align}
			\summonics_{\substack{\fm\\
					\fp\mid\fm\Rightarrow\fp\in\mathcal B\cup\mathcal S}}
			\frac{\mu_q(\fm)\gamma_{\fa,\ft}(K(\phi[\fm]))}
			{[K(\phi[\fm])K(\Carl[\ft]):K]}
			&=
			\frac{1}{|G|}
			\summonics_{\substack{\fm\text{ square-free}\\
					\fp\mid\fm\Rightarrow\fp\in\mathcal B\cup\mathcal S}}
			\mu_q(\fm)
			\Big|X_{\fa,\ft}\cap\bigcap_{\fp\mid\fm}H_\fp\Big|\\
			&=
			\frac{1}{|G|}
			\Big|X_{\fa,\ft}\setminus
			\bigcup_{\fp\in\mathcal B\cup\mathcal S}H_\fp\Big| =
			\frac{|\Omega_{\phi,\ft,\fa,\mathcal S}|}
			{[L_{\phi,\ft,\mathcal S}:K]}.
		\end{align}
	\end{proof}
	
	\subsection{Proof of Theorem \ref{thm:non-CM-density}}

		Every square-free $\fm \in A$ may be written uniquely as $\fm=\fm_0\fm_1$,
		where $\fm_0$ is supported on $\mathcal B\cup\mathcal S$ and $\fm_1$ is supported on
		$\mathcal G\setminus\mathcal S$.
		In particular, $\gcd(\fm_0,\fm_1)=1$
		and by Lemma \ref{lem:N-phi},
		\begin{align}\label{m0m1-compositum}
			K(\phi[\fm])
			=
			K(\phi[\fm_0])K(\phi[\fm_1]).
		\end{align}
		Since $\mathcal S$ is an entanglement set for $(\phi,\ft)$, we have by Lemma \ref{finite-entanglement},
		\begin{align}\label{m1-cap-L}
			K(\phi[\fm_1])\cap L_{\phi,\ft,\mathcal S}=K.
		\end{align}
		Since $K(\phi[\fm_0])K(\Carl[\ft])
			\subseteq L_{\phi,\ft,\mathcal S}$,
		we obtain from \eqref{m1-cap-L} that
		\begin{align}\label{m1-disjoint-m0C}
			K(\phi[\fm_1])
			\cap
			K(\phi[\fm_0])K(\Carl[\ft])
			=K.
		\end{align}
		Therefore,
		\begin{align}\label{degree-factorization}
			[K(\phi[\fm])K(\Carl[\ft]):K]
			=
			[K(\phi[\fm_0])K(\Carl[\ft]):K]
			[K(\phi[\fm_1]):K].
		\end{align}
		We also claim that
		\begin{align}\label{gamma-factorization-intersection}
			K(\phi[\fm])\cap K(\Carl[\ft])
			=
			K(\phi[\fm_0])\cap K(\Carl[\ft]).
		\end{align}
		Indeed, by \eqref{m1-cap-L},
	$
			K(\phi[\fm_1])\cap K(\phi[\fm_0])=K,
	$
		and by \eqref{m1-disjoint-m0C},
	$
			K(\phi[\fm_1])
			\cap
			K(\phi[\fm_0])K(\Carl[\ft])
			=K.
$
		Thus, applying Lemma \ref{linearly-disjoint-1} with $L_1, L_2, L_3$ equal to 
	$K(\phi[\fm_1]), K(\phi[\fm_0]), K(\Carl[\ft])$ respectively, 
		and using \eqref{m0m1-compositum}, we obtain \eqref{gamma-factorization-intersection}. Hence,
		\begin{align}\label{gamma-factorization}
			\gamma_{\fa,\ft}(K(\phi[\fm]))
			=
			\gamma_{\fa,\ft}(K(\phi[\fm_0])).
		\end{align}
		The series defining $\delta_{\phi,\ft,\fa}$ is absolutely convergent, as proved at the end of subsection \ref{subsec:CM-main}. Therefore, using \eqref{degree-factorization}, \eqref{gamma-factorization}, and the multiplicativity of $\mu_q$, we may write
		\begin{align}
			\delta_{\phi,\ft,\fa}
			&=
			\summonics_{\substack{\fm_0\\
					\fp\mid\fm_0\Rightarrow\fp\in\mathcal B\cup\mathcal S}}
			\summonics_{\substack{\fm_1\\
					\fp\mid\fm_1\Rightarrow\fp\in\mathcal G\setminus\mathcal S}}
			\frac{
				\mu_q(\fm_0)\mu_q(\fm_1)
				\gamma_{\fa,\ft}(K(\phi[\fm_0]))
			}
			{
				[K(\phi[\fm_0])K(\Carl[\ft]):K]
				[K(\phi[\fm_1]):K]
			}\\
			&=
			\Bigg(
			\summonics_{\substack{\fm_0\\
					\fp\mid\fm_0\Rightarrow\fp\in\mathcal B\cup\mathcal S}}
			\frac{
				\mu_q(\fm_0)
				\gamma_{\fa,\ft}(K(\phi[\fm_0]))
			}
			{
				[K(\phi[\fm_0])K(\Carl[\ft]):K]
			}
			\Bigg)
			\Bigg(
			\summonics_{\substack{\fm_1\\
					\fp\mid\fm_1\Rightarrow\fp\in\mathcal G\setminus\mathcal S}}
			\frac{\mu_q(\fm_1)}
			{[K(\phi[\fm_1]):K]}
			\Bigg).
			\label{density-two-factors}
		\end{align}
		By Lemma \ref{finite-factor-Omega}, the first factor in \eqref{density-two-factors} is equal to
		\begin{align}\label{first-factor-Omega}
			\frac{|\Omega_{\phi,\ft,\fa,\mathcal S}|}{[L_{\phi,\ft,\mathcal S}:K]}.
		\end{align}
		We now consider the second factor. If $\fm_1$ is square-free and supported on
		$\mathcal G\setminus\mathcal S$, then by Lemma \ref{lem:N-phi} and the Chinese Remainder Theorem,
$
			[K(\phi[\fm_1]):K]
			=
			\prod_{\fp\mid\fm_1}
			[K(\phi[\fp]):K].
$
		Hence, again using absolute convergence,
		\begin{align}\label{second-factor-product}
			\summonics_{\substack{\fm_1\\
					\fp\mid\fm_1\Rightarrow\fp\in\mathcal G\setminus\mathcal S}}
			\frac{\mu_q(\fm_1)}
			{[K(\phi[\fm_1]):K]}
			=
			\prod_{\fp\notin\mathcal B\cup\mathcal S}
			\left(
			1-\frac{1}{[K(\phi[\fp]):K]}
			\right).
		\end{align}
		Combining \eqref{density-two-factors}, \eqref{first-factor-Omega}, and \eqref{second-factor-product}, we obtain
		\eqref{density-finite-factor-intro}.
		
		It remains to prove the positivity criterion. By part (1) of Lemma \ref{lem:N-phi},
		$[K(\phi[\fp]):K]>1$
		for every $\fp\notin\mathcal B\cup\mathcal S$. Moreover, the absolute convergence established earlier gives
$
			\sum'_{\fp\notin\mathcal B\cup\mathcal S}
			\frac{1}{[K(\phi[\fp]):K]}
			<\infty.
	$
		Therefore,
		\begin{align}
			\prod_{\fp\notin\mathcal B\cup\mathcal S}
			\left(
			1-\frac{1}{[K(\phi[\fp]):K]}
			\right)>0.
		\end{align}
		It then follows that $\delta_{\phi,\ft,\fa}>0$ if and only if $\Omega_{\phi,\ft,\fa,\mathcal S}\neq\varnothing$.
		This proves the theorem.

	\subsection{Proof of Corollary \ref{cor:positivity-final-explicit}}
	
	Let
	$\mathcal S
	\coloneqq
	\{\fp\in\mathcal G:\fp\mid\ft\}$ and $\fs$ be the product of elements in $\mathcal S$. We first note that $\mathcal S$ is an
	entanglement set for $(\phi,\ft)$.
	Indeed, since $\phi^1\simeq_K\Carl$, we have
	$K(\Carl[\ft])
	=
	K(\phi^1[\ft])
	\subseteq K(\phi[\ft])$.
	As $\fs\mid\ft$, it follows that
	$D_{\mathcal S}K(\Carl[\ft])
	\subseteq K(\phi[\ft])$ and thus it suffices to show
	\begin{align}\label{explicit-entanglement}
		K(\phi[\fm])\cap
		D_{\mathcal B}K(\phi[\ft])=K.
	\end{align}
	For every square-free monic
	$\fm$ supported on $\mathcal G\setminus\mathcal S$,
	we have $\gcd(\fm,\ft)=1$ and
	$\gcd(\fm\ft,N(\phi))=1$.
	Hence, by Lemma \ref{lem:N-phi},
	$
	K(\phi[\fm])\cap K(\phi[\ft])=K,
	$
	and
	$
	K(\phi[\fm\ft])\cap D_{\mathcal B}
	=
	K(\phi[\ft])\cap D_{\mathcal B}=K.
	$
	Thus, applying Lemma \ref{linearly-disjoint-1}
	with $L_1=K(\phi[\fm])$,
	$L_2=K(\phi[\ft])$ and $L_3=D_{\mathcal B}$,
	we obtain \eqref{explicit-entanglement}. 
	Thus, $\mathcal S$ is an entanglement set.
	By Theorem \ref{thm:non-CM-density} and Lemma \ref{lem:N-phi}, we have
	\begin{align}
		\delta_{\phi,\ft,\fa}
		=
		\frac{|\Omega_{\phi,\ft,\fa,\mathcal S}|}
		{[L_{\phi,\ft,\mathcal S}:K]}
		\prod_{\fp\nmid N(\phi)\ft}
		\left(1-\frac{1}{|\GL_r(A/\fp A)|}\right).
		\label{density-cor16}
	\end{align}
	Since $D_{\mathcal S}K(\Carl[\ft])
	\subseteq K(\phi[\ft])$ and
	$\gcd(\ft,N(\phi))=1$, Lemma \ref{lem:N-phi}
	gives
	$D_{\mathcal B}\cap
	D_{\mathcal S}K(\Carl[\ft])=K$.
	Hence
	$
	\Gal(L_{\phi,\ft,\mathcal S}/K)
	\simeq
	\Gal(D_{\mathcal B}/K)
	\times
	\Gal(D_{\mathcal S}K(\Carl[\ft])/K)
	$. We hence define
	\begin{align}
		\Omega_{\mathcal B}
		&\coloneqq
		\big\{
		\sigma\in\Gal(D_{\mathcal B}/K):
		\sigma|_{K(\phi[\fp])}\neq\mathbbm{1}
		\text{ for every }\fp\in\mathcal B
		\big\},\\
		\Omega_{\mathcal S,\ft,\fa}
		&\coloneqq
		\big\{
		\sigma\in\Gal(D_{\mathcal S}K(\Carl[\ft])/K):
		\sigma(\lambda_\ft)=\Carl_\fa(\lambda_\ft),
		\sigma|_{K(\phi[\fp])}\neq\mathbbm{1}
		\text{ for every }\fp\in\mathcal S
		\big\},
	\end{align}
	where $\lambda_\ft$ generates the
	$A/\ft A$-module $\Carl[\ft]$. Consequently,
	\begin{align}\label{Omega-product}
		\frac{|\Omega_{\phi,\ft,\fa,\mathcal S}|}
		{[L_{\phi,\ft,\mathcal S}:K]}
		=
		\frac{|\Omega_{\mathcal B}|}
		{[D_{\mathcal B}:K]}
		\cdot
		\frac{|\Omega_{\mathcal S,\ft,\fa}|}
		{[D_{\mathcal S}K(\Carl[\ft]):K]}.
	\end{align}
	
	Put
	$\fb\coloneqq\prod_{\fp\in\mathcal B}\fp$ and for every $\fp\in\mathcal B$, put
	$H_\fp\coloneqq
	\{\sigma\in\Gal(D_{\mathcal B}/K):
	\sigma|_{K(\phi[\fp])}=\mathbbm{1}\}$.
	For every square-free $\fd\mid\fb$,
	Lemma \ref{lcm} gives
	$
	\bigcap_{\fp\mid\fd}H_\fp
	=
	\Gal(D_{\mathcal B}/K(\phi[\fd]))$.
	Thus, by inclusion-exclusion, 
	\begin{align}
		\frac{
			|\Omega_{\mathcal B}|}
		{[D_{\mathcal B}:K]}
		=
		\sum_{\fd\mid\fb}
		\frac{\mu_q(\fd)}{[K(\phi[\fd]):K]}.
		\label{C-phi-factor}
	\end{align}
	We call this factor $C_\phi$. Since $\phi^1\simeq_K\Carl$
	and $q>2$, there exists
	$c\in\F_q^*\setminus\{1\}$ whose
	corresponding automorphism on
	$K(\Carl[\fb])\subseteq D_{\mathcal B}$
	acts nontrivially on $K(\Carl[\fp])$
	for every $\fp\in\mathcal B$.
	Extending this automorphism to
	$D_{\mathcal B}$ shows that $\Omega_{\mathcal B}$ is nonempty.
	Hence, $C_\phi>0$. Thus, by \eqref{density-cor16}, \eqref{Omega-product} and \eqref{C-phi-factor}, it suffices to show that
	\begin{align}\label{Omega-second-factor}
		\frac{|\Omega_{\mathcal S,\ft,\fa}|}
		{[D_{\mathcal S}K(\Carl[\ft]):K]} = \frac{1}{\varphi(\ft)}
		\prod_{\substack{\fp\mid\ft\\ \fp \mid \fa-1}}
		\left(
		1-\frac{1}
		{|\SL_r(A/\fp A)|}
		\right)>0.
	\end{align}

	Since $K(\Carl[\ft])/K$ is abelian,
	by Lemma \ref{lem : max-abelian} and
	$\phi^1\simeq_K\Carl$, we have
	\begin{align}
		D_{\mathcal S}\cap K(\Carl[\ft])
		=K(\Carl[\fs])\label{DSC-intersection}
	\end{align}
	Since $\gcd(\fs,N(\phi))=1$, Lemma
	\ref{lem:N-phi} gives
	$[D_{\mathcal S}:K]
	=\prod_{\fp\in\mathcal S}|\GL_r(A/\fp A)|$.
	Hence,
	\begin{align}
		[D_{\mathcal S}K(\Carl[\ft]):K]
		=\frac{[D_{\mathcal S}:K]\varphi(\ft)}
		{\varphi(\fs)}
		=\varphi(\ft)
		\prod_{\fp\in\mathcal S}|\SL_r(A/\fp A)|.
		\label{DS-C-degree}
	\end{align}
	Let $\sigma_{\fa,\ft}\in\Gal(K(\Carl[\ft])/K)$
	be the automorphism satisfying
	$\sigma_{\fa,\ft}(\lambda_\ft)
	=\Carl_\fa(\lambda_\ft)$.
	For every $\fp\in\mathcal S$, the restriction map
	$
	\Gal(K(\phi[\fp])/K)
	\longrightarrow
	\Gal(K(\Carl[\fp])/K)
	$
	is surjective, with kernel of cardinality
	$|\SL_r(A/\fp A)|$, by Lemma \ref{lem:N-phi}
	and the identification $\phi^1\simeq_K\Carl$.
	Thus, there are exactly
	$|\SL_r(A/\fp A)|$ automorphisms of
	$K(\phi[\fp])$ extending
	$\sigma_{\fa,\ft}|_{K(\Carl[\fp])}$. 
	If $\fa\not\equiv1\pmod{\fp}$, every such
	extension is nontrivial. If
	$\fa\equiv1\pmod{\fp}$, these extensions
	form the group
	$\Gal(K(\phi[\fp])/K(\Carl[\fp]))$,
	and hence exactly one of them is the identity.
	Thus, in either case, the number of
	nontrivial extensions is
	$
	|\SL_r(A/\fp A)|
	-\mathbf{1}_{\fa\equiv1\pmod{\fp}}$.
	
	By Lemma \ref{lem:N-phi}, the fields
	$\{K(\phi[\fp])\}_{\fp\in\mathcal S}$
	are jointly linearly disjoint over $K$.
	Thus, these extensions can be chosen
	independently and combined into an automorphism
	$\tau\in\Gal(D_{\mathcal S}/K)$.
	Moreover, by \eqref{DSC-intersection},
	$\tau$ and $\sigma_{\fa,\ft}$ agree on
	$D_{\mathcal S}\cap K(\Carl[\ft])$,
	and hence extend uniquely to an automorphism
	of $D_{\mathcal S}K(\Carl[\ft])$.
	Conversely, every automorphism in
	$\Omega_{\mathcal S,\ft,\fa}$ arises in this way.
	Consequently,
	\begin{align}
		|\Omega_{\mathcal S,\ft,\fa}|
		=
		\prod_{\fp\in\mathcal S}
		\left(
		|\SL_r(A/\fp A)|
		-\mathbf{1}_{\fa\equiv1\pmod{\fp}}
		\right).\label{Omega-sf-use}
	\end{align}
	Therefore, the equality in \eqref{Omega-second-factor} follows by \eqref{DS-C-degree} and \eqref{Omega-sf-use}.
	
	Recall that $\fa \mid N(\phi)$ every $\fa \in A$ with $\deg \fa =1$. Hence the positivity of \eqref{Omega-second-factor} is immediate since $\gcd(\ft,N(\phi))=1$ and thus $|\SL_r(A/\fp A)|>1$ for every $\fp \mid \ft$.

	\normalcolor
	
	\subsection{Proof of Example \ref{prop:non-Carlitz-positive-intro}}
	
	Let $\phi_0$ be the rank-$2$ Drinfeld $A$-module over $K$ defined by
$\label{def:phi0-example}
		{\phi_0}_T=T+\tau-T^4\tau^2.
$
	By \cite[Theorem 1.2]{Zywina2011}, for every nonzero monic $\fm\in A$,
	\begin{align}\label{phi0-full-image-example}
		\Gal(K(\phi_0[\fm])/K)
		\simeq
		\GL_2(A/\fm A).
	\end{align}
	The associated rank-$1$ Drinfeld module is
$
		{\phi_0^1}_T
		=
		T+T^4\tau.
$
	Note that $T^{-1}\Carl_TT = T^{-1}(T+\tau)T=T+T^4\tau$ as well
	and hence $\phi_0^1 \simeq \Carl$, and
	\begin{align}\label{phi01-Carl-fields}
		K(\phi_0^1[\fm])
		=
		K(\Carl[\fm])
	\end{align}
	for every nonzero monic $\fm\in A$.
	Let $u$ be a nonzero $T$-torsion point of $\Carl$. Then $u^4=-T$ and $K(u)=K(\Carl[T])$
	We consider the Drinfeld module $\phi \coloneqq u^{-1}\phi_0u$ given by
$
		\phi_T = T+u^4\tau-T^4u^{24}\tau^2.
$
	Using $u^4 =-T$, we have
$
		\phi_T
		=
		T-T\tau-T^{10}\tau^2.
$
	We state some relevant properties of $\phi$ .
	
	\begin{lemma}\label{lem:non-Carlitz-example}
		Let $\phi$ be the Drinfeld module defined in Example
		\ref{prop:non-Carlitz-positive-intro}. Then the following hold.
		\begin{enumerate}
			\item $\End_{\overline K}(\phi)=A$.
			
			\item The associated rank-$1$ Drinfeld module given by
			$\phi^1_T=T+T^{10}\tau$ is such that
$     \phi^1$ is not sign-normalized,   $	\phi^1\not\simeq_K\Carl$, 
 and 
			$	K(\phi^1[T])=K(\Carl[T])$.

			\item For every nonzero monic $\fm\in A$, $\Gal(K(\phi[\fm])/K)
				\simeq
				\GL_2(A/\fm A)$.
			
			\item For every monic $\fm\in A$,
			$r_{\fm,T}=1$.

		\end{enumerate}
	\end{lemma}
	
	\begin{proof}
		Recall that $\phi=u^{-1}\phi_0u$, where $u^4=-T$. Since $\phi$ and $\phi_0$ are isomorphic over
		$\overline K$, we have an isomorphism of $A$-algebras
		\begin{align}\label{same-endo}
			\End_{\overline K}(\phi)
			\cong_A\End_{\overline K}(\phi_0),
		\end{align}
		Now, it follows from
		\cite[Theorem 7.5.17(5)]{Papikian}, that since the $j$-invariant (see \cite[Definition 3.8.3]{Papikian})
		$j(\phi_0)
			=
			-1^6/T^4
			=
			-T^{-4}$,
		is not integral over $A$, we have $\End_{\overline K}(\phi_0)=A$. Thus, by \eqref{same-endo},
	$
			\End_{\overline K}(\phi)=A.
$
		This proves part (1).
		
		We next prove part (2). By \eqref{def0:phi1},  the associated rank-$1$
		Drinfeld module is given by
		$\phi^1_T=T+T^{10}\tau$, 
		and is not sign-normalized since $T^{10} \notin \F_5^*$. 
		If $\phi^1\simeq_K\Carl$, then there exists some $v \in K^*$ such that $\phi^1_T = v^{-1}\Carl_Tv=T+v^4\tau$,
		and hence $T^{10}$ must be a fourth power in
		$K^*$, which is not possible. Hence $\phi^1\not\simeq_K\Carl$.
		
		Now, put $\eta\coloneqq T^{-2}u^{-1}$.
		Since $u^4=-T$, we have
		$\eta^4=-T^{-9}$.
		Therefore, $\phi^1_T(\eta)
			=
			T\eta+T^{10}\eta^5		=
			\eta\left(T+T^{10}\eta^4\right)
			=
			0$.
		Thus, $\eta$ is a nonzero $T$-torsion point of $\phi^1$. Since
		$\phi^1[T]\simeq A/(T)$ is a free $A/(T)$-module of rank $1$,
		$K(\phi^1[T])
			=
			K(\eta)$.
		On the other hand,
		$K(\eta)=K(u)$,
		and $u$ is a nonzero $T$-torsion point of $\Carl$. Hence
		$K(\phi^1[T])
			=
			K(u)
			=
			K(\Carl[T])$.
		This proves part (2).
		
		We now prove part (3). Since $\phi=u^{-1}\phi_0u$ and
		$u\in K(\Carl[T])$, we have
		\begin{align}\label{example-compositum}
			K(\phi[\fm])K(\Carl[T])
			=
			K(\phi_0[\fm])K(\Carl[T])
		\end{align}
		for every nonzero monic $\fm\in A$. Suppose first that $T\mid\fm$. By part (2),
$
			K(\Carl[T])
			=
			K(\phi^1[T])
			\subseteq
			K(\phi[\fm]),
$
		which implies that the left hand side of \eqref{example-compositum} is equal to $K(\phi[\fm])$. On the other hand, \eqref{phi01-Carl-fields} gives
$
			K(\Carl[T])
			=
			K(\phi_0^1[T])
			\subseteq
			K(\phi_0[\fm]),
	$
		which implies that the right hand side of \eqref{example-compositum} is equal to $K(\phi_0[\fm])$.
		Hence, we  obtain
$
			K(\phi[\fm])=K(\phi_0[\fm]),
$
		and by \eqref{phi0-full-image-example},
$
			\Gal(K(\phi_0[\fm])/K)
			\simeq
			\GL_2(A/\fm A).
$
		Now suppose that $\gcd(\fm,T)=1$. Applying \eqref{phi0-full-image-example} and using CRT gives
		\begin{align}\label{mcapT}
			K(\phi_0[\fm])\cap K(\phi_0[T])=K.
		\end{align}
		Since
$
			K(\Carl[T])
			=
			K(\phi_0^1[T])
			\subseteq
			K(\phi_0[T]),
$
		it follows from \eqref{mcapT} that
$
			K(\phi_0[\fm])\cap K(\Carl[T])=K.
$
		Therefore, by \eqref{example-compositum},
		\begin{align}
			|\GL_2(A/\fm A)|
			&=
			[K(\phi_0[\fm]):K]
			=
			[K(\phi[\fm])K(\Carl[T]):K(\Carl[T])]\\
			&=
			[K(\phi[\fm]):
			K(\phi[\fm])\cap K(\Carl[T])]\\
			&\leq
			[K(\phi[\fm]):K] \leq
			|\GL_2(A/\fm A)|.
		\end{align}
		Hence
	$
			[K(\phi[\fm]):K]
			=
			|\GL_2(A/\fm A)|,
$
		which proves part (3).
		
		It remains to prove part (4). The second assertion of
		\cite[Theorem 1.2]{Zywina2011} implies that, for every nonzero monic
		$\fm\in A$,
		\begin{align}\label{phi0-geometric-example}
			K(\phi_0[\fm])\cap\overline{\F_5}
			=
			\F_5.
		\end{align}
		Using \eqref{example-compositum}, the inclusion
$
			K(\Carl[T])
			=
			K(\phi_0^1[T])
			\subseteq
			K(\phi_0[T]),
$
		and Lemma \ref{lcm}, we obtain
		\begin{align}
			K(\phi[\fm])K(\Carl[T])
			&=
			K(\phi_0[\fm])K(\Carl[T])
			\subseteq
			K(\phi_0[\lcm(\fm,T)]).
		\end{align}
		Thus, by \eqref{phi0-geometric-example},
	$
			K(\phi[\fm])K(\Carl[T])
			\cap
			\overline{\F_5}
			=
			\F_5.
$
		Consequently, $r_{\fm,T}=1$
		for every monic $\fm\in A$. This proves part (4).
	\end{proof}
	
	Now, by part (3) of Lemma \ref{lem:non-Carlitz-example} and the Chinese Remainder Theorem, we may take
	$M(\phi)=1$ (see Lemma \ref{lem:M-phi}).
	Hence,
	\begin{align}\label{example-N-B}
		N(\phi)
		=
		\prod_{\deg\fp=1}\fp,
		\quad
		\mathcal B
		=
		\{\fp \text{ primes}:\deg\fp=1\}, \quad
		\mathcal G
		=
		\{\fp \text{ primes}:\deg\fp\neq1\}.
	\end{align}
	Furthermore, part (2) of Lemma \ref{lem:non-Carlitz-example} gives
	\begin{align}\label{example-Ct-inside-DB}
		K(\Carl[T])
		=
		K(\phi^1[T])
		\subseteq
		K(\phi[T])
		\subseteq
		D_{\mathcal B}.
	\end{align}
	Also, by Lemma \ref{lem:N-phi}, for every square-free monic $\fm$ supported on
	$\mathcal G$,
$
		K(\phi[\fm])\cap D_{\mathcal B}=K.
$
	Together with \eqref{example-Ct-inside-DB}, this gives
$
		K(\phi[\fm])
		\cap
		D_{\mathcal B}K(\Carl[T])
		=
		K.
$
	Thus, we may take $\mathcal S=\varnothing$
	as an entanglement set for $(\phi,T)$, and
$\label{example-L=DB}
		L_{\phi,T,\varnothing}
		=
		D_{\mathcal B}.
$
	Therefore, by Theorem \ref{thm:non-CM-density}, in order to prove
	Example \ref{prop:non-Carlitz-positive-intro}, it remains only to
	show that
$
		\Omega_{\phi,T,\fa,\varnothing}\neq\varnothing
$
	for every $\fa\in(A/(T))^*$.
	
	\begin{proposition}\label{prop:example-Omega-nonempty}
		For every $\fa\in(A/(T))^*$, one has
	$
			\Omega_{\phi,T,\fa,\varnothing}\neq\varnothing.
$
	\end{proposition}
	
	\begin{proof}
		Fix $\fa\in(A/(T))^*$. Put
	$
			\fb\coloneqq\prod_{\fp\in\mathcal B}\fp.
$
		Then, by Lemma \ref{lcm},
		$D_{\mathcal B}=K(\phi[\fb])$. Recall that $L_{\phi,T,\varnothing}=D_{\mathcal B}$. Hence, part (3) of Lemma \ref{lem:non-Carlitz-example}, together with the
		Chinese Remainder Theorem, gives
		\begin{align}\label{example-DB-product}
			\Gal(L_{\phi,T,\varnothing}/K)
			\simeq
			\prod_{\fp\in\mathcal B}
			\Gal(K(\phi[\fp])/K).
		\end{align}
		For every $\fp \in \mathcal B$ with $\fp \neq T$, choose any $\mathbbm 1 \neq \sigma_\fp \in \Gal(K(\phi[\fp])/K)$. Now, since 
		$K(\Carl[T])
			=
			K(\phi^1[T])
			\subseteq
			K(\phi[T])$,
		the restriction map
		$\Gal(K(\phi[T])/K)
			\longrightarrow
			\Gal(K(\Carl[T])/K)$
		is surjective. Hence, we may choose $\sigma_T\in\Gal(K(\phi[T])/K)$
		such that $\sigma_T(\lambda_T)=\Carl_\fa(\lambda_T)$.
		We may further require $\sigma_T\neq\mathbbm 1$. Indeed, this is automatic if
		$\Carl_\fa(\lambda_T)\neq\lambda_T$, while if $\Carl_\fa(\lambda_T)=\lambda_T$, we choose any element
$
			\mathbbm 1 \neq \sigma_T \in \Gal(K(\phi[T])/K(\Carl[T])),
$
		which is nontrivial since
		\begin{align}
			[K(\phi[T]):K]
			=
			|\GL_2(\F_5)|
			>
			4
			=
			[K(\Carl[T]):K].
		\end{align}
		Thus, by construction, $\sigma_T(\lambda_T)=\Carl_\fa(\lambda_T)$
		and
		$\sigma_\fp\neq\mathbbm 1$
		for every $\fp\in\mathcal B$.
		By \eqref{example-DB-product}, we then obtain $\sigma \in \Gal(L_{\phi,T,\varnothing}/K)$ such that $\sigma|_{K(\phi[\fp])}=\sigma_\fp$ for every $\fp \in \mathcal B$.	
		Since
		\begin{align}
			\Omega_{\phi,T,\fa,\varnothing}
			=
			\{
			\sigma\in
			\Gal(K(\phi[\fb])/K):
			\sigma|_{K(\Carl[T])}(\lambda_T)=\Carl_\fa(\lambda_T),
			\sigma|_{K(\phi[\fp])}\neq\mathbbm{1}
			\text{ for every }\fp\in\mathcal B
			\},
		\end{align}
		it then follows  that 
		$\sigma\in\Omega_{\phi,T,\fa,\varnothing}$.
		Hence
		$\Omega_{\phi,T,\fa,\varnothing}\neq\varnothing$.
	\end{proof}
	
	Thus, by Proposition \ref{prop:example-Omega-nonempty} and Theorem
	\ref{thm:non-CM-density},
	$\delta_{\phi,T,\fa}>0$
	for every $\fa\in(A/(T))^*$. The remaining assertions of Example
	\ref{prop:non-Carlitz-positive-intro} follow from Lemma
	\ref{lem:non-Carlitz-example}. This completes the proof. 
	
	\section*{Acknowledgment}
Arnab Saha	is supported by the ARG-MATRICS Grant (Grant No. ANRF/ARGM/2025/001294/ MTR). 

	\section*{Data availability statement}
Data sharing is not applicable to this article as no datasets were generated or analyzed during the current study.
	
	\section*{Conflict of interest  statement}
	On behalf of all authors, the corresponding author states that there is no conflict of interest.
	
	\bibliographystyle{amsplain}
	\bibliography{cyclicity}

\end{document}